\documentclass[11pt]{amsart} 
\usepackage{amssymb}

\newcommand{\Cold}{k\DGmod}

\newcommand{\limto}{{\displaystyle\lim_{\longrightarrow}}}
\newcommand{\rightlim}{\mathop{\limto}}

\newcommand{\leftlim}{\mathop{\displaystyle\lim_{\longleftarrow}}}
\newcommand{\limfromn}{\leftlim\limits_{\raise3pt\hbox{$n$}}}
\newcommand{\limton}{\rightlim\limits_{\raise3pt\hbox{$n$}}}

\newcommand{\rightlimit}[1]{\mathop{\lim\limits_{\longrightarrow}}\limits%
                   _{\raise3pt\hbox{$\scriptstyle #1$}}}

\newcommand{\rightlimitquot}[1]{\mathop{\mathop{\mbox{``lim''}}\limits_{\longrightarrow}}\limits%
                   _{\raise3pt\hbox{$\scriptstyle #1$}}}

\newcommand{\leftlimit}[1]{\mathop{\lim\limits_{\longleftarrow}}\limits%
                   _{\raise3pt\hbox{$\scriptstyle #1$}}}

\newcommand{\ato}{{\mathop{a}\limits_{\to}}}
\newcommand{\afrom}{{\mathop{a}\limits_{\leftarrow}}}

\newcommand{\cAto}{{\mathop{\cA}\limits_{\to}}}
\newcommand{\cAprimeto}{{\mathop{\cA'}\limits_{\to}}}
\newcommand{\cAfrom}{{\mathop{\cA}\limits_{\leftarrow}}}
\newcommand{\cAtofrom}{{\mathop{\cA}\limits_{\leftrightarrow}}}
\newcommand{\tildetocA}{\stackrel{\leadsto}{\cA}}
\newcommand{\tildefromcA}{\stackrel{\triangleleft\!\!\sim}{\cA}}
\newcommand{\tildetofromcA}{\stackrel{\leftrightsquigarrow}{\cA}}
\newcommand{\tildecAto}{{\mathop{\tilde\cA}\limits_{\to}}}

\newcommand{\cBto}{{\mathop{\cB}\limits_{\to}}}
\newcommand{\cBfrom}{{\mathop{\cB}\limits_{\leftarrow}}}

\newcommand{\cCto}{{\mathop{\cC}\limits_{\to}}}
\newcommand{\cCfrom}{{\mathop{\cC}\limits_{\leftarrow}}}

\newcommand{\rightDelta}{{\mathop{\Delta}\limits^{\to}}}
\newcommand{\leftDelta}{{\mathop{\Delta}\limits^{\leftarrow}}}
\newcommand{\leftrightDelta}{{\mathop{\Delta}\limits^{\leftrightarrow}}}

\newcommand{\Fto}{{\mathop{F}\limits_{\to}}}

\newcommand{\cRto}{{\mathop{\cR}\limits_{\to}}}
\newcommand{\cRfrom}{{\mathop{\cR}\limits_{\leftarrow}}}

\newcommand{\Tto}{{\mathop{T}\limits_{\to}}}

\newcommand{\ppito}{{\mathop{\ppi}\limits_{\to}}}
\newcommand{\ppifrom}{{\mathop{\ppi}\limits_{\leftarrow}}}

\newcommand{\neswarrow}{\swarrow\!\!\!\!\!\nearrow}

\newcommand{\cA}{{\mathcal A}}
\newcommand{\cB}{{\mathcal B}}
\newcommand{\cC}{{\mathcal C}}

\newcommand{\cQ}{{\mathcal Q}}
\newcommand{\cR}{{\mathcal R}}
\newcommand{\cK}{{\mathcal K}}

\newcommand{\cD}{{\mathcal D}}

\newcommand{\cM}{{\mathcal M}}

\newcommand{\cS}{{\mathcal S}}
\newcommand{\cT}{{\mathcal T}}

\newcommand{\bI}{{\bf I}}
\newcommand{\bJ}{{\bf J}}  
\newcommand{\BN}{{\mathbb N}}

\newcommand{\BQ}{{\mathbb Q}}
\newcommand{\bR}{{\bf R}}
\newcommand{\BR}{{\mathbb R}}
\newcommand{\BZ}{{\mathbb Z}}

\newcommand{\iso}{\buildrel{\sim}\over{\longrightarrow}}
\newcommand{\mono}{\hookrightarrow}
\newcommand{\qeq}{\buildrel{\approx}\over{\longrightarrow}}
\newcommand{\qeqleft}{\buildrel{\approx}\over{\longleftarrow}}

\newcommand{\ppi}{\xi}

\newcommand{\bcdot}{{\textstyle\cdot}}
\newcommand{\Ltensor}{\buildrel L\over\otimes}

\newcommand{\dg}{\mbox{\bf DG}}
\newcommand{\DG}{{\mathbb {DG}}}

\newcommand{\MC}{\mbox{\rm MC}}   
\newcommand{\Mmod}{\mbox{\rm -mod}}   
\newcommand{\DGalg}{\mbox{\rm DGalg}}   
\newcommand{\DGcoalg}{\mbox{\rm DGcoalg}}   
\newcommand{\DGmod}{\mbox{\rm -DGmod}}   
\DeclareMathOperator{\dgmod}{{-DGmod}}

\DeclareMathOperator{\Alg}{{Alg}}
\DeclareMathOperator{\Barr}{{Bar}}
\DeclareMathOperator{\Card}{{Card}}
\DeclareMathOperator{\CohFunct}{{CohoFunct}}
\DeclareMathOperator{\Cone}{{Cone}}

\DeclareMathOperator{\Funct}{{Funct}}

\DeclareMathOperator{\gr}{{gr}}
\DeclareMathOperator{\Ho}{{Ho}}
\DeclareMathOperator{\Hodot}{{Ho}^{\bcdot}}

\DeclareMathOperator{\naive}{{naive}}
\DeclareMathOperator{\Ob}{{Ob}}
\DeclareMathOperator{\pretr}{{pre-tr}}

\DeclareMathOperator{\stup}{{stup}}
\DeclareMathOperator{\Stand}{{Stand}}
\DeclareMathOperator{\tr}{{tr}}
\DeclareMathOperator{\Tr}{{Tr}}

\DeclareMathOperator{\Coker}{{Coker}}
\DeclareMathOperator{\discr}{{discr}}

\DeclareMathOperator{\End}{{End}}

\DeclareMathOperator{\Ext}{{Ext}}

\DeclareMathOperator{\HHom}{{\bf Hom}}
\DeclareMathOperator{\Hom}{{Hom}}
\DeclareMathOperator{\ind}{{ind}}
\DeclareMathOperator{\Ind}{{Ind}}

\DeclareMathOperator{\id}{{id}}
\DeclareMathOperator{\im}{{Im}}
\DeclareMathOperator{\Mor}{{Mor}}
\DeclareMathOperator{\Mmor}{\cM\it{or}}

\DeclareMathOperator{\Res}{{Res}}

\DeclareMathOperator{\Top}{{top}}

\theoremstyle{definition}

\numberwithin{equation}{section}

\begin{document}

\thanks{Partially supported by NSF grant DMS-0100108}

\title{DG quotients of DG categories}

\author{Vladimir Drinfeld}

\address{Dept. of Math., Univ. of Chicago, 5734 University 
Ave., Chicago, IL 60637, USA}

\email{drinfeld@math.uchicago.edu}

\maketitle

\begin{abstract}
Keller introduced a notion of quotient of a
differential  graded category modulo a full differential
graded subcategory which agrees with Verdier's notion of
quotient of a triangulated category modulo a triangulated
subcategory. This work is an attempt to further develop his
theory.

\end{abstract}

\noindent {\bf Conventions.}
We fix a commutative ring $k$ and write $\otimes$ instead of
$\otimes_k$ and ``DG category" instead of `` differential
graded $k$-category". If $\cA$ is a DG category we write ``DG
module over $\cA$'' instead of ``DG functor from $\cA$ to
the DG category of complexes of $k$-modules'' (more details
on the DG module terminology can be found in \S\ref{DGapp}).
Unless stated otherwise, all categories are assumed to be
small.  Triangulated categories are systematically viewed as
$\BZ$-graded categories (see \ref{gradedstruct}). A
triangulated subcategory $\cC'$ of a triangulated subcategory
$\cC$ is required to be full, but we do not require it to be
strictly full (i.e., to contain all objects of $\cC$
isomorphic to an object of $\cC'$). In the definition of
quotient of a triangulated category we do not require the
subcategory to be thick (see \ref{triang1}-\ref{triang2}).

\section{Introduction}

\subsection{}
It has been clear to the experts since the 1960's
that Verdier's notions of derived category and triangulated
category \cite{V1,V2} are not quite satisfactory: when you
pass to the homotopy category you forget too much. This is
why Grothendieck developed his derivator theory \cite{G,Ma}.

A different approach was suggested by Bondal
and Kapranov \cite{BK}. According to \cite{BK} one should
work with pretriangulated DG categories rather than with
triangulated categories in Verdier's sense (e.g., with the
DG category of bounded above complexes of projective modules
rather than the bounded above derived category of modules).
Hopefully the part of homological algebra most relevant for
algebraic geometry will be rewritten using DG categories or
rather the more flexible notion of $A_{\infty}$-category due
to Fukaya and Kontsevich (see \cite{F,FS,
K1,K2,Ke3,Ke3',KS,L-H,L}), which goes back to Stasheff's
notion of $A_{\infty}$-algebra \cite{St1,St2}.

One of the basic tools developed by Verdier \cite{V1,V2} is
the notion of quotient of a triangulated category by a 
triangulated subcategory. Keller \cite{Ke4''} has started to
develop a theory of quotients in the DG setting. This work
is an attempt to further develop his theory. I tried to make
this article essentially self-contained, in particular it
can be read independently of \cite{Ke4''}.

The notion of quotient in the setting of
$A_{\infty}$-categories is being developed by
Kontsevich -- Soibelman
\cite{KS} and Lyubashenko -- Ovsienko \cite{L-O}). 

\subsection{}
The basic notions related to that of
DG category are recalled in \S2.
Let $\cA$ be a DG category and $\cB\subset\cA$ a full DG
subcategory. Let $\cA^{\tr}$ denote the triangulated category
associated to $\cA$ (we recall its definition in
\ref{pretr}). A {\it DG quotient\,} (or simply a {\it 
quotient\,}) of $\cA$ modulo $\cB$ is a diagram of DG
categories and DG functors 
\begin{equation} \label{DGquot}
 \cA\qeqleft\tilde\cA\stackrel\ppi\to\cC
\end{equation}
such that the DG functor $\tilde\cA\to\cA$ is a 
quasi-equivalence (see \ref{homotnotions} for the
definition), the functor 
$\Ho (\tilde\cA )\to\Ho (\cC )$ is essentially
surjective, and the functor
$\tilde\cA^{\tr}\to\cC^{\tr}$ induces an equivalence
$\cA^{\tr}/\cB^{\tr}\to\cC^{\tr}$. Keller \cite{Ke4''} proved
that a DG quotient always exists (recall that our DG
categories are assumed to be small, otherwise even the
existence of $\cA^{\tr}/\cB^{\tr}$ is not clear). We recall
his construction of the DG quotient in \S\ref{constr2}, and
give a new construction in \S\ref{constr1}.

The new construction is reminiscent of but easier than
Dwyer-Kan localization \cite{DK1, DK2, DK3}. It is very
simple under a certain flatness assumption (which is
satisfied automatically if one works over a field): one just
kills the objects of $\cB$ (see
\ref{constr}). Without this assumption one has to first
replace $\cA$ by a suitable resolution (see \ref{flatres}).

The idea of Keller's original construction of the DG quotient
(see \S\ref{constr2}) is to take the orthogonal complement of
$\cB$ as a DG quotient, but as the orthogonal complement of
$\cB$ in $\cA$ is not necessarily big enough he takes the
complement  not in $\cA$ but in its ind-version $\cAto$
studied by him in \cite{Ke4}.
The reason why it is natural to consider the orthogonal
complement in $\cAto$ is explained in \ref{explanation}.
Of course, instead of $\cAto$ one can use the
pro-version $\cAfrom$.

Keller's construction using $\cAto$  (resp. $\cAfrom$) is
convenient for considering right (resp. left) derived DG
functors (see \S\ref{Derivedsect}).

In \ref{5.1} we show that the DG quotient of $\cA$ modulo
$\cB$ is ``as unique as possible'', so one can speak of {\it
thhe\,} DG quotient of $\cA$ modulo $\cB$
(``thhe'' is the homotopy version of ``the'').
In \ref{univ} and \ref{uniqueness} we give another
explanation of uniqueness. Unfortunately, both explanations
are somewhat clumsy.

\subsection{Hom complexes of the DG quotient}     
\label{homdescription}

We are going to describe them first as objects of the derived
category of $k$-modules (see \ref{homdescription1}), then in
a stronger sense (see \ref{homdescription2}). We will do it
by successive approximation starting with less precise and
less technical statements.

\subsubsection{}      \label{homdescription1}
Each construction of the DG quotient shows that if
$X,Y\in\Ob\cA$, $\tilde X,\tilde Y\in\Ob\tilde\cA$, 
$\tilde X\mapsto X$, $\tilde Y\mapsto Y$ then the complex 
\begin{equation} \label{homderived}
\Hom_{\cC}(\ppi (\tilde X),\ppi (\tilde Y)) 
\end{equation}
viewed as an object of the derived category of complexes of
$k$-modules  is canonically isomorphic to 
\begin{equation} \label{thecone}
\Cone (h_Y\Ltensor_{\cB}\tilde h_X\to\Hom (X,Y)),
\end{equation}
where $h_Y$ is the right DG $\cB$-module defined by
$h_Y(Z):=\Hom (Z,Y)$, $Z\in\cB$, and
$\tilde h_X$ is the left DG $\cB$-module defined by $\tilde
h_X(Z):=\Hom (X,Z)$, $Z\in\cB$. One can compute
$h_Y\Ltensor_{\cB}\tilde h_X$ using a semi-free resolution of
$h_Y$ or $h_X$ (see \ref{semifreedom} for the definition of
``semi-free''), and this corresponds to Keller's
construction of the DG quotient. If $h_Y$ or 
$\tilde h_X$ is homotopically flat over
$k$ (see \ref{hoflat} for the definition of ``homotopically
flat") then one can compute $h_Y\Ltensor_{\cB}\tilde h_X$
using the bar resolution, and this corresponds to the new
construction of the DG quotient (see \ref{raznoe}(i)). 

\subsubsection{}           \label{homdescription2}
Let $D(\cA )$ denote the derived category of {\it right\,} DG
modules over $\cA$. By \ref{derived} the functor
$D(\cA )\to D(\tilde\cA )$ is an equivalence, so 
for fixed $\tilde Y\in\Ob\tilde\cA$ the
complex (\ref{homderived}) defines an object of $D(\cA )$.
This object is canonically isomorphic to (\ref{thecone}).
Quite similarly, for fixed $\tilde X\in\Ob\tilde\cA$ the
complex (\ref{homderived}) viewed as an object of 
$D (\tilde\cA^{\circ})$ is canonically isomorphic to
(\ref{thecone}). If $\tilde\cA$ is homotopically flat over
$k$ (see \ref{hoflat}) then (\ref{homderived})
and (\ref{thecone}) are canonically isomorphic in 
$D (\tilde\cA\otimes_k\tilde\cA^{\circ})$ (see
\ref{raznoe}(i)). (Without the homotopical flatness
assumption they are canonically isomorphic as objects of the
category $D(\cA\Ltensor\cA^{\circ})$ defined in
\ref{triangrem}.)

\subsubsection{}            \label{homdescription3}
Let (\ref{homderived})$_Y$ denote (\ref{homderived}) viewed
as an object of $D(\cA )$. The morphism
$\mbox{(\ref{thecone})}_Y\to\mbox{(\ref{homderived})}_Y$
mentioned in \ref{homdescription1} and \ref{homdescription2}
is uniquely characterized by the following property: the
composition
$h_Y:=\Hom (?,Y)\to
\mbox{(\ref{thecone})}_Y\to\mbox{(\ref{homderived})}_Y$
equals the obvious morphism 
$\Hom (?,Y)\to\mbox{(\ref{homderived})}_Y$. To prove the
existence and uniqueness of such a morphism we may assume
that $\tilde\cA=\cA$ and the DG functor $\tilde\cA\to\cA$
equals $\id_{\cA}$. Rewrite the DG $\cA^{\circ}$-module
$X\mapsto h_Y\Ltensor_{\cB}\tilde h_X$ as
$L\Ind\cdot\Res h_Y$ (here $\Res :D(\cA )\to D(\cB )$ is the
restriction functor and $L\Ind$ is its left adjoint, i.e.,
the derived induction functor) and notice that 
$\Hom (L\Ind\cdot\Res h_Y, M)=0$ for every DG
$\cA^{\circ}$-module $M$ with $\Res M=0$, in particular for
$M=\mbox{(\ref{homderived})}_Y$. As
$\Res\mbox{(\ref{homderived})}_Y=0$ and 
$\Res\cdot L\Ind\simeq\id$, the fact that our morphism
$\mbox{(\ref{thecone})}_Y\to\mbox{(\ref{homderived})}_Y$ is
an isomorphism is equivalent to the implication
(i)$\Rightarrow$(ii) in the following proposition.

\subsection{\bf Proposition}    
\label{resolutionlemmanew} 
{\it Let $\ppi :\cA\to\cC$ be a DG functor and
$\cB\subset\cA$ a full DG subcategory such that the objects
of $\ppi (\cB )$ are contractible and 
$\Ho (\ppi ):\Ho (\cA )\to\Ho (\cC )$ is essentially
surjective. Then the following properties are equivalent:

(i) $\ppi :\cA\to\cC$ is a DG quotient of $\cA$ modulo $\cB$;

(ii) for
every $Y\in\cA$ the DG $\cA^{\circ}$-module 
\begin{equation} \label{thecone2}
X\mapsto\Cone (\Hom_{\cA}(X,Y)\to
\Hom_{\cC}(\ppi (X),\ppi (Y)) )
\end{equation}
is in the essential image of the derived induction
functor $L\Ind :D(\cB )\to D(\cA )$;

($ii^{\circ}$) for
every $X\in\cA$ the DG $\cA$-module
\begin{displaymath} 
Y\mapsto\Cone (\Hom_{\cA}(X,Y)\to
\Hom_{\cC}(\ppi (X),\ppi (Y)) )
\end{displaymath}
is in the essential image of 
$L\Ind :D(\cB^{\circ})\to D(\cA^{\circ})$.}

\medskip

The proof is contained in \ref{proofof04}.

\medskip

\noindent {\bf Remark.} A DG $\cA^{\circ}$-module $M$
belongs to the essential image of the derived induction
functor $L\Ind :D(\cB )\to D(\cA )$ if and only if the
morphism $L\Ind\Res M\to M$ is a quasi-isomorphism.

\subsection{On Keller's construction of the DG quotient}
\label{explanation}
As explained in \ref{corollaryproof}, the next proposition
follows directly from \ref{resolutionlemmanew}. The symbol
$\Hodot$ below denotes the graded homotopy category (see
\ref{homotnotions}).

\subsubsection{\bf Proposition}  \label{corollary}
{\it Let $\ppi :\cA\to\cC$ be a DG quotient of $\cA$ modulo
$\cB$ and let $\ppi^* :D(\cC )\to D(\cA )$ be the
corresponding restriction functor. Then 

{(\rm a)} the composition 
$\Hodot (\cC )\mono D(\cC )\to D(\cA )$ is fully faithful; 

{(\rm b)} an object of $D(\cA )$ belongs to its essential
image if and only if it is isomorphic to 
$\Cone (L\Ind\Res a\to a)$ for some 
$a\in\Hodot (\cA )\subset D(\cA )$, where
$L\Ind$ (resp. $\Res$) is the derived induction (resp.
restriction) functor corresponding to $\cB\mono\cA$.}

\medskip

{\bf Remark.} In fact, the whole functor 
$D(\cC )\to D(\cA )$ is fully faithful (see \ref{univ}(ii)
or \ref{exactnessoftoandfrom}(ii))

\subsubsection{}  
So if $\ppi :\cA\to\cC$ is a DG quotient then $\Hodot (\cC)$
identifies with a full subcategory of $D(\cA )$. But 
$D(\cA )=\Hodot (\cAto )$, where $\cAto$ is the DG category
of semi-free DG $\cA^{\circ}$-modules (see
\ref{semifreedom}). Thus $\Hodot (\cC)$ identifies with the
graded homotopy category of a certain DG subcategory of
$\cAto$. This is the DG quotient $\cA\nearrow\cB$
from \S\ref{constr2}.

\subsection{Universal property of the DG quotient}
\label{univintro}

\subsubsection{\bf 2-category of DG categories} 
There is a reasonable way to organize all
(small) DG categories into a 2-category ${\bf DGcat}$, i.e.,
to associate to each two DG categories $\cA_1, \cA_2$ a {\it
category of quasi-functors\,} $T(\cA_1,\cA_2)$ and to define
weakly associative composition functors 
$T(\cA_1,\cA_2)\times T(\cA_2,\cA_3)\to T(\cA_1,\cA_3)$ so
that for every DG category $\cA$ there is a weak unit object
in $T(\cA ,\cA )$. Besides, each $T(\cA_1,\cA_2)$ is equipped
with a graded $k$-category structure, and if $\cA_2$ is
pretriangulated in the sense of \ref{pretr} then
$T(\cA_1,\cA_2)$ is equipped with a triangulated structure.
We need ${\bf DGcat}$ to formulate the universal property
\ref{univ} of the DG quotient. The definition of 
${\bf DGcat}$ will be recalled  in~\S\ref{2catofDG}. Here are
two key examples. 

\medskip

\noindent {\bf Examples.} (i) Let $\cK$ be a DG
model of the derived category of complexes of $k$-modules
(e.g., $\cK=$ the DG category of semi-free DG $k$-modules).
Then $T(\cA ,\cK )$ is the derived category of DG
$\cA$-modules. (If $\cK$ is not small then $T(\cA ,\cK )$ is
defined to be the direct limit of $T(\cA ,\cK')$ for all
small full DG subcategories $\cK'\subset\cK$).

(ii) If $\cA_0$ is the
DG category with one object whose endomorphism DG algebra
equals $k$ then $T(\cA_0,\cA)$ is the graded homotopy
category $\Hodot (\cA)$. 

\medskip

It is clear from the definition of $T(\cA_1 ,\cA_2)$ (see
\S\ref{2catofDG}) or from Example (ii) above that $\Phi\in
T(\cA_1 ,\cA_2)$ induces a graded functor $\Hodot
(\cA_1)\to\Hodot (\cA_2)$ and thus $\Hodot$ becomes a
(non-strict) 2-functor from ${\bf DGcat}$ to that of
graded categories. It is also clear from \S\ref{2catofDG}
that one has a bigger 2-functor $\cA\mapsto\cA^{\tr}$ from 
${\bf DGcat}$ to the 2-category of triangulated categories
(with triangulated functors as 1-morphisms).

A DG functor $F:\cA_1\to\cA_2$ defines an object 
$\Phi_F\in T(\cA_1 ,\cA_2)$ (see \ref{Konts1}). Thus one gets
a 2-functor ${\bf DGcat}_{\naive}\to {\bf DGcat}$, where
${\bf DGcat}_{\naive}$ is the 2-category with DG
categories as objects, DG functors as 1-morphisms, and 
degree zero morphisms of DG functors as 2-morphisms. If $F$
is a quasi-equivalence then $\Phi_F$ is invertible. So a
diagram $\cA_1\qeqleft\tilde\cA_1\stackrel F\to\cA_2$ still
defines an object of $T(\cA_1 ,\cA_2)$. All isomorphism
classes of objects of $T(\cA_1 ,\cA_2)$ come from such
diagrams (see \ref{Konts2} and \ref{DGresol}). 

\subsubsection{\bf Main Theorem}  \label{univ}
{\it Let $\cB$ be a full DG
subcategory of a DG category $\cA$. For all pairs 
$(\cC ,\ppi )$, where $\cC$ is a DG category and 
$\ppi\in T(\cA ,\cC)$, the following properties are
equivalent:

(i) the functor $\Ho (\cA )\to\Ho (\cC )$
corresponding to $\ppi$ is essentially surjective, and the
functor $\cA^{\tr}\to\cC^{\tr}$ corresponding to $\ppi$
induces an equivalence $\cA^{\tr}/\cB^{\tr}\to\cC^{\tr}$;

(ii) for every DG category $\cK$ the functor 
$T(\cC ,\cK )\to T(\cA ,\cK)$ corresponding to $\ppi$ is
fully faithful and $\Phi\in T(\cA ,\cK)$ belongs to its
essential image if and only if the image of $\Phi$ in 
$T(\cB ,\cK)$ is zero.

A pair $(\cC ,\ppi )$ satisfying (i)-(ii) exists and is
unique in the sense of $\,{\bf DGcat}$.}

\medskip

A weaker version of the universal property was proved by
Keller, who worked not with the 2-category ${\bf DGcat}$ but
with the category whose morphisms are 2-isomorphism classes
of 1-morphisms of ${\bf DGcat}$ (see Theorem 4.6,
Proposition 4.1, and Lemma 4.2 of \cite{Ke4''}) . Theorem
\ref{univ} will be proved in \ref{univproof}  using the
following statement, which easily follows (see
\ref{preservedbytensoringproof}) from Proposition
\ref{resolutionlemmanew}.

\subsubsection{\bf Proposition} \label{preservedbytensoring}
{\it Let $\ppi :\cA\to\cC$ be a quotient of a DG
category $\cA$ modulo a full DG subcategory $\cB$. If a  DG
category $\cK$ is homotopically flat over $k$ then 
$\ppi\otimes\id_{\cK}:\cA\otimes\cK\to\cC\otimes\cK$ is a 
quotient of the DG category $\cA\otimes\cK$ modulo
$\cB\otimes\cK$.}

\subsection{More on uniqueness} \label{uniqueness}
Let $(\cC_1 ,\ppi_1)$ and $(\cC_2 ,\ppi_2)$, 
$\ppi_i\in T(\cA ,\cC_i )$, be DG quotients of $\cA$ modulo
$\cB$. Then one has an object $\Phi\in T(\cC_1,\cC_2 )$
defined up to unique isomorphism. In fact, the graded
category $T(\cC_1,\cC_2 )$ comes from a certain DG
category (three choices of which are mentioned
in \ref{DGmodels}) 
and one would like to lift $\Phi$ to a homotopically
canonical object of this DG category. The following argument
shows that this is possible under reasonable assumptions. If
$\cC_1$ and $\cC_2$ are homotopically flat over $k$ in the
sense of \ref{hoflat}  these assumptions hold for the Keller
model (see
\ref{DGmodels},   
in particular (\ref{quasitensor2})\,). 

Suppose that $T(\cA ,\cC_i)$ (resp. $T(\cC_1 ,\cC_2)$) is
realized as the graded homotopy category of a DG category
$\DG (\cA ,\cC_i)$ (resp. $\DG (\cC_1 ,\cC_2)$) and suppose
that the graded functor
\begin{displaymath} 
T(\cA ,\cC_1)\times T(\cC_1 ,\cC_2)\times 
T(\cA ,\cC_2)^{\circ}\to\mbox{\{Graded }
k\mbox{-modules\}}
\end{displaymath}
defined by $(F_1,G,F_2)\mapsto\bigoplus_n\Ext^n(F_2,GF_1)$
is lifted to a DG functor
\begin{equation} \label{quasitensor}
\Psi :\DG (\cA ,\cC_1)\times \DG (\cC_1 ,\cC_2)\times 
\DG (\cA ,\cC_2)^{\circ}\to k\DGmod ,
\end{equation}
where $k\DGmod$ is the DG category of complexes of
$k$-modules. We claim that once $\ppi_i$, $i\in\{ 1,2\}$, is 
lifted to an object of $\DG (\cA ,\cC_i)$
one can lift $\Phi\in T(\cC_1,\cC_2 )$ to an object of
$\DG (\cC_1,\cC_2 )$ in a homotopically canonical way.
Indeed, once $\ppi_i$ is lifted to an object of 
$\DG (\cA,\cC_i)$ the DG functor (\ref{quasitensor}) yields
a DG functor $\psi :\DG (\cC_1 ,\cC_2)\to k\DGmod$ such that
the corresponding graded functor $T(\cC_1 ,\cC_2)\to$
\{Graded $k$-modules\} is corepresentable (it is
corepresentable by $\Phi$). Such a functor defines a
homotopically canonical object of $\DG (\cC_1 ,\cC_2)$ (see
\ref{stupid}-\ref{fixquasirep}).

\subsection{What do DG categories form?}  \label{experts}
To formulate uniqueness of the DG quotient in a more elegant
and precise way than in \ref{uniqueness} one probably has to
spell out the relevant structure on the class of all DG
categories (which is finer than the structure of 
2-category). I hope that this will be done by the experts.
Kontsevich and Soibelman are working on this subject. They
introduce in \cite{KS,KoSo3} a notion of homotopy
$n$-category so that a homotopy $1$-category is same as an
$A_{\infty}$-category (the notion of homotopy category is
defined in \cite{KoSo3} with respect to some category of
``spaces'', and in this description of the results of
\cite{KoSo3} we assume that ``space''=``complex of
$k$-modules''). They show that homotopy $1$-categories form a
homotopy $2$-category and they hope that homotopy
$n$-categories form a homotopy $(n+1)$-category. They also
show that the notion of  homotopy $n$-category is closely
related to the little $n$-cubes operad. E.g., they prove in
\cite{KoSo1,KoSo3} that endomorphisms of the identity
1-morphism of an object of a homotopy $2$-category form an
algebra over the chain complex of the little squares operad
(Deligne's conjecture). As DG categories are
$A_{\infty}$-categories we will hopefully understand what DG
categories form as soon as Kontsevich and Soibelman publish
their results.

In the available texts they assume that the ground ring $k$
is a field. Possibly  the case of an arbitrary ground ring
$k$ is not much harder for experts, but a non-expert like
myself becomes depressed when he comes to the conclusion
that DG models of the  triangulated category $T(\cA ,\cK )$
are available only if you first replace $\cA$ or $\cK$ by
a resolution which is homotopically flat over $k$ 
(see \ref{DGmodels}).

\subsection{Structure of the article} In \S\ref{DGrecall} we
recall the basic notions related to DG categories. In
\S\S\ref{constr1},\ref{constr2} we give the two
constructions of the quotient DG category. In
\S\ref{Derivedsect} and \S\ref{DerivedIsect} we discuss
the notion of derived DG functor. The approach of
\S\ref{Derivedsect} is based on Keller's construction of
the DG quotient, while the approach of \S\ref{DerivedIsect}
is based on any DG quotient satisfying a certain flatness
condition, e.g., the DG quotient from \S\ref{constr1}. In
\S\ref{commutsect} we give an explanation of the uniqueness
of DG quotient. In
\S\S\ref{proofmain}-\ref{univproofsection} 
we prove the theorems formulated in
\S\S\ref{constr1}-\ref{DerivedIsect}.  

Finally, \S\S\ref{triang}-\ref{2catofDG}  
are appendices; hopefully they make this article essentially
self-contained. 

\subsection{} I am very grateful to R.~Bezrukavnikov who
asked me how to introduce the notion of quotient in the
framework of DG categories and drew my attention to Keller's
article \cite{Ke4}. I thank A.~Neeman,
P.~Deligne, V.~Hi\-nich, M.~Kapranov,  V.~Lyuba\-shenko, Y.~Soibelman,
A.~Tikaradze and especially  A.~Beilinson and M.~Kontsevich
for stimulating discussions. I thank
B.~Keller, S.~Majid, M.~Mandell,
J.~P.~May, J.~Stasheff, and D.~Yetter for useful references.

\section{DG categories: recollections and notation}
\label{DGrecall}

\subsection{}  
We fix a commutative ring $k$ and write $\otimes$ instead of
$\otimes_k$ and ``DG category" instead of `` differential
graded $k$-category". So a DG category is a category $\cA$ in which the
sets $\Hom (X,Y)$, $X,Y\in\Ob\cA$, are provided with the
structure of a $\BZ$-graded $k$-module and a differential 
$d:\Hom (X,Y)\to\Hom (X,Y)$ of degree $1$ so that for every
$X,Y,Z\in\Ob\cA$ the composition map 
$\Hom (X,Y)\times\Hom (Y,Z)\to\Hom (X,Z)$ comes from a
morphism of complexes
$\Hom (X,Y)\otimes\Hom (Y,Z)\to\Hom (X,Z)$. Using the super
commutativity isomorphism $A\otimes B\iso B\otimes A$ in the
category of DG $k$-modules one defines for every 
DG category $\cA$ the dual DG category $\cA^{\circ}$ with
$\Ob\cA^{\circ}=\Ob\cA$,
$\Hom_{\cA^{\circ}}(X,Y)=\Hom_{\cA}(Y,X)$ (details can be
found in \S 1.1 of \cite{Ke4}).

The {\it tensor product} of DG categories $\cA$ and
$\cB$ is defined as follows: 

(i) $\Ob (\cA\otimes\cB):=\Ob\cA\times\Ob\cB$; for
$a\in\Ob\cA$ and
$b\in\Ob\cB$ the corresponding object of $\cA\otimes\cB$ is
denoted by $a\otimes b$;

(ii) $\Hom (a\otimes b, a'\otimes b'):=\Hom (a,a')\otimes
\Hom (b,b')$ and the composition map is defined by
$(f_1\otimes g_1)(f_2\otimes g_2):= (-1)^{pq}f_1f_2\otimes
g_1g_2$, $p:=\deg g_1$,  $q:=\deg f_2$.

\subsection{Remark} 
Probably the notion of DG category was introduced around 1964
(G.~M.~Kelly \cite{Ke2} refers to it as a new notion used in
\cite{Ke1} and in an unpublished work by Eilenberg and
Moore).

\subsection{}  \label{homotnotions}
Given a DG category $\cA$ one defines a graded category
$\Hodot (\cA )$ with $\Ob\Hodot (\cA )=\Ob\cA$ by
replacing each $\Hom$ complex by the direct sum of its
cohomology groups. We call
$\Hodot (\cA )$ the {\it  graded homotopy category\,} of
$\cA$. Restricting ourselves to the $0$-th cohomology of
the $\Hom$ complexes we get the {\it homotopy
category\,} $\Ho (\cA)$.

A DG functor $F$ is said to be a {\it quasi-equivalence\,}
if $\Hodot (F):\Hodot (\cA )\to\Hodot (\cB )$ is fully
faithful and $\Ho (F)$ is essentially surjective. We will
often use the notation $\cA\qeq\cB$ for a quasi-equivalence
from $\cA$ to $\cB$. The following two notions are less
reasonable. $F:\cA\to\cB$ is said to be a {\it
quasi-isomorphism\,} if $\Hodot (F)$  is an isomorphism. We
say that $F:\cA\to\cB$ is a  {\it DG equivalence\,} if it is
fully faithful and for every object $X\in\cB$ there is a
closed isomorphism of degree 0 between $X$ and an object of
$F(\cA )$.

\subsection{} \label{pretr}
To a DG category $\cA$ Bondal and Kapranov
associate a triangulated category $\cA^{\tr}$ (or
$\Tr^+(\cA)$ in the notation of \cite{BK}). It is defined as
the homotopy category of a certain DG category
$\cA^{\pretr}$. The idea of the definition of $\cA^{\pretr}$
is to formally add to $\cA$ all cones, cones of morphisms
between cones, etc. Here is the precise definition from
\cite{BK}. The objects of $\cA^{\pretr}$ are ``one-sided
twisted complexes", i.e., formal expressions
$(\bigoplus_{i=1}^nC_i[r_i],q)$, where
$C_i\in\cA$, $r_i\in\BZ$, $n\ge 0$, $q=(q_{ij})$, 
$q_{ij}\in\Hom (C_j,C_i)[r_i-r_j]$ is homogeneous of degree
$1$, $q_{ij}=0$ for $i\ge j$, $dq+q^2=0$. If
$C,C'\in\Ob\cA^{\pretr}$, $C=(\bigoplus_{j=1}^nC_j[r_j],q)$,
$C'=(\bigoplus_{i=1}^mC'_i[r'_i],q')$ then the $\BZ$-graded
$k$-module $\Hom (C,C')$ is the space of matrices
$f=(f_{ij})$, $f_{ij}\in\Hom (C_j,C'_i)[r'_i-r_j]$, and the
composition map 
$\Hom (C,C')\otimes\Hom (C',C'')\to\Hom (C,C'')$ is matrix
multiplication. The differential 
$d:\Hom (C,C')\to\Hom (C,C')$ is defined by
$df:=d_{\naive}f+q'f-(-1)^lfq$ if $\deg f_{ij}=l$, where
$d_{\naive}f:=(df_{ij})$.

$\cA^{\pretr}$ contains $\cA$ as a full DG subcategory.
If $X,Y\in\cA$ and $f:X\to Y$ is a closed morphism of
degree 0 one defines $\Cone (f)$ to be the object 
$(Y\oplus X[1], q)\in\cA^{\pretr}$, where 
$q_{12}\in\Hom (X,Y)[1]$ equals $f$ and
$q_{11}=q_{21}=q_{22}=0$.

\medskip

\noindent {\bf Remark.} As explained in \cite{BK}, one has a
canonical fully faithful DG functor (the Yoneda embedding)
$\cA^{\pretr}\to\cA^{\circ}\DGmod$, where
$\cA^{\circ}\DGmod$ is the DG category of DG
$\cA^{\circ}$-modules ; a DG $\cA^{\circ}$-module is
DG-isomorphic to an object of $\cA^{\pretr}$ if and only if
it is finitely generated and semi-free in the sense of
\ref{semifreedom}. Quite similarly one can
identify $\cA^{\pretr}$ with the DG category dual to that of
finitely generated semi-free DG $\cA$-modules.

\medskip
 
A non-empty DG category $\cA$ is said to be {\it
pretriangulated\,} if for every $X\in\cA$, $k\in\BZ$ the
object
$X[k]\in\cA^{\pretr}$ is homotopy equivalent to an object of
$\cA$ and for every closed morphism $f$ in $\cA$ of degree 0
the object $\Cone (f)\in\cA^{\pretr}$ is homotopy equivalent
to an object of
$\cA$. We say that $\cA$ is {\it strongly pretriangulated\,}
(+-pretriangulated in the terminology of \cite{BK}) if same
is true with ``homotopy equivalent'' replaced by
``DG-isomorphic'' (a DG-isomorphism is an invertible closed 
morphism of degree 0). 

If $\cA$ is pretriangulated then every closed degree 0
morphism $f:X\to Y$ in $\cA$ gives rise to the usual
triangle $X\to Y\to\Cone (f)\to X[1]$ in 
$\Ho (\cA )$. Triangles of this type and those isomorphic
to them are called distinguished. Thus if $\cA$ is
pretriangulated then $\Hodot (\cA )$ becomes a triangulated
category (in fact, the Yoneda embedding idenitifies
$\Hodot (\cA )$ with a triangulated subcategory of
$\Hodot (\cA^{\circ}\DGmod)$).

If $\cA$ is pretriangulated (resp. strongly pretriangulated)
then every object of $\cA^{\pretr}$ is homotopy equivalent
(resp. DG-isomorphic) to an object of $\cA$. As explained
in \cite{BK}, the DG category $\cA^{\pretr}$ is always
strongly pretriangulated, so
$\cA^{\tr}:=\Hodot (\cA^{\pretr})$ is a triangulated
category.

\subsection{Proposition}   \label{reasonable}
{\it If a DG functor $F:\cA\to\cB$ is a quasi-equivalence
then same is true for the corresponding DG functor
$F^{\pretr}:\cA^{\pretr}\to\cB^{\pretr}$.}

\medskip

The proof is standard.

\subsection{Remark} Skipping the condition ``$q_{ij}=0$ for
$i\ge j$'' in the definition of $\cA^{\pretr}$ one gets
the definition of the DG category $\mbox{Pre-Tr}(\cA )$
considered by Bondal and Kapranov \cite{BK}. In Proposition
\ref{reasonable} one {\it cannot\,} replace $\cA^{\pretr}$
and $\cB^{\pretr}$ by $\mbox{Pre-Tr}(\cA )$ and
$\mbox{Pre-Tr}(\cB )$. E.g., suppose that $\cA$ and $\cB$
are DG algebras (i.e., DG categories with one object), namely
$\cA$ is the de Rham algebra of a $C^{\infty}$
manifold $M$ with trivial real cohomology and nontrivial
$\pi_1$, $\cB=\BR$, and $F:\cA\to\cB$ is the evaluation
morphism corresponding to a point of $M$. Then
$\mbox{Pre-Tr}(F):\mbox{Pre-Tr}(\cA )\to\mbox{Pre-Tr}(\cB )$
is not a quasi-equivalence. To show this notice that
$K_0(M)\otimes\BQ=\BQ$, so there exists a vector bundle $\xi$
on $M$ with an integrable connection $\nabla$ such that
$\xi$ is trivial but $(\xi ,\nabla )$ is not. $\xi$-valued 
differential forms form a DG $\cA$-module $M$ which is free
as a graded $\cA$-module. Considering $M$ as an object of
$\mbox{Pre-Tr}(\cA )$ we see that $\mbox{Pre-Tr}(F)$
is not a quasi-equivalence.

\subsection{Derived category of DG  modules} \label{derived}
Let $\cA$ be a DG category. Following \cite{Ke4} we denote
by $D(\cA )$ the derived category of DG
$\cA^{\circ}$-modules, i.e., the Verdier quotient of the
homotopy category of DG $\cA^{\circ}$-modules by the
triangulated subcategory of acyclic DG $\cA^{\circ}$-modules.
According to Theorem~10.12.5.1 of \cite{BL} (or Example 7.2
of \cite{Ke4}) if a DG functor $\cA\to\cB$ is a
quasi-equivalence then the restriction functor $D(\cB )\to
D(\cA )$ and its left adjoint functor (the derived induction
functor) are equivalences. This also follows from
\ref{qu-eq} because $D(\cA )$ can be identified with the
homotopy category of semi-free DG
$\cA^{\circ}$-modules (see \ref{semifreedom}).

\subsection{} \label{fiberprod}
Given DG functors $\cA'\to\cA\leftarrow\cA''$ one defines
$\cA'\times_{\cA}\cA''$ to be the fiber product in the
category of DG categories. This is the most naive
definition (one takes the fiber product both at the
level of objects and at the level of morphisms). More
reasonable versions are discussed in \S\ref{diagapp}.

\subsection{} \label{Mmor}

To a DG category $\cA$ we associate a new DG category
$\Mmor\cA$, which is equipped with a DG functor 
$\Cone :\Mmor\cA\to\cA^{\pretr}$. The objects of $\Mmor\cA$
are triples $(X,Y,f)$, where $X,Y\in\Ob\cA$ and $f$ is a
closed morphism $X\to Y$ of degree $0$. At the level of
objects $\Cone (X,Y,f)$ is the cone of $f$. We define
$\Hom ((X,Y,f),(X',Y',f'))$ to be the subcomplex 
\begin{displaymath} 
\{ u\in\Hom (\Cone (f),\Cone (f'))\, |\,\pi'ui=0\},
\end{displaymath}
where $i:Y\to\Cone (f)$ and $\pi':\Cone (f')\to X'[1]$ are
the natural morphisms. At the level of morphisms,
$\Cone :\Hom ((X,Y,f),(X',Y',f'))\to\Hom (\Cone (f),\Cone
(f'))$ is defined to be the natural embedding. Composition of
the morphisms of $\Mmor\cA$ is defined so that 
$\Cone :\Mmor\cA\to\cA^{\pretr}$ becomes a DG functor. There
is an obvious DG functor $\Mmor\cA\to\cA\times\cA$ such that
$(X,Y,f)\mapsto (X,Y)$.

\subsection{} \label{Mmorstup}
Given a DG category $\cA$ one has the ``stupid'' DG category
$\Mmor_{\stup}\cA$ equipped with a DG functor 
$F:\Mmor_{\stup}\cA\to\cA\times\cA$: it has the same objects
as $\Mmor\cA$ (see \ref{Mmor}), 
$\Hom ((X,Y,f),(X',Y',f'))$ is the subcomplex 
\begin{displaymath} 
\{ (u,v)\in\Hom (X,X')\times\Hom (Y,Y')\, |\, f'u=vf\},
\end{displaymath}
$F(X,Y,f):=(X,Y)$, $F(u,v)=(u,v)$, and composition of
the morphisms of $\Mmor_{\stup}\cA$ is defined so that 
$F:\Mmor_{\stup}\cA\to\cA\times\cA$ becomes a DG functor.
There are canonical DG functors
$\Phi :\Mmor_{\stup}\cA\to\Mmor\cA$ and
$\Psi :\Mmor\cA\to\Mmor_{\stup}\cA$ such that 
$\Phi (X,Y,f):=(X,Y,f)$, $\Psi (X,Y,f):=(Y,\Cone (f),i)$,
where $i:Y\to\Cone (f)$ is the natural morphism. So one gets
the DG functor
\begin{equation}   \label{1/3}
\Phi\Psi :\Mmor\cA\to\Mmor\cA
\end{equation}

\section{A new construction of the DG quotient}
\label{constr1}

\subsection{Construction} \label{constr}
Let $\cA$ be a DG category and $\cB\subset\cA$ a full DG
subcategory. We denote by $\cA/\cB$ the DG category obtained
from $\cA$ by adding for every object $U\in\cB$ a  morphism
$\varepsilon_U:U\to U$ of degree
$-1$ such that $d(\varepsilon_U)=\id_U$   (we add neither
new objects nor new relations between the morphisms).

\medskip

So for $X,Y\in\cA$ we have an isomorphism of graded
$k$-modules (but not an isomorphism of complexes)
\begin{equation} \label{0}
 \bigoplus_{n=0}^{\infty}\Hom_{\cA/\cB}^n(X,Y)\iso
 \Hom_{\cA/\cB}(X,Y)
\end{equation}
where $\Hom_{\cA/\cB}^n(X,Y)$ is the direct sum of tensor
products 
$\Hom_{\cA}(U_n,U_{n+1})\otimes k[1]\otimes
\Hom_{\cA}(U_{n-1},U_n)\otimes k[1]\ldots\otimes
k[1]\otimes\Hom_{\cA}(U_0,U_1)$, 
$\;U_0:=X$, $U_{n+1}:=Y$, $U_i\in\cB$ for $1\le i\le n$ (in
particular, $\Hom_{\cA/\cB}^0(X,Y)=\Hom_{\cA}(X,Y)$) ; the
morphism  (\ref{0}) maps 
$f_n\otimes\varepsilon\otimes f_{n-1}\ldots\otimes\varepsilon
\otimes f_0$ 
to $f_n\varepsilon_{U_n}f_{n-1}\ldots\varepsilon_{U_1}f_0$,
where $\varepsilon$ is the canonical generator of $k[1]$.
Using the formula $d(\varepsilon_U)=\id_U$ one can easily
find the differential on the l.h.s. of (\ref{0})
corresponding to the one on the r.h.s. The image of
$\bigoplus_{n=0}^N\Hom_{\cA/\cB}^n(X,Y)$ is a subcomplex of 
$\Hom_{\cA/\cB}(X,Y)$, so we get a filtration on
$\Hom_{\cA/\cB}(X,Y)$. The map (\ref{0}) induces an
isomorphism of complexes 
\begin{equation} \label{0gr}
 \bigoplus_{n=0}^{\infty}\Hom_{\cA/\cB}^n(X,Y)\iso
 \gr\Hom_{\cA/\cB}(X,Y)
\end{equation}

\subsection{Example} If $\cA$ has a single object $U$
with $\End_{\cA}U=R$ then $\cA/\cA$ has a single object $U$
with $\End_{\cA/\cA}U=\tilde R$, where the DG algebra
$\tilde R$ is obtained from the DG algebra $R$ by adding a
new generator $\varepsilon$ of degree $-1$ with $d\varepsilon
=1$. As a DG $R$-bimodule, $\tilde R$ equals 
$\Cone (\Barr (R)\to R)$, where $\Barr (R)$ is the
bar resolution of the DG $R$-bimodule $R$. Both descriptions
of $\tilde R$ show that it has zero cohomology.

A more interesting example can be found in \ref{vogt}.

\subsection{} \label{hoflat}
The triangulated functor 
$\cA^{\tr}\to (\cA/\cB)^{\tr}$ maps $\cB^{\tr}$ to zero
and therefore induces a triangulated functor 
$\Phi:\cA^{\tr}/\cB^{\tr}\to (\cA/\cB)^{\tr}$. Here
$\cA^{\tr}/\cB^{\tr}$ denotes Verdier's quotient (see
\S\ref{triang}). We will prove that if $k$ is a field then
$\Phi$ is an equivalence. For a general ring $k$ this is
true under an additional assumption. E.g., it is enough to
assume that $\cA$ is homotopically flat over $k$ (we prefer
to use the name ``homotopically flat'' instead of
Spaltenstein's name ``K-flat'' which is probably due to the
notation $K(\cC )$ for the homotopy category of complexes in
an additive category $\cC$). A DG category $\cA$ is said to
be  {\it homotopically flat} over $k$ if for every
$X,Y\in\cA$ the complex $\Hom(X,Y)$ is homotopically flat
over $k$ in Spaltenstein's sense \cite{Sp}, i.e., for every
acyclic complex $C$ of $k$-modules $C\otimes_k\Hom(X,Y)$ is
acyclic. In fact, homotopical flatness of $\cA$ can be
replaced by one of the following weaker assumptions:
\begin{equation} \label{FL1}
\mbox{$\Hom (X,U)$ is homotopically flat over $k$ for all
$X\in\cA$, $U\in\cB$;}
\end{equation}
\begin{equation} \label{FL2}
\mbox{$\Hom (U,X)$ is homotopically flat over $k$ for all
$X\in\cA$, $U\in\cB$.}
\end{equation}
Here is our first main result.

\subsection{Theorem} \label{main}
{\it Let $\cA$ be a DG category and $\cB\subset\cA$ a full DG
subcategory. If either (\ref{FL1}) or (\ref{FL2}) holds then
$\Phi:\cA^{\tr}/\cB^{\tr}\to (\cA/\cB)^{\tr}$ is an
equivalence.}

\medskip

The proof is contained in \S\ref{proofmain}.

\subsection{} \label{flatres}
If (\ref{FL1}) and (\ref{FL2}) are
not satisfied one can construct a diagram (\ref{DGquot}) by
choosing a homotopically flat resolution $\tilde\cA\qeq\cA$
and putting $\cC:=\tilde\cA/\tilde\cB$, where
$\tilde\cB\subset\tilde\cA$ is the full subcategory of
objects whose image in $\cA$ is homotopy equivalent to an
object of $\cB$. Here ``homotopically flat resolution''
means that $\tilde\cA$ is homotopically flat and the DG
functor $\tilde\cA\to\cA$ is a quasi-equivalence (see
\ref{homotnotions}). The existence of homotopically flat
resolutions  of $\cA$ follows from Lemma \ref{DGresol}.

\subsection{Remarks} \label{raznoe}

(i) If (\ref{FL1}) or (\ref{FL2}) holds then one can compute
(\ref{thecone}) using the bar resolution of the DG 
$\cB$-module $\tilde h_X$ or the DG $\cB^{\circ}$-module
$h_Y$. The corresponding complex representing the object
(\ref{thecone}) of the derived category is precisely
$\Hom_{\cA/\cB}(X,Y)$.

(ii) Let $\tilde\cA$ and $\tilde\cB$ be as in \ref{flatres}
and suppose that (\ref{FL1}) or (\ref{FL2}) holds for both
$\cB\subset\cA$ and $\tilde\cB\subset\tilde\cA$. Then the DG
functor $\tilde\cA/\tilde\cB\to\cA/\cB$ is a
quasi-equivalence, i.e., it induces an equivalence of the
corresponding homotopy categories. This follows from Theorem
\ref{main}. One can also directly show that if
$X,Y\in\Ob(\cA/\cB)=\Ob\cA$ are the images of
$\tilde X,\tilde Y\in\Ob(\tilde\cA/\tilde\cB)=\Ob\tilde\cA$
then the morphism $\Hom_{\tilde\cA/\tilde\cB}(\tilde X,
\tilde Y)\to\Hom_{\cA/\cB}(X,Y)$ is a quasi-isomorphism (use
(\ref{0gr}) and notice that the morphism
$\Hom_{\tilde\cA/\tilde\cB}^n(\tilde X, \tilde Y) \to
\Hom_{\cA/\cB}^n(X,Y)$ is a quasi-isomorphism for every
$n$; this follows directly from the definition of
$\Hom^n$ and the fact that (\ref{FL1}) or (\ref{FL2}) holds
for $\cB\subset\cA$ and $\tilde\cB\subset\tilde\cA$).

(iii) Usually the DG category $\cA/\cB$ is huge. E.g., if
$\cA$ is the DG category of all complexes from some universe
$U$ and $\cB\subset\cA$ is the subcategory of acyclic
complexes then the complexes $\Hom_{\cA/\cB}(X,Y)$,
$X,Y\in\cA$, are not $U$-small for obvious reasons (see
\cite{GV}, \S 1.0 for the terminology) even though
$(\cA/\cB)^{\tr}$ is a $U$-category. But it follows from
Theorem \ref{main} that whenever $(\cA/\cB)^{\tr}$ is a
$U$-category there exists an $A_{\infty}$-category $\cC$
with $U$-small $\Hom$ complexes equipped with an
$A_{\infty}$-functor $\cC\to\cA/\cB$ which is a
quasi-equivalence (so one can work with $\cC$ instead of
$\cA/\cB$).

(iv) The DG category $\cA/\cB$ defined in \ref{constr}
depends on the ground ring $k$, so the full notation should
be $(\cA/\cB)_k$. Given a morphism $k_0\to k$ we have a
canonical functor $F:(\cA/\cB )_{k_0}\to (\cA/\cB )_k$. If
(\ref{FL1}) or  (\ref{FL2}) holds for both $k_0$ and $k$ then
the functor $(\cA/\cB )_{k_0}\to (\cA/\cB )_k$ is
a quasi-isomorphism by \ref{main}.

\subsection{Example}    \label{vogt}
\subsubsection{}       \label{vogt1}
Let $\cA_0$ be the DG category with two objects $X_1,X_2$
freely generated by a morphism $f:X_1\to X_2$ of degree 0
with $df=0$ (so $\Hom (X_i,X_i)=k$, $\Hom (X_1,X_2)$ is the
free module $kf$ and $\Hom (X_2,X_1)=0$). Put
$\cA:=\cA_0^{\pretr}$. Let $\cB\subset\cA$ be the full DG
subcategory with a single object $\Cone (f)$. Instead of
describing the whole DG quotient $\cA/\cB$ we will describe
only the full DG subcategory $(\cA/\cB )_0\subset\cA/\cB$
with objects $X_1$ and $X_2$ (the DG functor 
$(\cA/\cB )_0^{\pretr}\to (\cA/\cB )^{\pretr}$ is a DG
equi\-valence in the sense of \ref{homotnotions}, so
$\cA/\cB$ can be considered as a full DG subcategory of
$(\cA/\cB )_0^{\pretr}$). Directly using the definition of
$\cA/\cB$ (see \ref{constr}) one shows that $(\cA/\cB )_0$
equals the DG category $\cK$ freely generated by our original
$f:X_1\to X_2$ and also a morhism $g:X_2\to X_1$ of degree
0, morphisms $\alpha_i :X_i\to X_i$ of degree $-1$, and a
morphism $u:X_1\to X_2$ of degree $-2$ with the differential
given by $df=dg=0$, $d\alpha_1=gf-1$, $d\alpha_2=fg-1$,
$du=f\alpha_1-\alpha_2f$. On the other hand, one has the
following description of $\Hodot ((\cA/\cB )_0$).

\subsubsection{\bf Lemma}     \label{contractiblequot}
{\it $\Ext^n_{\cA/\cB}(X_i,X_j)=0$ for $n\ne 0$,
$\Ext^0_{\cA/\cB}(X_i,X_i)=k$, and
$\Ext^0_{\cA/\cB}(X_1,X_2)$,
$\Ext^0_{\cA/\cB}(X_2,X_1)$ are free $k$-modules
generated by $f$ and $f^{-1}$.}

\medskip

As $(\cA/\cB )_0=\cK$ one gets the following corollary. 

\subsubsection{\bf Corollary}     \label{contractibility}
{\it $\cK$ is a resolution of the $k$-category
$\bI_2$ generated by the category $\bJ_2$ with 2 objects and
precisely one morphism with any given source and target.}

\medskip

Clearly $\cK$ is semi-free in the sense of \ref{semifree}.

\subsubsection{\bf Proof of Lemma \ref{contractiblequot}}
\label{contractiblequotproof}
By \ref{main}, $\Hodot (\cA/\cB)=\cA_0^{\tr}/\cB^{\tr}$. As
$X_2\in (\cB^{\tr})^{\perp}$ the map
$\Ext^n_{\cA}(X_i,X_2)\to
\Ext^n_{\cA_0^{\tr}/\cB^{\tr}}(X_i,X_2)$, $i=1,2$,
is an isomorphism by \ref{known}. Therefore
$\Ext^n_{\cA/\cB}(X_i,X_2)$ is as stated in the lemma. But
$f:X_1\to X_2$ becomes an isomorphism in $\Ho (\cA/\cB)$, so
$\Ext^n_{\cA/\cB}(X_i,X_1)$ is also as stated.
\hfill\qedsymbol

\medskip

\subsubsection{\bf Modification of the proof}
In the above proof we used Theorem \ref{main} and \ref{known}
to show that
$\varphi :\Ext^n_{\cA}(X_i,X_2)\to\Ext^n_{\cA/\cB}(X_i,X_2)$
is an isomorphism. In fact, this follows directly from
(\ref{0gr}), which is an immediate consequence of the
definition of $\cA/\cB$. Indeed, $\varphi$ is induced by the
canonical morphism
$\alpha:\Hom_{\cA}(X_i,X_2)\to
\Hom_{\cA/\cB}(X_i,X_2)$. 
By (\ref{0gr}) $\alpha$ is injective and
$L:=\Coker\alpha_i$ is the union of an increasing sequence
of subcomplexes 
$0=L_0\subset L_1\subset\ldots$ such that
$L_n/L_{n-1}=\Hom_{\cA^{\pretr}/\cB}^n(X_i,X_2)$
for $n\ge 1$. Finally,
$\Hom_{\cA^{\pretr}/\cB}^n(X_i,X_2)$ is acyclic for all
$n\ge 1$ because the complex $\Hom_{\cA}(U,X_2)$, $U:=\Cone
(f:X_1\to X_2)$, is contractible. 

\subsubsection{\bf Remarks} \label{vogtrem}
(i) The DG category $\cK$ from \ref{vogt1} and the fact that
it is a resolution of $\bI_2$ were known to Kontsevich
\cite{K2}. One can come to the definition of
$\cK$ as follows. The naive guess is that already the DG
category $\cK'$ freely generated by $f,g,\alpha_1, \alpha_2$
as above is a resolution of $\bI_2$, but one discovers a
nontrivial element $\nu\in\Ext^{-1}(X_1,X_2)$  by
representing $fgf-f$ as a coboundary in two different ways
(notice that $f(gf-1)=fgf-f=(fg-1)f$). Killing $\nu$ one gets
the DG category $\cK$, which already turns out to be a
resolution of $\bJ_2$.

(ii) The DG category $\cK$ from \ref{vogt1} has a
topological analog $\cK_{\Top}$. This is a topological
category with two objects $X_1,X_2$ freely generated by
morphisms $f\in\Mor (X_1, X_2)$, 
$g\in\Mor (X_2, X_1)$, continuous maps 
$\alpha_i:[0,1]\to\Mor (X_i, X_i)$, and a continuous map
$u:[0,1]\times [0,1]\to\Mor (X_1,X_2)$ with defining
relations $\alpha_i(0)=\id_{X_i}$, $\alpha_1(1)=gf$,
$\alpha_2(1)=fg$, $u(t,0)=f\alpha_1(t)$,
$u(t,1)=\alpha_2(t)f$, $u(0,\tau )=f$, $u(1,\tau )=fgf$. It
was considered by Vogt \cite{V}, who was inspired by an
article of R.~Lashof. The spaces
$\Mor_{\cK_{\Top}}(X_i,X_j)$ are contractible. This can be
easily deduced from \ref{contractibility} using a cellular
decomposition of $\Mor_{\cK_{\Top}}(X_i,X_j)$ such that the
composition maps
\begin{displaymath} 
\Mor_{\cK_{\Top}}(X_i,X_j)\times\Mor_{\cK_{\Top}}(X_j,X_k)
\to\Mor_{\cK_{\Top}}(X_i,X_k)
\end{displaymath}
are cellular and the DG category one gets by replacing the
topological spaces $\Mor_{\cK_{\Top}}(X_i,X_j)$ by their
cellular chain complexes equals $\cK$.

\section{The DG categories $\cAto$ and $\cAfrom$.
Keller's construction of the DG quotient.}
\label{constr2}

The DG category $\cA/\cB$ from \S\ref{constr1} depends on the
ground ring $k$ (see \ref{raznoe}(iv)). Here we
describe Keller's construction of a quotient DG category,
which does not depend at all on $k$ (if you like, assume
$k=\BZ$). The construction makes use of the DG
category $\cAto$ studied by him in \cite{Ke4}, which may be
considered as a DG version of the category of ind-objects.
There is also a dual construction based on $\cAfrom$ (a DG
version of the category of pro-objects).

\subsection{}       \label{indproversions}
If $\cA$ is a DG category we denote by $\cAto$ the DG
category of semi-free DG $\cA^{\circ}$-modules (see
\ref{semifreedom} for the definition of ``semi-free'').
 The notation $\cAto$ has been chosen
because one can think of objects of $\cAto$ as a certain kind
of direct limits of objects of $\cA^{\pretr}$ (see
\ref{smallversion}). We put 
$\cAfrom:=(\mathop{\cA^{\circ}}\limits_{\to})^{\circ}$. Of 
course, the DG categories $\cAto$ and $\cAfrom$ are not
small. They are strongly pretriangulated in the sense of
\ref{pretr}, and $\Ho (\cAto )=\cAto^{\tr}$ identifies with
the derived category $D(\cA )$ of DG $\cA^{\circ}$-modules
(see \ref{semifreedom}). We have the fully faithful DG
functors $\cAfrom\leftarrow\cA\to\cAto$. Given a DG functor
$\cB\to\cA$ one has the induction DG functors
$\cBto\to\cAto$ and $\cBfrom\to\cAfrom$ (see
\ref{induction}). In particular, if $\cB\subset\cA$ is a full
subcategory then $\cBto$, $\cBfrom$ are identified with full
DG subcategories of  $\cAto$, $\cAfrom$.

\subsection{Remark} \label{smallversion}
Here is a small version of $\cAto$. Fix
an infinite set $I$ and consider the following DG category
$\cA^{\pretr\to}_I$ (which coincides with the DG category
$\cA^{\pretr}$ from \ref{pretr} if $I=\BN$). To define an
object of $\cA^{\pretr\to}_I$ make the following changes in
the definition of an object of
$\cA^{\pretr}$. First, replace $\bigoplus_{i=1}^nC_i[r_i]$ by
$\bigoplus_{i\in I}C_i[r_i]$ and require the cardinality of
$\{ i\in I|C_i\ne 0\}$ to be strictly less then that of $I$.
Second, replace the triangularity condition on $q$ by the
existence of an ordering of $I$ such that $q_{ij}\ne 0$ only
for $i<j$ and $\{ i\in I|i<j\}$ is finite for every $j\in I$
(in other words, for $j\in I$ let $I_{<j}$ denote the set
of $i\in I$ for which there is a finite sequence 
$i_0,\ldots , i_n\in I$ with $n>0$, $i_0=j$, $i_n=i$ such
that $q_{{i_{k+1}{i_k}}}\ne 0$, then for every $j\in I$ 
the set $I_{<j}$ should be finite and should not contain
$j$). Morphisms of $\cAto$ are defined to be matrices
$(f_{ij})$ as in \ref{pretr} such that
$\{ i\in I|f_{ij}\ne 0\}$ is finite for every $j\in I$. The
DG functor $\cA\to\cAto$ extends in the obvious way to a
fully faithful DG functor $\cA^{\pretr\to}_I\to\cAto$.

One also has the DG category $\cA^{\pretr\leftarrow}_I:=
((\cA^{\circ})^{\pretr\to}_I)^{\circ}$ and the fully faithful
DG functor $\cA^{\pretr\to}_I\to\cAfrom$.

\subsection{Remark}  \label{qu-eq}

A quasi-equivalence $F:\cA\qeq\cB$ induces quasi-equivalences
$\cAto\qeq\cBto$, $\cAfrom\qeq\cBfrom$, 
$\cA^{\pretr\to}_I\qeq\cB^{\pretr\to}_I$,
$\cA^{\pretr\leftarrow}_I\qeq\cB^{\pretr\leftarrow}_I$ (the
fact that $\cAto\to\cBto$ is a quasi-equivalence was
mentioned in \ref{derived}). This is a  consequence of the
following lemma.

\subsection{Lemma}   \label{Keller}
{\it A triangulated subcategory of $\Ho(\cAto )$ containing
$\Ho(\cA )$ and closed under (infinite) direct sums coincides
with $\Ho (\cAto )$. A triangulated subcategory of 
$\Ho (\cA^{\pretr\to}_I)$ containing $\Ho(\cA )$ and closed
under direct sums indexed by sets $J$ such that 
$\Card J<\Card I$ coincides with $\Ho (\cA^{\pretr\to}_I)$.}

\medskip

This was proved by Keller (\cite{Ke4}, p.69). Key idea: if
one has a sequence of DG $\cA^{\circ}$-modules
$M_i$ and morphisms $f_i:M_i \to M_{i+1}$ then one has an
exact sequence 
$0\to M {\buildrel{1-f}\over{\longrightarrow}} M\to
\limto M_i\to 0$, 
where $M :=\mathop{\bigoplus}_i M_i$ and 
$f:M\to M$ is induced by the $f_i$'s. 

\subsection{} \label{formu'}

Now let $\cB\subset\cA$ be a full DG subcategory. Let
$\cB^{\perp}$ (resp. $^{\perp}\cB$) denote the full DG
subcategory of $\cAto$ (resp. of $\cAfrom$) that consists of
objects $X$ such that for every $b\in\cB$ the complex $\Hom
(b,X)$ (resp. $\Hom (X,b)$) is acyclic. Recall that $D(\cA
)=\Ho (\cAto )=\cAto^{\tr}$.

\subsection{Proposition}   \label{exactnessoftoandfrom}
{\it Let $\ppi :\cA\to\cC$ be a 
quotient of a DG category $\cA$ modulo $\cB\subset\cA$. Then 

(i) $\ppito :\cAto\to\cCto$  is a quotient of $\cAto$
modulo $\cBto$;

(i') $\ppifrom :\cAfrom\to\cCfrom$  is a quotient of
$\cAfrom$ modulo $\cBfrom$;

(ii) the restriction functor $D(\cC )\to D(\cA )$ is fully
faithful, and its essential image consists precisely of
objects of $D(\cA )$ annihilated by the
restriction functor 
$\rho : D(\cA )\to D(\cB )$; the functor 
$D(\cA )/D(\cC )\to D(\cB )$ induced by $\rho$ is an
equivalence.}

\medskip

See \ref{exactnessoftofromproof} for the proof.

\subsection{Proposition}   \label{admis}
{\it ($i$) The essential image of $\cBto^{\tr}$ in
$\cAto^{\tr}$ is right-admissible in the sense of
\ref{admissible}. 

($ii$) The right orthogonal complement of
$\cBto^{\tr}$ in $\cAto^{\tr}$ equals $(\cB^{\perp})^{\tr}$.

($iii$) The functor
$(\cB^{\perp})^{\tr}\to\cAto^{\tr}/\cBto^{\tr}$ is an
equivalence. 

($iv$) The functor
$\cA^{\tr}/\cB^{\tr}\to\cAto^{\tr}/\cBto^{\tr}$ is fully
faithful.

($i^{\circ}$)-($iv^{\circ}$) Statements ($i$)-($iv$) remain
true if one replaces
$\cAto^{\tr}$ and $\cBto^{\tr}$ by $\cAfrom^{\tr}$ and
$\cBfrom^{\tr}$, ``right'' by  ``left'', and $\cB^{\perp}$ by
$^{\perp}\cB$.}

\medskip

The proof will be given in \ref{proofmain2}.

\subsection{Remark}
Keller \cite{Ke4''} derives Proposition
\ref{exactnessoftoandfrom}(i) from Neeman's theorem on
compactly generated triangulated categories (Theorem 2.1 of
\cite{N3}). Statements (i) and (iv) of Proposition
\ref{admis} are particular cases of Lemmas 1.7 and 2.5 of
Neeman's work \cite{N3}.

\subsection{} \label{formu}
Now let $\cA\nearrow\cB\subset\cB^{\perp}$ be the full DG
subcategory of objects $X\in\cB^{\perp}$ such that for some
$a\in\cA$ and some closed morphism $f:a\to X$ of degree 0
the cone of $f$ is homotopy equivalent to an object of
$\cBto$. Let $\cA\swarrow\cB\subset\null^{\perp}\cB$ be the
full DG subcategory of objects $X\in\null^{\perp}\cB$ such
that for some $a\in\cA$ and some closed morphism $f:X\to a$
of degree $0$ the cone of $f$ is homotopy equivalent to an
object of $\cBfrom$. By Proposition \ref{admis} we have the
fully faithful functor
$\cA^{\tr}/\cB^{\tr}\to\cAto^{\tr}/\cBto^{\tr}=
(\cB^{\perp})^{\tr}=\Ho (\cB^{\perp})$, and its essential
image equals $(\cA\nearrow\cB )^{\tr}$. So we get an
equivalence 
\begin{equation} \label{form}
\cA^{\tr}/\cB^{\tr}\iso (\cA\nearrow\cB )^{\tr}
\end{equation}
and a similar equivalence
$\cA^{\tr}/\cB^{\tr}\iso (\cA\swarrow\cB )^{\tr}$.

\subsection{} \label{cAl}

Let us construct a diagram (\ref{DGquot}) with
$\cC=\cA\nearrow\cB$ such that the corresponding functor
$\cA^{\tr}\to(\cA\nearrow\cB )^{\tr}$ induces (\ref{form})
(so $\cA\nearrow\cB$ will become a DG quotient of $\cA$
modulo $\cB$). The DG category $\tilde\cA =\tildetocA$ is
defined as follows. First consider the DG category
$\Mmor\cAto$ (see \ref{Mmor}). Its objects are triples
$(a,Y,g)$, where $a,Y\in\cAto$ and $g$ is a closed
morphism $a\to Y$ of degree 0. We define
$\tildetocA\subset\Mmor\cAto$ to be the full DG subcategory
of triples $(a,Y,g)$ such that $a\in\cA$, 
$Y\in\cA\nearrow\cB\subset\cAto$,  and 
$\Cone (a\stackrel g\to Y)$ is homotopy equivalent to
an object of $\cBto$. The DG functors
$\cA\leftarrow\tildetocA\to\cA\nearrow\cB$ are defined by
$(a,Y,g)\mapsto a$ and $(a,Y,g)\mapsto Y$. 

\subsection{Remarks} \label{tildeprime}
(i) Let $\tildetocA'\subset\Mmor\cAto$
be the full DG subcategory of triples $(P,a,f)\in\Mmor\cAto$
such that $P\in\cBto$, $a\in\cA$, and 
$\Cone (P\stackrel f\to a)\in\cB^{\perp}$. The DG functor 
(\ref{1/3}) (with $\cA$ replaced by $\cAto$) induces a
quasi-equivalence $\tildetocA'\qeq\tildetocA$, so one can use
$\tildetocA'$ instead of $\tildetocA$.

(ii) It follows from the definition of (\ref{1/3}) that the
image of the DG functor $\tildetocA'\to\tildetocA$ is
contained in
$\tildetocA_{\stup}:=\tildetocA\cap\Mmor_{\stup}\cAto$
(see \ref{Mmorstup} for the definition of
$\Mmor_{\stup}\cAto\subset\Mmor\cAto$).

\subsection{} \label{cAr}
Dualizing the construction from \ref{cAl} one gets the full
DG subcategory $\tildefromcA\subset\Mmor\cAfrom$ which
consists of triples $(\bar Y,a,\bar g)$ such that 
$\bar Y\in\cA\swarrow\cB$, $a\in\cA$, and
$\Cone (\bar Y\stackrel{\bar g}\to  a)$ is homotopy 
equivalent to an object of $\cBfrom$. Dualizing
\ref{tildeprime} one gets a DG category $\tildefromcA'$
equipped with a quasi-equivalence
$\tildefromcA'\qeq\tildefromcA$;
$\tildefromcA'\subset\Mmor\cAfrom$ is the full DG subcategory
of triples $(a,\bar P,\bar f)$ such that $a\in\cA$, 
$\bar P\in\cBfrom$, and $\Cone (\bar f)\in\null^\perp\cB$.
The diagrams
$\cA\leftarrow\tildefromcA\to\cA\swarrow\cB$ and 
$\cA\leftarrow\tildefromcA'\to\cA\swarrow\cB$
are also DG quotients of $\cA$ modulo $\cB$. The image of
the DG functor $\tildefromcA'\to\tildefromcA$ is contained in
$\tildetofromcA_{\stup}:=
\tildefromcA\cap\Mmor_{\stup}\cAfrom$.

\subsection{}         \label{3.9}
One can also include the diagrams constructed in \ref{cAl}
and \ref{cAr} into a canonical commutative diagram of DG
categories and DG functors
\begin{equation} \label{kvadrat}
        \begin{array}{ccccc}
\cA&=&\cA&=&\cA\\
\uparrow&&\uparrow&&\uparrow\\
\tildefromcA&\qeqleft&\tildetofromcA&\qeq
&\tildetocA\\
\downarrow&&\downarrow&&\downarrow\\
\cA\swarrow\cB&\qeqleft&\cA\neswarrow\cB&\qeq
&\cA\nearrow\cB\\
     \end{array}
\end{equation}
in which each column is a DG quotient of $\cA$ modulo $\cB$.
The DG category $\tildetofromcA$ is defined to be the
fiber product $\tildefromcA'\times_{\cA}\tildetocA'$, where
$\tildefromcA'$ and $\tildetocA'$ were defined in
\ref{tildeprime}, \ref{cAr} (recall that ``fiber product'' is
understood in the most naive sense, see \ref{fiberprod}).
To define $\cA\neswarrow\cB$ we use the DG category
$\cAtofrom$ such that
$\Ob\cAtofrom:=\Ob\cAfrom\bigsqcup\Ob\cAto$, $\cAto$ and
$\cAfrom$ are full DG subcategories of $\cAtofrom$, and for
$Y\in\Ob\cAto$, $\bar Y\in\Ob\cAfrom$ one has 
$\Hom (Y,\bar Y):=0$,
$\Hom (\bar Y,Y):=Y\otimes_{\cA}\bar Y$ (recall that $Y$ is a
DG $\cA^{\circ}$-module and $\bar Y$ is a DG
$\cA$-module, so $Y\otimes_{\cA}\bar Y$ is well defined, see
\ref{tensorproduct}). For $a\in\cA$ we denote by $\afrom$
(resp. $\ato$) the image of $a$ in $\cAfrom$ (resp.
$\cAto$); we have the ``identity'' morphism
$e=e_a:\afrom\to\ato$. Now define
$\cA\neswarrow\cB\subset\Mmor\cAtofrom$ to be the full DG
subcategory of triples $(\bar Y,Y,f)\in\Mmor\cAtofrom$ such
that $\bar Y\in\cA\swarrow\cB\subset\cAfrom$,
$Y\in\cA\nearrow\cB\subset\cAto$,
and $f:\bar Y\to Y$ can be represented as a composition
$\bar Y\buildrel{g}\over{\to}\afrom
\buildrel{e}\over{\to}\ato\buildrel{h}\over{\to}Y$, 
$a\in\cA$, so that $\Cone (g)$ is homotopy equivalent to an
object of $\cBfrom$ and $\Cone (h)$ is homotopy equivalent
to an object of $\cBto$ ($g$ and $h$ are closed morphisms of
degree 0). 

The DG functors $\cA\neswarrow\cB\to\cA\nearrow\cB$ and
$\cA\neswarrow\cB\to\cA\swarrow\cB$ send 
$(\bar Y,Y,f)\in\cA\neswarrow\cB$ respectively to $Y$ and
$\bar Y$. The DG functor
$\tildetofromcA\to\cA\neswarrow\cB\subset\Mmor\cAtofrom$ is
defined to be the composition
\begin{displaymath} 
\tildetofromcA:=\tildefromcA'\times_{\cA}\tildetocA'\to
\tildefromcA_{\stup}\times_{\cA}\tildetocA_{\stup}
\buildrel{F}\over{\to}\Mmor\cAtofrom
\end{displaymath}
where the DG functors $\tildefromcA'\to\tildefromcA_{\stup}$
and $\tildetocA'\to\tildetocA_{\stup}$ were defined in
\ref{tildeprime}-\ref{cAr} and
$F:\tildefromcA_{\stup}\times_{\cA}\tildetocA_{\stup}\to
\Mmor\cAtofrom$ is the composition DG functor: at the level
of objects, if $u=(a,Y,g:a\to Y)\in\Mmor_{\stup}\cAto$ and 
$\bar u=(\bar Y,a,\bar g:\bar Y\to a)\in
\Mmor_{\stup}\cAfrom$, $a\in\cA$, 
then $F(\bar u,u)=(\bar Y,Y,g\bar g)$; there is no problem
to define the DG functor $F$ at the level of
morphisms because we are working with the
``stupid'' versions $\tildefromcA_{\stup}$,
$\tildetocA_{\stup}$, $\Mmor_{\stup}$ (the ``non-stupid''
composition
$\tildefromcA\times_{\cA}\tildetocA\to\Mmor\cAtofrom$ is
defined as an $A_{\infty}$-functor rather than as a DG
functor).

\section{Derived DG functors}   \label{Derivedsect}

We will define a notion of right derived functor in the DG
setting modeled on Deligne's definition in the triangulated
setting. One can easily pass from right derived DG functors
to left ones by considering the dual DG categories.

\subsection{Deligne's definition} \label{remind} 

Let $G:\cT\to\cT'$ be a triangulated functor between
triangulated categories and $\cS\subset\cT$ a 
triangulated subcategory. Denote by
$\CohFunct (\cT')$ the category of $k$-linear cohomological
functors from $(\cT')^{\circ}$ to the category of
$k$-modules.
$RG$ is defined to be the functor 
$\cT/\cS\to\CohFunct (\cT')$ defined by
\begin{equation}   \label{shorthand}
RG(Y):=\rightlimitquot{(Y\to Z)\in Q_Y} G(Z),
\end{equation}
which is a shorthand for
\begin{equation}  \label{tr-derived}
RG(Y)(X):=\rightlimit{(Y\to Z)\in Q_Y}\Hom (X,G(Z)),
\quad Y\in\cT,\, X\in\cT'\,.
\end{equation}
Here $Q_Y$ is the filtering category of
$\cT$-morphisms $f:Y\to Z$ such that $\Cone (f)$ is
isomorphic to an object of $\cS$.

$RG$ has the following universal property. Let 
$\pi :\cT\to\cT/\cS$ denote the canonical functor and $\nu
:\cT'\to\CohFunct (\cT' )$ the Yoneda embedding. Let 
$\Phi :\cT/\cS\to\CohFunct (\cT' )$ be a graded functor (see
\ref{gradedstruct} for a discussion of the meaning of
``graded''). Then there is a canonical isomorphism
\begin{equation}  \label{universality}
     \Hom (RG,\Phi )=\Hom (\nu G,\Phi\pi )
\end{equation}
functorial in $\Phi$ (here $\Hom$ is the set of morphisms of
graded functors). In particular, if 
$RG(\cT/\cS )\subset\cT'$ then $RG:\cT/\cS\to\cT'$ is a
derived functor in Verdier's sense \cite{V1,V2}.

Let $(\cT/\cS )_G$ be the category of triples
$(Y,X,\varphi )$, where $Y\in\cT/\cS$, $X\in\cT'$, 
$\varphi :X\iso RG(Y)$. The functor 
\begin{equation}  \label{partialderived}
(\cT/\cS )_G\to\cT', \quad (Y,X,\varphi )\mapsto X
\end{equation}
is also denoted by $RG$. We have an equivalence 
$(Y,X,\varphi )\mapsto Y$ between $(\cT/\cS )_G$ and a full
subcategory of $\cT/\cS$ (the full subcategory of objects
$Y\in\cT/\cS$ such that $RG(Y)$ is defined as an object of
$\cT'$).

\medskip

\noindent {\bf Remark.} Deligne (cf. Definition 1.2.1 of
\cite{De}) considers $RG$ as a functor from $\cT/\cS$ to the
category of ind-objects $\ind (\cT')$ rather than to the
category $\CohFunct (\cT')$. In fact, this does not matter.
First of all, the image of the functor $RG$ defined by
(\ref{tr-derived}) is contained in the full subcategory of
ind-representable functors $(\cT')^{\circ}\to k\Mmod$, which
is canonically identified with $\ind (\cT')$ (see \S 8.2 of
\cite{GV}). This is enough for our purposes, but in fact
since $\cT'$ is small {\it every\,} $H\in\CohFunct (\cT')$
is ind-representable by a well known lemma (see, e.g., Lemma
7.2.4 of \cite{N2}), which is a version of Brown's theorem
\cite{Br,Br2}. Proof: by Theorem 8.3.3 of \cite{GV} it
suffices to check that the category $\cT'/H:=\{(X,u)|
X\in\cT',u\in H(X)\}$ is filtering.

\subsection{} \label{RFnew}
Let $\cA$ be a DG category and $\cB\subset\cA$ a full DG
subcategory. Let $F$ be a DG
functor from $\cA$ to a DG category~$\cA'$. To define the
right derived DG functor $RF$ we use the DG quotient
$\cA\nearrow\cB$ from \ref{formu}. By definition, 
$RF:\cA\nearrow\cB\to\cAto'$ is the restriction of the DG
functor $\Fto :\cAto\to\cAto'$ to the DG subcategory
$\cA\nearrow\cB\subset\cB^{\perp}\subset\cAto$. A
2-categorical reformulation of this definition is given in
Remark (ii) from \ref{Induality}.

\subsection{} \label{comparenew}
Let us show that the definition of $RF$ from \ref{RFnew}
agrees with Deligne's definition of the right derived
functor of a triangulated functor between triangulated
categories (see \ref{remind}). 

Suppose we are in the situation of \ref{RFnew}. We have the
DG functor $RF:\cA\nearrow\cB\to\cAto'$  and the
corresponding triangulated functor
$(RF)^{\tr}:(\cA\nearrow\cB)^{\tr}\to (\cAto')^{\tr}$.
Using (\ref{form}) we can rewrite it as 
$(RF)^{\tr}:\cA^{\tr}/\cB^{\tr}\to (\cAto')^{\tr}$. On the
other hand, we have the triangulated  functor
$F^{\tr}:\cA^{\tr}\to (\cA')^{\tr}$ and its derived functor 
$RF^{\tr}:\cA^{\tr}/\cB^{\tr}\to\CohFunct ((\cA')^{\tr})$
(see \ref{remind}). Finally, one has the functor 
$H^{0}:(\cAto')^{\tr}\to\CohFunct ((\cA')^{\tr})$ defined
as follows: a right DG $\cA'$-module $M\in\cAto'$ uniquely
extends to a right DG $(\cA')^{\pretr}$-module $\tilde M$
(cf. \ref{Moritaeq}),   
and $H^{0}(M)$ is defined to be the zeroth cohomology of
$\tilde M$ (or equivalently $H^{0}(M)$ is the cohomological
functor $N\mapsto\Hom_{(\cAto')^{\tr}}(N,M)$, 
$N\in (\cA')^{\tr}\subset (\cAto')^{\tr}$).

Finally, using that
$(\cA')^{\circ}\DGmod=((\cA')^{\pretr})^{\circ}\DGmod$
(see \ref{Moritaeq}) one gets the functor 
$H^{0}:((\cA')^{\circ}\DGmod)^{\tr}
\to\CohFunct ((\cA')^{\tr})$. 

We are going to construct an isomorphism
$RF^{\tr}\iso H^{0}(RF)^{\tr}$. To this end, consider the
diagram
\begin{equation} \label{kvadr}
        \begin{array}{ccccc}
\tildetocA&\qeq&\cA&\buildrel F\over\to&\cA'\\
\downarrow&&\downarrow&&\downarrow\\
\cA\nearrow\cB&\mono&\cAto&\to&\cAto'
     \end{array}
\end{equation}
(see \ref{cAl} for the definition of $\tildetocA$). Its left
square is not commutative, but there is a canonical morphism
from the composition $\tildetocA\to\cA\mono\cAto$ to the
composition $\tildetocA\to\cA\nearrow\cB\mono\cAto$. So we
get a canonical morphism $\varphi$ from the composition
$\tildetocA^{\tr}\to (\cA')^{\tr}\to
(\cAto')^{\tr}\to\CohFunct (\cA')$ to
the composition 
$\tildetocA^{\tr}\to (\cA\nearrow\cB)^{\tr}\to
(\cAto')^{\tr}\to\CohFunct (\cA')$. 
By \ref{cAl} we can identify $\tildetocA^{\tr}$ with
$\cA^{\tr}$ and $(\cA\nearrow\cB)^{\tr}$ with
$\cA^{\tr}/\cB^{\tr}$, so $\varphi$ induces a
morphism                                  
\begin{equation} \label{themorphismnew}
RF^{\tr}\to H^{0}(RF)^{\tr}
\end{equation}
by the universal property (\ref{universality}) of $RF^{\tr}$.

\subsection{Proposition.}  \label{comparisonnew}
{\it The morphism (\ref{themorphismnew}) is an isomorphism.}

\medskip 

See \S\ref{comparisonproof} for a proof.

\subsection{}  \label{RF2new}

Define the DG category $(\cA\nearrow\cB )_F$ to be the
(naive) fiber product of $\cA'\times (\cA\nearrow\cB )$ and
$\rightDelta_{\cAprimeto}$ over $\cAprimeto\times\cAprimeto$,
where $\rightDelta_{\cAprimeto}$ is the ``diagonal'' DG
category defined in \ref{diag} and $\cA\nearrow\cB$ is mapped
to $\cAprimeto$ by $RF$. So the objects of 
$(\cA\nearrow\cB )_F$ are triples $(Y,X,\varphi )$, where
$Y\in\cA\nearrow\cB$, $X\in\cA'$, and $\varphi :X\to RF(Y)$
is a homotopy equivalence.  The DG functor 
$(\cA\nearrow\cB )_F\to\cA'$ defined by 
$(Y,X,\varphi )\mapsto X$ is also called the right derived DG
functor of $F$ and denoted by $RF$. 

Now consider the triangulated functor
$G=F^{\tr}:\cA^{\tr}\to(\cA')^{\tr}$.
It follows from \ref{comparisonnew} that 
$((\cA\nearrow\cB )_F)^{\tr}$ identifies with the
triangulated category
$(\cA^{\tr}/\cB^{\tr})_G$ from \ref{remind} and 
$(RF)^{\tr}:((\cA\nearrow\cB )_F)^{\tr}\to (\cA')^{\tr}$
identifies with Deligne's derived functor 
$RG:(\cA^{\tr}/\cB^{\tr})_G\to (\cA')^{\tr}$.

\subsection{}  \label{RF2new'}
The definition of $(\cA\nearrow\cB )_F$ used
$\rightDelta_{\cAprimeto}$. There are also versions of
$(\cA\nearrow\cB )_F$ using the DG
categories $\leftDelta_{\cAprimeto}$ and
$\leftrightDelta_{\cAprimeto}$ from \ref{diag}. They will be
denoted respectively by $(\cA\nearrow\cB )_{\leftarrow F}$
and $(\cA\nearrow\cB )_{\leftrightarrow F}$. E.g., the
objects of $(\cA\nearrow\cB )_{\leftarrow F}$ are triples
$(Y,X,\psi )$, where $Y\in\cA\nearrow\cB$, $X\in\cA'$, and
$\psi :RF(Y)\to X$ is a homotopy equivalence. We have the
right  derived DG functors 
$RF:(\cA\nearrow\cB )_{\leftarrow F}\to\cA'$ and
$RF:(\cA\nearrow\cB )_{\leftrightarrow F}\to\cA'$.
Sometimes we will
write $(\cA\nearrow\cB )_{\to F}$ instead of
$(\cA\nearrow\cB )_F$. The DG functors 
$(\cA\nearrow\cB )_{\to F}\leftarrow
(\cA\nearrow\cB )_{\leftrightarrow F}
\to(\cA\nearrow\cB )_{\leftarrow F}$ are quasi-equivalences
by \ref{diagolemma}, and one has a canonical commutative
diagram
\begin{equation} \label{kvadratt-1}
        \begin{array}{ccccc}
(\cA\nearrow\cB )_{\to F}&\qeqleft &
(\cA\nearrow\cB )_{\leftrightarrow F} &
\qeq & (\cA\nearrow\cB )_{\leftarrow F} \\
\!\!\!\!\!\!\!\! RF\downarrow&&
\!\!\!\!\!\!\!\! RF\downarrow&&\,\,\,\,\,\,\,\,\,
\downarrow RF
\\
\cA&=&\cA&=&\cA \\
        \end{array}
\end{equation}

\section{Some commutative diagrams} \label{commutsect}

\subsection{Uniqueness of DG quotient} \label{5.1}

Let $\cA$ be a DG category and $\cB\subset\cA$ a full DG
subcategory. Given a quotient (\ref{DGquot}) of $\cA$
modulo $\cB$ we will ``identify'' it with the quotient
$\cA\qeqleft\tildetocA\buildrel{F}\over{\to}\cA\nearrow\cB$
from \ref{cAl}. More precisely, we will construct a canonical
commutative diagram of DG categories
\begin{equation} \label{kvadratt0}
        \begin{array}{ccccc}
\tildetocA&\qeqleft&\underline{\cA}&\qeq&\tilde\cA \\
\!\!\!\! F\downarrow&&\downarrow&&\,\,\,\,\downarrow\ppi \\
\cA\nearrow\cB&\qeqleft&\underline{\cC}&\qeq&\cC \\
        \end{array}
\end{equation}
(the symbols $\qeq$, $\qeqleft$ denote quasi-equivalences).
To this end, notice that the derived DG functor 
$R\ppi :(\tilde\cA\nearrow\tilde\cB )_{\ppi}\to\cC$
defined in \ref{RF2new} and the projection 
$(\tilde\cA\nearrow\tilde\cB)_{\ppi}\to\tilde\cA\nearrow
\tilde\cB$ are quasi-equivalences (here $\tilde\cB$ is the
preimage of $\cB$ in $\tilde\cA$). Put
$\underline{\cC}:=(\tilde\cA\nearrow\tilde\cB )_{\ppi}$.
Define the DG functor $\underline{\cC}\to\cC$ to equal
$R\ppi$ and the DG functor
$\underline{\cC}\to\cA\nearrow\cB$ to be the composition 
$\underline{\cC}=(\tilde\cA\nearrow\tilde\cB)_{\ppi}\to
\tilde\cA\nearrow\tilde\cB\to\cA\nearrow\cB$. We put
$\underline{\cA}:=\stackrel{\leadsto}{\tilde\cA}$, i.e.,
$\underline{\cA}$ is the analog of $\tildetocA$ with 
$(\cA ,\cB )$ replaced by $(\tilde\cA ,\tilde\cB )$. The DG
functor $\underline{\cA}\to\tilde\cA$ is the analog of
$\tildetocA\to\cA$. The DG functor
$\underline{\cA}\to\tildetocA$ is induced by the DG functors
$\tilde\cA\to\cA$ and $\tilde\cB\to\cB$. Finally,
$\underline{\cA}\to\underline{\cC}$ is the DG functor 
${\stackrel{\leadsto}{\tilde\cA}}\to\underline{\cC}$
defined by $(a,Y,g)\mapsto (Y,\ppi (a), \ppito (g))$
(here $a\in\tilde\cA$, 
$Y\in\tilde\cA\nearrow\tilde\cB\subset\tildecAto$,
and $g:a\to Y$ is a closed morphism of degree 0 whose cone is
homotopy equivalent to an object of $\cBto$; recall that an
object of $\underline{\cC}$ is a triple $(Y,X,\varphi )$,
where $Y\in\tilde\cA\nearrow\tilde\cB\subset\tildecAto$,
$X\in\cC$, and $\varphi$  is a homotopy equivalence from $X$
to $R\ppi (Y)$, i.e., the image of $Y$ under 
$\ppito:\tildecAto\to\cCto$).

\subsection{More diagrams (to be used in
\S\ref{DerivedIsect})}

\subsubsection{} \label{resolvedDG}

Now let us consider the case that $\tilde\cA=\cA$ and
the DG functor $\tilde\cA\to\cA$ equals $\id_{\cA}$, so our
quotient (\ref{DGquot}) is just a DG category $\cC$
equipped with a DG functor $\ppi :\cA\to\cC$. Then 
diagram (\ref{kvadratt0}) becomes
\begin{equation} \label{kvadratt}
        \begin{array}{ccccccc}
\tildetocA&=&\tildetocA&\qeq&\cA &&\\
\!\!\!\! F\downarrow&&\downarrow&&\,\,\,\,\downarrow\ppi &&\\
\cA\nearrow\cB&\qeqleft&\underline{\cC}&\qeq&\cC &,&
\underline{\cC}:=(\cA\nearrow\cB )_{\ppi}\\
        \end{array}
\end{equation}
Here the DG functors
$\cA\leftarrow\tildetocA\to\cA\nearrow\cB$ are
same as in (\ref{kvadrat}).

In \ref{slightlybetter} we will use a slightly different
canonical commutative diagram of DG categories
\begin{equation} \label{kvadratt2}
        \begin{array}{ccccc}
\cA\nearrow\cB&\qeqleft&\overline{\cC}&\qeq&\cC  \\
\cap&&\downarrow&&\,\,\,\,\,\downarrow\ppi^*  \\
\cAto&\leftarrow&\cA^{\circ}\mbox{-resDGmod}&\qeq&
\cA^{\circ}\DGmod \\
     \end{array}
\end{equation}
in which $\ppi^*$ is defined by 
$\ppi^*c(a):=\Hom (\ppi (a),c)$. Here is the construction.

Let us start with the lower row of (\ref{kvadratt2}).
Consider the DG category $\Mmor (\cA^{\circ}\DGmod )$ (see
\ref{Mmor} for the definition of $\Mmor$). Its objects are
triples $(Q,M,f)$, where $Q,M\in\cA^{\circ}\DGmod$ and
$f:Q\to M$ is a closed morphism of degree 0. We define 
$\cA^{\circ}\mbox{-resDGmod}\subset
\Mmor (\cA^{\circ}\DGmod )$ to be the full DG subcategory of
triples $(Q,M,f)$ such that $Q\in\cAto$ and $f$ is a
quasi-isomorphism (so $Q$ is a semi-free resolution of $M$).
In other words, $\cA^{\circ}\mbox{-resDGmod}$ is the DG 
category of {\it resolved DG $\cA^{\circ}$-modules.} The DG
functors $\cA^{\circ}\mbox{-resDGmod}\to\cAto$ and
$\cA^{\circ}\mbox{-resDGmod}\to\cA^{\circ}\DGmod$ are
defined by $(Q,M,f)\mapsto Q$ and $(Q,M,f)\mapsto M$.

We define $\overline{\cC}$ to be the DG category
$(\cA\nearrow\cB )_{\leftarrow\ppi}$ from \ref{RF2new'}. So
the objects of
$\overline{\cC}$ are triples $(Y,X,\psi )$, where
$Y\in\cA\nearrow\cB$, $X\in\cC$, and $\psi :R\ppi (Y)\to X$
is a homotopy equivalence in $\cCto$. The  upper row of
(\ref{kvadratt2}) is defined just as the lower row of
(\ref{kvadratt0}).

The DG functor
$\overline{\cC}\to\cA^{\circ}\mbox{-resDGmod}
\subset\Mmor (\cA^{\circ}\DGmod )$ is defined as follows. To
$(Y,X,\psi )\in\overline{\cC}$ one assigns 
$(Y,\ppi^*X ,\chi)\in\Mmor (\cA^{\circ}\DGmod )$, where 
$\chi :Y\to\ppi^*X$ corresponds to $\psi :R\ppi (Y)\to X$ by
adjointness. This assignment extends in the obvious way to a
DG functor from $\overline{\cC}$ to 
$\Mmor (\cA^{\circ}\DGmod )$.  To show that its image is
contained in 
$\cA^{\circ}\mbox{-resDGmod}$ we have to prove that 
$\chi :Y\to\ppi^*X$ is
a quasi-isomorphism. This follows from the next lemma.

\subsubsection{\bf Lemma}  \label{resolutionlemma}
{\it The natural morphism
$Y\to\ppi^*\ppito (Y)=\ppi^*R\ppi(Y)$,
$Y\in\cB^{\perp}\subset\cAto\subset\cA^{\circ}\DGmod$, is
a quasi-isomorphism.}

\begin{proof}
We will identify $\Ho(\cAto )$ with the derived category
$D(\cA )$ of $\cA^{\circ}$-mo\-dules (so both $Y$ and
$\ppi^*\ppito (Y)$ will be considered as objects of
$\Ho(\cAto )$. The essential image of $\Ho(\cBto )$ in
$\Ho(\cAto )$ will be again denoted by $\Ho(\cBto )$.

It suffices to show that
\begin{equation} \label{conuschtonado}
\Cone (Y\to\ppi^*\ppito (Y))\in\Ho(\cBto )
\end{equation}
for every $Y\in\Ho (\cAto )$ (then for $Y\in\cB^{\perp}$ one
has $\Cone (Y\to\ppi^*\ppito (Y))\in
\Ho(\cBto )\cap\Ho(\cB^{\perp})=0$). 
Proposition \ref{resolutionlemmanew} says that
(\ref{conuschtonado}) holds for $Y\in\Ho (\cA )$. Objects
$Y\in\Ho (\cAto )$ for which (\ref{conuschtonado}) holds
form a triangulated subcategory closed under (infinite)
direct sums. So (\ref{conuschtonado}) holds for all $Y\in\Ho
(\cAto )$ by Lemma \ref{Keller}. 
\end{proof}

\subsubsection{}  \label{glue}

Now let $\underline{\overline{\cC}}$ denote the DG
category $(\cA\nearrow\cB )_{\leftrightarrow\ppi}$ defined in
\ref{RF2new'}. Using the quasi-equivalences 
$\overline{\cC}\qeqleft\underline{\overline{\cC}}\qeq
\underline{\cC}$ one can ``glue'' (\ref{kvadratt} ) and
(\ref{kvadratt2}) and get a canonical commutative diagram
of DG categories
\begin{equation} \label{kvadratt4}
        \begin{array}{ccccc}
\tildetocA&\qeqleft&\tildetocA\times_{\underline{\cC}}
\underline{\overline{\cC}}&\qeq&\cA \\
\downarrow&&\downarrow&&\,\,\downarrow\ppi\\
\cA\nearrow\cB&\qeqleft&\underline{\overline{\cC}}&\qeq&\cC 
\\
\cap&&\downarrow&&\,\,\,\,\,\downarrow\ppi^*  \\
\cAto&\leftarrow&\cA^{\circ}\mbox{-resDGmod}&\qeq
&\cA^{\circ}\DGmod 
\\
     \end{array}
\end{equation}
(the DG functor
$\tildetocA\times_{\underline{\cC}}
\underline{\overline{\cC}}\to\tildetocA$ is a
quasi-equivalence by \ref{diagolemma}, and the DG functor
$\tildetocA\times_{\underline{\cC}}
\underline{\overline{\cC}}\to\cA$ is the composition
$\tildetocA\times_{\underline{\cC}}
\underline{\overline{\cC}}\to\tildetocA\to\cA$, so it is
also a quasi-equivalence).

\section{More on derived DG functors.}
\label{DerivedIsect}

\subsection{} 

Let $\ppi :\cA\to\cC$ be a quotient of a DG category
$\cA$ by a full DG subcategory $\cB\subset\cA$ (so in
(\ref{DGquot}) $\tilde\cA=\cA$ and the DG functor
$\tilde\cA\to\cA$ equals $\id_{\cA}$). Let $F$ be a DG
functor from $\cA$ to a DG category $\cA'$. Under a suitable 
flatness assumption (e.g., if $\cC$ is the DG quotient
$\cA/\cB$ from \S\ref{constr1} and (\ref{FL2}) holds) we
will define notions of the right derived DG functor of $F$,
which correspond to derived triangulated functors
(\ref{tr-derived}) and (\ref{partialderived}). They are
essentially equivalent to those from \ref{RFnew} and
\ref{RF2new} but are based on $\cC$ rather than the DG
quotient $\cA\nearrow\cB$ from \ref{formu}. One can easily
pass from right derived DG functors to left ones by
considering the dual DG categories.

\subsection{}   \label{transversality}
Consider the DG functor
\begin{equation}   \label{ppistarreminder}
\ppi^*:\cC\to\cA^{\circ}\DGmod,\quad
\ppi^*c(a):=\Hom (\ppi (a),c) 
\end{equation}
From now on we assume that the diagram 
$\cC\buildrel\ppi\over\leftarrow\cA\buildrel F\over\to\cA'$
satisfies the following flatness condition: for
all $c\in\Ob\cC$
\begin{equation} \label{transvers}
\mbox{the morphisms
$\ppi^*c\otimes_{\cA}\cA'\to\ppi^*c\Ltensor_{\cA}\cA'$
are quasi-isomorphisms.}
\end{equation}
This condition is satisfied if $\cC$ is the DG quotient
$\cA/\cB$ from \S\ref{constr1} and (\ref{FL2}) holds: in
this case the DG
$\cA^{\circ}$-modules $\ppi^*c$, $c\in\cC$, are
homotopically flat by Lemma \ref{flatness2}(i).

\subsection{}  \label{RF}
We are going to define a DG version of the derived triangulated
functor (\ref{tr-derived}). As a first step, consider the DG
functor
\begin{equation} \label{Rid}
\bR F: \cC\to (\cA')^{\circ}\DGmod
\end{equation}
corresponding to the DG
$\cC\otimes(\cA')^{\circ}$-module $\cC\otimes_{\cA}\cA'$ (see
\ref{long}). (This is only a first step because the homotopy
category of the target of $\bR F$ is not the {\it derived\,}
category of DG $(\cA')^{\circ}$-modules). The isomorphism 
$\cC\otimes_{\cA}\cA'=\HHom_{\cC}\otimes_{\cA}\cA'$ (see
(\ref{long})) shows that $\bR F=\Ind_F\circ\ppi^*$, where
$\ppi^*:\cC\to\cA^{\circ}\DGmod$ is defined by
(\ref{ppistarreminder}) and
$\Ind_F:\cA^{\circ}\DGmod\to (\cA')^{\circ}\DGmod$ is the
induction DG functor (see \ref{induction}). 

The fiber product of $\cC$ and
$(\cA')^{\circ}\mbox{-resDGmod}$ over
$(\cA')^{\circ}\DGmod$ will be denoted by $\cC_{[F]}$ (see
\ref{resolvedDG} for the definition of
$(\cA')^{\circ}\mbox{-resDGmod}$). The DG functor
$\cC_{[F]}\to\cC$ is a quasi-equivalence. We define the
derived DG functor $RF:\cC_{[F]}\to\cAprimeto$ to be the
composition
$\cC_{[F]}\to (\cA')^{\circ}\mbox{-resDGmod}\to\cAprimeto$.
A 2-categorical reformulation of this definition will
be given in Remark (ii) from \ref{Induality}.

Let $\cC_{(F)}$ denote the preimage of
$\cA'\subset\cAprimeto$ under $RF$ (so $\cC_{(F)}$ is a full
DG subcategory of $\cC_{[F]}$). One has
$RF:\cC_{(F)}\to\cA'$.

In \ref{suffer1}-\ref{slightlybetter}  
we will show using (\ref{transvers}) that the above
definitions are reasonable: the DG functor
$RF:\cC_{[F]}\to\cAprimeto$  is essentially equivalent to
the DG functor $RF$ from \ref{RFnew}  and therefore agrees
with the derived triangulated functor (\ref{tr-derived}).
There is a similar relation between $RF:\cC_{(F)}\to\cA'$,
the DG functor from \ref{RF2new}, and the derived
triangulated functor (\ref{partialderived}).

\medskip

\noindent {\bf Remark.} If $k$ is a field or, more
generally, if 
\begin{equation} \label{SF}
\mbox{$\Hom (U,X)$ is a semi-free DG $k$-module for all
$X\in\cA$, $U\in\cB$.}
\end{equation}
then the image of $\bR F:\cC\to (\cA')^{\circ}\DGmod$ is
contained in the full subcategory $\cAprimeto$ of semi-free
DG $(\cA')^{\circ}$-modules (in the case $\cA'=\cA$,
$F=\id_{\cA}$ this is Lemma \ref{flatness2}(ii), and
the general case follows). So if (\ref{SF}) holds then one
does not have to consider $\cC_{[F]}$: one can simply define
$RF:\cC\to\cAprimeto$ to be the DG functor corresponding to
$\bR F$.

\subsection{}  \label{suffer1}
Assuming (\ref{transvers}) we will ``identify''
$RF:\cC_{[F]}\to\cA'$ with the DG functor
$RF:\cA\nearrow\cB\to\cAprimeto$ from \ref{RFnew}. More
precisely, here is a construction of a commutative diagram
\begin{equation} \label{kvadratt6}
        \begin{array}{ccccccc}
\cA\nearrow\cB&\qeqleft&\cC_1&\qeq&\cC_{[F]}&\qeq&\cC \\
\!\!\!\!\!\!\!\!\!\!
RF\downarrow&&\downarrow&&\,\,\,\,\,\,\,\downarrow  RF &&\\
\cAprimeto&=&\cAprimeto&=&\cAprimeto&&\\
     \end{array}
\end{equation}
Put $\cC_1:=\cC_{[\id_{\cA}]}$, so the objects of $\cC_1$
are triples $(c,Q,f)$, where $c\in\cC$,
$Q\in\cAto$, and $f:Q\to\ppi^*c$ is a quasi-isomorphism. The
derived DG functor $R\id_{\cA}:\cC_1\to\cA$, i.e., the DG
functor $\cC_1\to\cAto$ defined by $(c,Q,f)\mapsto Q$,
induces a quasi-equivalence
$\cC_1\qeq\cA\nearrow\cB\subset\cAto$ (see
\ref{corollary}). To define the DG functor
$\cC_1\to\cC_{[F]}$ notice that by the flatness
assumption (\ref{transvers}) the image of the composition 
\begin{displaymath} 
\cC_1=\cC_{[\id_{\cA}]}\to\cA^{\circ}\mbox{-resDGmod}
\mono\Mmor (\cA^{\circ}\DGmod )\to
\Mmor ((\cA')^{\circ}\DGmod )
\end{displaymath}
is contained in
$(\cA')^{\circ}\mbox{-resDGmod}$, so we get a DG functor
$\cC_1=\cC_{[\id_{\cA}]}\to
(\cA')^{\circ}\mbox{-resDGmod}$ whose composition with the
DG functor
$(\cA')^{\circ}\mbox{-resDGmod}$ $\to
(\cA')^{\circ}\DGmod$ equals (\ref{Rid}), i.e., we get a
DG functor $\cC_1\to\cC_{[F]}$.

\subsection{}   \label{slightlybetter}
In fact, one can construct a slightly better diagram
\begin{equation} \label{kvadratt8}
        \begin{array}{ccccccc}
\tildetocA&\qeqleft&\tilde\cA&&\qeq&&\cA \\
\downarrow&&\downarrow&&&&\,\,\downarrow\ppi\\
\cA\nearrow\cB&\qeqleft&\tilde\cC&\qeq&\cC_{[F]}&\qeq&\cC \\
\!\!\!\!\!\!\!\!\!\!
RF\downarrow&&\downarrow&&\,\,\,\,\,\,\,\downarrow  RF &&\\
\cAprimeto&=&\cAprimeto&=&\cAprimeto&&\\
     \end{array}
\end{equation}
To this end, first replace in (\ref{kvadratt6}) $\cC_1$  by
the DG category $\overline{\cC}$ from
(\ref{kvadratt2}) (the right square of (\ref{kvadratt2})
defines a DG functor
$\overline{\cC}\to\cC_1$, which is a
quasi-equivalence because $\overline{\cC}\to\cC$ and
$\cC_1\to\cC$ are). Next, put
$\tilde\cC:=\underline{\overline{\cC}}$ (see \ref{glue} for
the definition of $\underline{\overline{\cC}}$) and replace
$\overline{\cC}$ by $\tilde\cC$. Now the upper
two rows of (\ref{kvadratt4}) yield (\ref{kvadratt8}) with
$\tilde\cA:=\tildetocA\times_{\underline{\cC}}
\underline{\overline{\cC}}$.

\section{Proof of  Theorem \ref{main}.} \label{proofmain}

\subsection{}

We can suppose that (\ref{FL1}) holds (if (\ref{FL2}) holds
replace $\cA$ and $\cB$ by the dual categories). 
It suffices to show that $\Phi$ is fully faithful (this will
imply that $\im \Phi$ is a triangulated subcategory of
$(\cA/\cB)^{\tr}$, but on the other hand 
$\im \Phi\supset\cA/\cB$, so $\Phi$ is essentially
surjective). In other words, it suffices to prove that for
every $X,Y\in\cA^{\pretr}$ and every $i\in\BZ$ the
homomorphism
\begin{equation} \label{1new}
\Ext^i_{\cA^{\tr}/\cB^{\tr}}(X,Y)\to
\Ext^i_{(\cA/\cB)^{\tr}}(X,Y)
\end{equation}
is bijective. It is enough to prove this for $X,Y\in\cA$.

\subsection{}

By (\ref{append1}), the l.h.s. of (\ref{1new}) can be
computed as follows:
\begin{equation} \label{2new}
\Ext^i_{\cA^{\tr}/\cB^{\tr}}(X,Y)=
\rightlimit{(Y\to Z)\in Q_Y}H^i\Hom_{\cA^{\pretr}}(X,Z)\, ,
\end{equation}
where $Q_Y$ is the filtering category of
$\cA^{\tr}$-morphisms $f:Y\to Z$ such that $\Cone (f)$ is
$\cA^{\tr}$-isomorphic to an object of $\cB^{\tr}$. 

The r.h.s. of (\ref{1new}) can be written as
\begin{equation} \label{3new}
\Ext^i_{(\cA/\cB)^{\tr}}(X,Y)=
\rightlimit{(Y\to Z)\in Q_Y}
H^i\Hom_{\cA^{\pretr}/\cB}(X,Z)\, .
\end{equation}
To see this, first notice that the DG functor
$\cA/\cB\to (\cA/\cB)^{\pretr}$ is fully faithful, so
$\Ext^i_{(\cA/\cB)^{\tr}}(X,Y):=
H^i\Hom_{(\cA/\cB)^{\pretr}}(X,Y)=H^i\Hom_{\cA/\cB}(X,Y)$;
then notice that a morphism $Y\to Z$ from $Q_Y$ induces an
isomorphism 
\begin{displaymath}
H^i\Hom_{\cA/\cB}(X,Y)=H^i\Hom_{\cA^{\pretr}/\cB}(X,Y)\iso
H^i\Hom_{\cA^{\pretr}/\cB}(X,Z)
\end{displaymath}
because $\Hom_{\cA^{\pretr}/\cB}(X,U)$ is acyclic for every
$U\in\cB$ (acyclicity is clear since $U$ is homotopy
equivalent to $0$ as an object of $\cA^{\pretr}/\cB$).

\subsection{}

Consider (\ref{1new})  as a morphism from the r.h.s. of
(\ref{2new}) to the r.h.s. of (\ref{3new}). Clearly it is
induced by the morphisms
$\alpha_Z:\Hom_{\cA^{\pretr}}(X,Z)\to
\Hom_{\cA^{\pretr}/\cB}(X,Z)$, 
$Z\in\cA^{\pretr}$. By (\ref{0gr}) each $\alpha_Z$ is
injective and $L_Z:=\Coker\alpha_Z$ is the union of an
increasing sequence of subcomplexes 
$0=(L_Z)_0\subset (L_Z)_1\subset\ldots$ such that
$(L_Z)_n/(L_Z)_{n-1}=\Hom_{\cA^{\pretr}/\cB}^n(X,Z)$
for $n\ge 1$. So to prove that (\ref{1new}) is bijective it
suffices to show that
\begin{displaymath} 
 \rightlimit{(Y\to Z)\in
 Q_Y}H^i\Hom_{\cA^{\pretr}/\cB}^n(X,Z)=0, 
\quad n\ge 1 \, .
\end{displaymath}
For $n\ge 1$ the DG functor
$Z\mapsto\Hom_{\cA^{\pretr}/\cB}^n(X,Z)$ is a direct sum of
DG functors of the form 
$Z\mapsto F_{X,U}\otimes\Hom_{\cA^{\pretr}}(U,Z)$, $U\in\cB$,
where $F_{X,U}$ is a homotopically flat complex of
$k$-modules. Since
\begin{displaymath} 
 \rightlimit{(Y\to Z)\in
 Q_Y}H^i\Hom_{\cA^{\pretr}}(U,Z)=
\Ext^i_{\cA^{\tr}/\cB^{\tr}}(U,Z)=0, \quad U\in\cB
\end{displaymath}
it remains to prove the following lemma.

\subsection{Lemma.} \label{limit}
{\it Let $\{ C_{\alpha}\}$ be a filtering inductive system
of objects of the homotopy category of complexes of
$k$-modules (so each $C_{\alpha}$ is a complex, to each
morphism $\mu:\alpha\to\beta$ there corresponds a morphism
$f_{\mu}:C_{\alpha}\to C_{\beta}$ and $f_{\mu\nu}$ is
homotopy equivalent to $f_{\mu}f_{\nu}$). Suppose that
$\rightlimit{\alpha}H^i(C_{\alpha})=0$ for all $i$. Then for
every homotopically flat complex $F$ of $k$-modules
$\rightlimit{\alpha}H^i(C_{\alpha}\otimes F)=0$.}

\medskip

\noindent {\bf Remark.} This would be obvious if we had a
true inductive system of complexes, i.e., if $f_{\mu\nu}$
were equal to $f_{\mu}f_{\nu}$ (because in this case
$\rightlimit{\alpha}H^i(C_{\alpha})=H^i(C)$, 
$\rightlimit{\alpha}H^i(C_{\alpha}\otimes F)=
H^i(C\otimes F)$, $C:=\rightlimit{\alpha}C_{\alpha}$).
If there are countably many $\alpha$'s then Lemma \ref{limit}
is still obvious bacause we can replace the morphisms
$f_{\mu}$ by homotopy equivalent ones so that
$f_{\mu\nu}=f_{\mu}f_{\nu}$.

\medskip

The proof of Lemma \ref{limit} is based on the following
lemma due to Spaltenstein \cite{Sp}.

\subsection{Lemma}  \label{resolution}

{\it For every complex $F$ of $k$-modules there is a
quasi\--iso\-mor\-phism  $F'\to F$, where $F'$ is a filtering
direct limit of finite complexes of finitely generated free
$k$-modules.}

\begin{proof}
One can take $F'$ to be a semi-free resolution of $F$ (see
\S\ref{Hinich}). Here is a slightly different argument close
to the one from \cite{Sp}. Represent $F$ as a direct limit
of bounded above complexes $F_n$, $n\in\BN$. Let $P_n\to F_n$
be a surjective quasi-isomorphism, where $P_n$ is a bounded
above complex of free $k$-modules. The morphism 
$P_n\to F_{n+1}$ can be lifted to a morphism 
$P_n\to P_{n+1}$. We can take $F'$ to be the direct limit of
the complexes $P_n$ (because each $P_n$ is the union of a 
filtering family of finite complexes of finitely generated 
free $k$-modules).
\end{proof}

\subsection{Proof of Lemma \ref{limit}}  
Let $F$ be as in Lemma \ref{limit}. Choose $F'$
as in Lemma \ref{resolution}. Since Lemma \ref{limit}
holds for $F'$ instead of $F$ it suffices to show that the
map $H^i(C_{\alpha}\otimes F')\to H^i(C_{\alpha}\otimes F)$
is an isomorphism. As $\Cone (F'\to F)$ is homotopically
flat and acyclic this follows from Proposition 5.8 of
\cite{Sp}: if a complex $C$ is homotopically flat and acyclic
then $C\otimes C'$ is acyclic for every complex $C'$ (proof:
by Lemma \ref{resolution} one may assume that $C'$ is either
homotopically flat or acyclic). \hfill\qedsymbol

\section{Proof of Propositions \ref{resolutionlemmanew} and
\ref{comparisonnew}.}
\label{proof3}

\subsection{Proof of Proposition \ref{comparisonnew}}
\label{comparisonproof}

Let $Y\in\Ob\cA$. Then 
\begin{equation} \label{firstlimit}
RF^{\tr}(Y)=\rightlimitquot{(Y\to Z)\in Q_Y} F^{\tr}(Z),
\end{equation}
Here $Q_Y$ is the filtering category of
$\cA^{\tr}$-morphisms $g:Y\to Z$ such that $\Cone (g)$ is
isomorphic to an object of $\cB^{\tr}$.

To compute $RF^{\tr}(Y)$ choose a closed morphism $f:P\to Y$
of degree $0$ with $P\in\cBto$, $\Cone (f)\in\cB^{\perp}$
(i.e., choose a semi-free resolution of the DG
$\cB^{\circ}$-module $b\mapsto\Hom (b,Y), b\in\cB$). Then 
\begin{equation} \label{secondlimit}
H^{0}(RF)^{\tr}(Y)=\rightlimitquot{(W\to P)\in
Q'_P} F^{\tr}(\Cone (W\to Y)),
\end{equation}
where $Q'_P$ is the filtering category of
$\cBto$-morphismsms $W\to P$ with $W\in\cB^{\pretr}$. We have
the functor $\Phi :Q'_P\to Q_Y$ that sends $h:W\to P$ to
$g:Y\to\Cone (fh)$, and (\ref{themorphismnew}) is the
morphism  from the r.h.s. of (\ref{firstlimit}) to the
r.h.s. of (\ref{secondlimit}) corresponding to $\Phi$. It
remains to prove the following lemma.

\subsection{Lemma}           \label{cofinal}
{\it Let $f:P\to Y$ be a closed morphism of degree $0$ with
$Y\in\cA$, $P\in\cBto$, $\Cone (f)\in\cB^{\perp}$. Then the
above functor $\Phi :Q'_P\to Q_Y$ is cofinal.}

\begin{proof}
By the definition of cofinality (see \S8.1 of \cite{GV}), we
have to show that for every $(g:Y\to Z)\in Q_Y$ there
exists $(W\to P)\in Q'_P$ such that the $\cA^{\tr}$-morphism
$Y\to\Cone (W\to Y)$ can be factored through $g$. There is a
distinguished triangle 
$V\buildrel \psi\over\to Y\buildrel g\over\to Z\to V[1]$,
$V\in\cB^{\tr}$, so it suffices to show that $\psi$ is in the
image of the composition
\begin{equation} \label{comp}
\rightlimit{(W\to P)\in Q'_P}\Hom_{\cA^{\tr}}(V,W)\to
\Hom_{\cAto^{\tr}}(V,P)\to\Hom_{\cA^{\tr}}(V,Y).
\end{equation}
This is clear because both maps in (\ref{comp}) are
bijective (the second one is bijective because
$V\in\cB^{\tr}$ and $\Cone (f:P\to Y)\in\cB^{\perp}$).
\end{proof}

\subsection{Proof of Proposition \ref{resolutionlemmanew}}
\label{proofof04}
We will use the convention of \ref{indproversions}: $\cBto$
is identified with its essential image under the induction
DG functor $\cBto\to\cAto$.

To prove that (i)$\Rightarrow$(ii) choose a closed morphism
$f:P\to Y$ of degree $0$ with $P\in\cBto\subset\cAto$, $\Cone
(f)\in\cB^{\perp}$ (i.e., choose a semi-free resolution of
the DG $\cB^{\circ}$-module $b\mapsto\Hom (b,Y), b\in\cB$).
It suffices to show that (\ref{thecone2}) is quasi-isomorphic
to $P[1]$. To this end, consider the commutative diagram
\begin{equation} \label{kvadratt-2}
        \begin{array}{ccc}
\Hom (X,Y)&\stackrel {u_X}\longrightarrow&
\limto\Hom (X,\Cone (W\to Y))\\
v_X\downarrow&&\downarrow\alpha_X\\
\Hom (\ppi (X),\ppi (Y))&\stackrel{\beta_X}\longrightarrow&
\limto\Hom (\ppi (X),\ppi (\Cone (W\to Y))) \\
        \end{array}
\end{equation}
in which the direct limits are over $(W\to P)\in Q'_P$ (see
\ref{comparisonproof} for the definition of $Q'_P$).
Objects of $\ppi (\cB )$ are homotopic to zero, so
$\beta_X$ is a quasi-isomorphism. By (\ref{append1}) and
\ref{cofinal} $\alpha_X$ is also a quasi-isomorphism. So the
DG $\cA^{\circ}$-module $X\mapsto\Cone (v_X)$ is
quasi-isomorphic to the DG $\cA^{\circ}$-module
$X\mapsto\Cone (u_X)$, i.e., to
$P[1]$.

To prove that (ii)$\Rightarrow$(i) consider again the
commutative diagram (\ref{kvadratt-2}). The DG
$\cA^{\circ}$-module $X\mapsto\Cone (u_X)$ is
quasi-isomorphic to
$P[1]$, and $\beta_X$ is a quasi-isomorphism. So if the DG
$\cA^{\circ}$-module $X\mapsto\Cone (v_X)$ is
quasi-isomorphic to an object of $\cBto\subset\cAto$ then
the DG $\cA^{\circ}$-module 
\begin{equation} \label{ourcone}
X\mapsto\Cone (\alpha_X), \quad X\in\cA
\end{equation}
is quasi-isomorphic to some $M\in\cBto\subset\cAto$. Clearly
$M$ is quasi-isomorphic to the restriction of (\ref{ourcone})
to $\cB$. By (\ref{append1}) and \ref{cofinal} one has
\begin{displaymath} 
\rightlimit{(W\to P)\in Q'_P} H^i\Hom (X,\Cone (W\to
Y))=\Ext^i_{\cA^{\tr}/\cB^{\tr}} (X,Y),\quad X,Y\in\cA.
\end{displaymath}
So the restriction of (\ref{ourcone}) to
$\cB$ is acyclic. Therefore $\alpha_X$ is a
quasi-isomor\-phism for all $X\in\cA$. So the canonical map
$\Ext^i_{\cA^{\tr}/\cB^{\tr}} (X,Y)\to\Ext^i_{\cC^{\tr}}
(\ppi (X),$ $\ppi (Y))$ is an isomorphism for all
$X,Y\in\cA$, i.e., the functor
$\cA^{\tr}/\cB^{\tr}\to\cC^{\tr}$ induced by $\ppi$ is fully
faithful. Its essential image is a triangulated subcategory
containing $\Ho (\cC )$, so it equals $\cC^{\tr}$.
\hfill\qedsymbol

\section{Proof of Propositions
\ref{corollary}, \ref{exactnessoftoandfrom} and
\ref{admis}.}

\label{lastproofsect}

\subsection{Proof of Proposition \ref{admis}} 
\label{proofmain2}

Identify $\cAto^{\tr}=\Ho(\cAto )$ with $D(\cA )$ and 
$\cBto^{\tr}=\Ho(\cBto )$ with $D(\cB )$. Then
the embedding $\cBto^{\tr}\to\cAto^{\tr}$ identifies with
the derived induction functor, so it has a right adjoint,
namely the restriction functor. This proves (i). By
adjointness, $(\cBto^{\tr})^{\perp}\subset\Ho(\cAto )$ is the
kernel of the restriction functor, which proves (ii).
Statement (iii)  follows from (i) and (ii). To prove (iv) 
apply Lemma \ref{lemka} in the following situation:
$\cT_0=\cA^{\tr}$, $\cT=\cAto^{\tr}$,
$\cQ_0=\cB^{\tr}$, $\cQ=\cBto^{\tr}$.

\subsection{Proof of Proposition \ref{corollary}}
\label{corollaryproof} 
(a) is a particular case of \ref{exactnessoftoandfrom}(ii).
Here is a direct proof of (a). As $\ppi$ is essentially
surjective it suffices to show that the morphism 
$f:\Ext^n(\ppi (a),c)\to\Ext^n(\ppi^*\ppi (a),\ppi^*c)$
is an isomorphism for every $a\in\cA$ and $c\in\cC$.
Decompose $f$ as
$\Ext^n(\ppi (a),c)=\Ext^n(a,\ppi^*c)
\buildrel{f'}\over{\longrightarrow}
\Ext^n(\ppi^*\ppi (a),\ppi^*c)$, where $f'$ comes from the
morphism $\varphi :a\to\ppi^*\ppi (a)$. By
\ref{resolutionlemmanew}(ii), there is a
distinguished triangle 
\begin{equation} \label{disttriangle}
L\Ind (N)\to
a\buildrel{\varphi}\over{\longrightarrow}\ppi^*\ppi (a)
\to L\Ind (N)[1], \quad N\in D(\cB )\, ,
\end{equation}
where $L\Ind :D(\cB )\to D(\cA )$ is the derived induction 
functor $L\Ind :D(\cB )\to D(\cA )$. As $\ppi^*c$ is 
annihilated by the restriction functor 
$\Res :D(\cA )\to D(\cB )$ we see that
$\Ext^n(L\Ind (N),\ppi^*c)=0$, so $f'$ is an
isomorphism.

Applying $\Res$ to (\ref{disttriangle}) and using the
equalities $\Res\cdot\ppi^*=0$, 
$\Res\cdot L\Ind =\id_{D(\cB )}\,$ 
we get $N=\Res a$ and 
$\ppi^*\ppi (a)\simeq\Cone (L\Ind\Res a\to a)$. 
This im\-plies~(b).
\hfill\qedsymbol

\subsection{Proof of Proposition \ref{exactnessoftoandfrom}}
\label{exactnessoftofromproof}
The derived category of $\cA^{\circ}$-modules identifies with
$\Ho(\cAto )$. The derived induction functor 
$I:\Ho(\cAto )\to\Ho(\cCto )$ is left adjoint to the
restriction functor $R:\Ho(\cCto )\to\Ho(\cAto )$.

By \ref{admis} we can identify $\Ho(\cAto )/\Ho(\cBto )$ with
$\Ho (\cB^{\perp})=(\Ho (\cBto))^{\perp}$. Clearly 
$R(\Ho(\cCto ))\subset\Ho (\cB^{\perp})$. Let
$i:\Ho(\cB^{\perp})\to\Ho(\cCto )$ and
$r:\Ho(\cCto )\to\Ho(\cB^{\perp})$ be the functors
corresponding to $I$ and $R$. It suffices to show that they
are quasi-inverse equivalences. Clearly $i$ is left adjoint
to $r$. So we have the adjunction morphisms $\id\to ri$,
$ir\to\id$, and we have to show that they are isomorphisms.
By \ref{resolutionlemma} the morphism
$\id\to ri$ is an isomorphism. Therefore, the natural
morphism $r\to rir$ is an isomorphism, so the morphism
$rir\to r$ is an isomorphism (because the composition $r\to
rir\to r$ equals $\id$), and finally the morphism
$ir\to\id$ is an isomorphism (because $r$ is conservative,
i.e., if $f$ is a morphism in $\Ho(\cCto )$ such that $r(f)$
is an isomorphism then $f$ is an isomorphism).
\hfill\qedsymbol

\section{Proof of Proposition \ref{preservedbytensoring} and
  Theorem \ref{univ}}
\label{univproofsection}

\subsection{Proof of Proposition \ref{preservedbytensoring}}
\label{preservedbytensoringproof}

Let $M_Y$ denote the DG $\cA^{\circ}$-module
(\ref{thecone2}). Replacing $\ppi :\cA\to\cC$ by
$\ppi\otimes\id_{\cK}:\cA\otimes\cK\to\cC\otimes\cK$ one
gets a similar DG $\cA^{\circ}\otimes\cK^{\circ}$-module
$M_{Y\otimes Z}$ for every $Z\in\cK$. Clearly 
$M_{Y\otimes Z}=M_Y\otimes h_Z$, where $h_Z$ is the image of
$Z$ under the Yoneda embedding $\cK\mono\cK^{\circ}\DGmod$.
As $\cK$ is homotopically flat over $k$ property
\ref{resolutionlemmanew}(ii) for $\ppi :\cA\to\cC$ implies
property \ref{resolutionlemmanew}(ii) for
$\ppi\otimes\id_{\cK}:\cA\otimes\cK\to\cC\otimes\cK$. It
remains to use Proposition \ref{resolutionlemmanew}.
\hfill\qedsymbol

\subsection{Proof of Theorem \ref{univ}}
\label{univproof}

A pair $(\cC ,\ppi )$ satisfying \ref{univ}(ii) is clearly
unique in the sense of ${\bf DGcat}$, and in
\S\S\ref{constr1}-\ref{constr2} we proved the existence of DG
quotient, i.e., the existence of a pair $(\cC,\ppi )$
satsifying \ref{univ}(i). So it remains to show that 
\ref{univ}(i)$\Rightarrow$\ref{univ}(ii).

We will use the definition of $T(\cA ,\cK )$ from 
\ref{flatcase}-\ref{trianggen}. 
One can assume that $\cK$ is homotopically flat over $k$.
So $T(\cA ,\cK )\subset D(\cA^{\circ}\otimes\cK)$,
$T(\cB ,\cK )\subset D(\cB^{\circ}\otimes\cK )$,
$T(\cC ,\cK )\subset D(\cC^{\circ}\otimes\cK )$. We can also
assume that $\ppi\in T(\cA ,\cC )$ comes from a DG functor
$\ppi :\cA\to\cC$ (otherwise replace $\cA$ by one of its
semi-free resolutions and apply \ref{Konts2}). So if
\ref{univ}(i) holds one can apply \ref{preservedbytensoring} 
and \ref{exactnessoftoandfrom}. 
We see that the restriction functor 
$D(\cC^{\circ}\otimes\cK )\to D(\cA^{\circ}\otimes\cK )$ is
fully faithful, and its essential image consists precisely of
objects of $D(\cA^{\circ}\otimes\cK )$ annihilated by the
restriction functor 
$D(\cA^{\circ}\otimes\cK )\to D(\cB^{\circ}\otimes\cK )$.
Property \ref{univ}(ii) follows. \hfill\qedsymbol

\section{Appendix I: Triangulated categories.} \label{triang}

\subsection{Categories with $\BZ$-action and graded
categories} \label{gradedstruct}

Let $C$ be a category with a weak action of $\BZ$, i.e.,
a monoidal functor from $\BZ$ to the monoidal
category $\Funct (C,C)$ of functors $C\to C$ (here $\BZ$ is
viewed as a monoidal category: $\Mor (m,n):=\emptyset$ if
$m\ne n$, $\Mor (n,n):=\{\id_n\}$, $m\otimes n:=m+n$ for
$m,n\in\BZ$). For $c_1,c_2\in C$ put
$\Ext^n(c_1,c_2):=\Mor(c_1,F_n(c_2))$, where $F_n:\cC\to\cC$
is the functor corresponding to $n\in\BZ$. Using the
isomorphism $F_mF_n\iso F_{m+n}$ one gets the 
composition map
$\Ext^m(c_1,c_2)\times\Ext^n(c_2,c_3)\to\Ext^{m+n}(c_1,c_3)$,
so $\cC$ becomes a $\BZ$-graded category. This $\BZ$-graded
category has an additional property: for every $n\in\BZ$ and
$c\in C$ there exists an object $c[n]\in C$ with an
isomorphism $c[n]\iso c$ of degree $n$. Every $\BZ$-graded
category $C$ with this property comes from an essentially
unique weak action of $\BZ$ on $C$.

Suppose that each of the categories $C$ and $C'$ is equipped
with a weak action of $\BZ$. Consider $C$ and $C'$ as graded
categories. Then a graded functor $C\to C'$ (i.e., a functor
between the corresponding graded categories) is the same as a
functor $\Phi :C\to C'$ equipped with an isomorphism
$\Phi\Sigma\iso\Sigma'\Phi$, where $\Sigma\in\Funct (C,C)$
and $\Sigma'\in\Funct (C',C')$ are the images of $1\in\BZ$.

An additive $\BZ$-graded category $C$ is considered as a
plain (non-graded) category by considering elements of
$\bigoplus_n\Ext^n(c_1,c_2)$ (rather than those of
$\bigsqcup_n\Ext^n(c_1,c_2)$) as morphisms $c_1\to c_2$.

All this applies, in particular, to triangulated categories.

\subsection{Quotients} \label{triang1}
The quotient $\cT/\cT'$ of a
triangulated category $\cT$ by a  triangulated subcategory
$\cT'$ is defined to be the localization of $\cT$ by the
multiplicative set $S$ of morphisms $f$ such that $\Cone (f)$
is isomorphic to an object of $\cT'$. The category
$\cT/\cT'$ has a canonical triangulated structure; by
definition, the distinguished triangles of $\cT/\cT'$ are
those isomorphic to the images of the distinguished
triangles of $\cT$. This is due to Verdier \cite{V1, V2}. 

He also proved in \cite{V1, V2} that for every $Y\in\Ob\cT$
the category $Q_Y$ of $\cT$-morphisms $f:Y\to Z$ such that
$\Cone (f)$ is isomorphic to an object of $\cT'$ is
filtering, and for every $Y\in\Ob\cT$ one has an isomorphism
\begin{equation} \label{append1}
\rightlimit{(Y\to Z)\in
Q_Y}\Ext^i_{\cT}(X,Z)\iso\Ext^i_{\cT/\cT'}(X,Y)\, ,
\end{equation}
\subsection{ Remarks}    \label{triang2}
(i) Verdier requires $\cT'$ to be
thick (\'epaisse), which means according to \cite{V2} that
an object of $\cT$ which is (isomorphic to) a direct summand
of an object $\cT'$ belongs to $\cT'$. But the statements
from \ref{triang1} hold without the thickness assumption
because in \S II.2.2 of \cite{V2} (or in \S2.3 of Ch~1 of
\cite{V1}) the multiplicative set $S$ is not required to be
saturated (by Proposition 2.1.8 of \cite{V2}
thickness of $\cT'$ is equivalent to saturatedness of~$S$).

(ii) $\cT/\cT'=\cT/\cT''$, where $\cT''\subset\cT$ is the
smallest thick subcategory containing $\cT'$. So according
to \cite{V2} an object of $\cT$ has zero image in
$\cT/\cT'$ if and only if it belongs to $\cT''$.

(iii) The definitions of thickness from \cite{V1} and
\cite{V2} are equivalent: if $\cT'\subset\cT$ is thick in the
sense of \cite{V2} then according to \cite{V2} $\cT'$ is the
set of objects of $\cT$ whose image in $\cT/\cT'$ is zero, so
$\cT'$ is thick in the sense of \cite{V1}. Direct proofs
of the equivalence can be found in \cite{R} (Proposition 1.3
on p.~305) and \cite{N} (Criterion 1.3 on p.~390).

\subsection{} \label{known}

Let $\cQ$ be a  triangulated subcategory of a
triangulated category $\cT$. Let $\cQ^{\perp}\subset\cT$ be
the right orthogonal complement of $\cQ$, i.e.,
$\cQ^{\perp}$ is the full subcategory of $\cT$ formed by
objects $X$ of $\cT$ such that $\Hom_{\cT}(Y,X)=0$ for all
$Y\in\Ob\cQ$. Then the
morphism $\Hom_{\cT}(Y,X)\to\Hom_{\cT/\cQ}(Y,X)$ is an
isomorphism for all $X\in\Ob\cQ$, $Y\in\Ob\cT$ (see \S6 of
Ch.~I of \cite{V1} and Proposition II.2.3.3 of \cite{V2}). In
particular, the functor $\cQ^{\perp}\to\cT/\cQ$ is fully
faithful. This is a particular case ($\cT_0=\cQ^{\perp}$,
$\cQ_0=0$) of the following lemma. 

\subsection{Lemma}    \label{lemka}

{\it Let $\cQ,\cT_0,\cQ_0$ be  triangulated subcategories of
a triangulated category $\cT$, $\cQ_0\subset\cQ\cap\cT_0$.
Suppose that every morphism from an object of $\cT_0$ to an
object of $\cQ$ factors through an object of $\cQ_0$. Then
the functor $\cT_0/\cQ_0\to\cT/\cQ$ is fully faithful.}

\begin{proof}
The functor $\cT_0/\cQ_0\to\cT/\cQ_0$ is fully faithful by
(\ref{append1}). Our factorization condition implies that
$\Hom_{\cT/\cQ_0}(X,Y)=0$ for all $X\in\Ob\cT_0$, 
$Y\in\Ob Q$. In other words, $\cT_0/\cQ_0$ is contained in
the right orthogonal complement of $\cQ/\cQ_0$ in
$\cT/\cT_0$, so by \ref{known} the functor
$\cT_0/\cQ_0\to(\cT/\cQ_0)/(\cQ/\cQ_0)=\cT/\cQ$ is fully
faithful.
\end{proof}

\subsection{Admissible subcategories}   \label{admissible}
Suppose that a triangulated subcategory $\cQ\subset\cT$ is
strictly full (``strictly'' means that every object
of $\cT$ isomorphic to an object of $\cQ$ belongs to $\cQ$).
Let $\cQ^{\perp}\subset\cT$ (resp. $^{\perp}\cQ\subset\cT$)
be the right (resp. left) orthogonal complement of $\cQ$,
i.e., the full subcategory of $\cT$ formed by objects $X$ of
$\cT$ such that $\Hom (Y,X)=0$ (resp. $\Hom (X,Y)=0$) for
all $Y\in\Ob\cQ$. According to \S1 of \cite{BK2}, 
$\cQ$ is said to be {\it right-admissible} if for each
$X\in\cT$ there exists a distinguished triangle
$X'\to X\to X''\to X'[1]$ with $X'\in\cQ$ and
$X''\in\cQ^{\perp}$ (such a triangle is unique up to unique
isomorphism). As $\cQ^{\perp}$ is thick,
$\cQ$ is right-admissible if and only if the functor
$\cQ\to\cT/\cQ^{\perp}$ is eesentially surjective. $\cQ$ is
said to be {\it left-admissible} if
$\cQ^{\circ}\subset\cT^{\circ}$ is right-admissible. There
is a one to one correspondence between right-admissible
subcategories $\cQ\subset\cT$ and left-admissible
subcategories $\cQ'\subset\cT$, namely
$\cQ'=\cQ^{\perp}$, $\cQ =\null^{\perp}\cQ'$.
According to \S1 of \cite{BK2} and Ch.~1, \S2.6 of \cite{V1}
right-admissibility is equivalent to each of  the following
conditions:

(a) $\cQ$ is thick and the functor $\cQ^{\perp}\to\cT/\cQ$
is essentially surjective (and therefore an equivalence);

(b) the inclusion functor $\cQ\mono\cT$ has a right adjoint;

(c) $Q$ is thick and the functor $\cT\to\cT/\cQ$ has a
right adjoint;

(d) $\cT$ is generated by $\cQ$ and $\cQ^{\perp}$ (i.e., if
$\cT'\subset\cT$ is a strictly full triangulated subcategory
containing $\cQ$ and $\cQ^{\perp}$ then $\cT' =\cT$).

\medskip

\noindent {\bf Remark.} A left or right adjoint of a
triangulated functor is automatically triangulated (see
\cite{Ke-Vo} or Proposition 1.4 of \cite{BK2}).

\section{Appendix II: Semi-free resolutions.} \label{Hinich}

\subsection{Definition} \label{2.1}
A DG $R$-module $F$ over a DG ring $R$ is {\it free} if it is
isomorphic to a direct sum of DG modules of the form $R[n]$,
$n\in\BZ$. A DG $R$-module $F$ is {\it semi-free} if the
following equivalent conditions hold:

1) $F$ can be represented as the union of an increasing
sequence of DG sumbodules $F_i$, $i=0,1,\ldots$, so that
$F_0=0$ and each quotient $F_i/F_{i-1}$ is free;

2) $F$ has a homogeneous $R$-module basis $B$ with the
following property: for a subset $S\subset B$ let $\delta(S)$
be the smallest subset $T\subset B$ such that $d(S)$ is
contained in the $R$-linear span of $T$, then for every
$b\in B$ there is an $n\in\BN$ such that
$\delta^n(\{ b\})=\emptyset$.

A complex of $k$-modules is semi-free if it is semi-free as a
DG $k$-module.

\subsection{Remarks}   \label{rema}
(i) A bounded above complex of free $k$-modules is semi-free.

(ii) Semi-free DG modules were explicitly introduced in
\cite{AH} (according to the terminology of \cite{AH}, a DG
module over a DG algebra $R$ is free if it is freely
generated, as an $R$-module, by homogeneous elements
$e_{\alpha}$ such that $de_{\alpha}=0$, so semi-free is
weaker than free). In fact, the notion of semi-free DG
module had been known to topologists long before \cite{AH}
(see, e.g., \cite{GM}). Semi-free DG modules are also
called ``cell DG modules'' (Kriz--May \cite{KM}) and
``standard cofibrant DG modules'' (Hinich \cite{H}).
In fact, Hinich shows in \S\S 2--3 of \cite{H} that DG
modules over a fixed DG algebra form a closed  model
category with weak equivalences being quasi-isomorphisms and
fibrations being surjective maps. He shows that a DG module
$C$ is cofibrant  (i.e., the morphism $0\to C$ is cofibrant)
if and only if it is a direct summand of a semi-free DG
module. 

(iii) As noticed in \cite{AFH} and \cite{H}, a semi-free
DG module $F$ is {\it homotopically projective}, which means
that for every acyclic DG module $N$ every morphism $f:F\to
N$ is homotopic to~$0$ (we prefer to use the name
``homotopically projective'' instead of Spaltenstein's name
``K-projective''). Indeed, if $\{ F_i\}$ is a filtration on
$F$ satisfying the condition from \ref{2.1} then every
homotopy between $f|_{F_{i-1}}$ and $0$ can be extended to a
homotopy between $f|_{F_i}$ and $0$. This also follows from
Lemma \ref{Keller} applied to the triangulated subcategory
$\cT_N$ of semi-free DG $R$-modules $F$ such that the
complex $\Hom (F,N)$ is acyclic ($\cT_N$ is closed under
arbitrary direct sums and contains $R$).

(iv) By (iii) and Lemma \ref{free} the functor from the
homotopy category of semi-free DG $R$-modules to the derived
category of $R$-modules is an equivalence.

\subsection{Lemma} \label{free}
{\it For every DG module $M$ over a DG algebra $R$ there is a
quasi\--iso\-mor\-phism $f:F\to M$ with $F$ a semi-free DG
$R$-module. One can choose $f$ to be surjective.}

\medskip

The pair $(F,f)$ is constructed in \cite{AH} as the direct
limit of $(F_i,f_i)$ where 
$0=F_0\mono F_1\mono F_2\mono$~\dots, each quotient
$F_i/F_{i-1}$ is free, $f_i:F_i\to M$,
$f_i|_{M_{i-1}}=f_{i-1}$. Given $F_{i-1}$ and
$f_{i-1}:F_{i-1}\to M$ one finds a morphism $\pi :P\to\Cone
(f_{i-1})[-1]$ such that
$P$ is free and $\pi$ induces an epimorphism of the
cohomology groups. $\pi$ defines a morphism 
$f_i:F_i:=\Cone (P\to F_{i-1})\to M$ such that
$f_i|_{F_{i-1}}=f_{i-1}$. The map 
$\Cone (f_{i-1})\to\Cone (f_i)$ induces a zero map of the
cohomology groups, so $\Cone (f)$ is acyclic, i.e., $f$ is
a quasi\--iso\-mor\-phism.

\medskip

\noindent {\bf Remark.} One can reformulate the above proof
of the lemma without using the ``linear'' word ``cone'' (it
suffices to replace ``category'' by ``module'' in the proof
of Lemma \ref{DGresol}).

\subsection{}   \label{semifree}
Hinich \cite{H} proved a version of Lemma \ref{free} for DG
algebras, i.e., DG categories with one object. The case of a
general DG category is similar.

\medskip

\noindent {\bf Definition.} Let $\cA$ be a 
DG category $\cA$ equipped with a DG functor $\cK\to\cA$.
We say that $\cA$ is {\it semi-free over $\cK$\,} if $\cA$
can be represented as the union of an increasing sequence of
DG subcategories $\cA_i$, $i=0,1,\ldots$, so that
$\Ob\cA_i=\Ob\cA$, $\cK$ maps isomorphically onto $\cA_0$,
and for every $i>0$ $\cA_i$ as a graded $k$-category
over $\cA_{i-1}$ (i.e., with forgotten differentials in the
$\Hom$ complexes) is freely generated over $\cA_{i-1}$ by a
family of homogeneous morphisms
$f_{\alpha}$ such that $df_{\alpha}\in\Mor\cA_{i-1}$.

\medskip

\noindent {\bf Definition.} A DG category $\cA$ is {\it
semi-free\,} if it is semi-free over $\cA_{\discr}$, where
$\cA_{\discr}$ is the DG category with
$\Ob\cA_{\discr}=\Ob\cA$ such that the endomorphism DG
algebra of each object of $\cA_{\discr}$ equals $k$ and
$\Hom_{\cA_{\discr}}(X,Y)=0$ if $X,Y$ are different objects
of $\cA_{\discr}$.

\medskip

\noindent {\bf Remarks.} 
1) Semi-free DG categories with one object were introduced in
\cite{H} under the name of ``standard cofibrant" DG
algebras. In fact, Hinich shows in \S\S 2, 4 of \cite{H}
that DG algebras form a closed  model category with weak
equivalences being quasi-isomorphisms and fibrations being
surjective maps. He shows that a DG algebra $R$ is
cofibrant  (i.e., the morphism
$k\to C$ is cofibrant) if and only if $R$ is a retract
of a semi-free DG algebra.

2) $\BZ_-$-graded semi-free DG algebras were considered as
early as 1957 by Tate \cite{T}, and $\BZ_+$-graded ones were
considered in 1973 by Sullivan \cite{Su1,Su2}. Hinich
\cite{H} explained following \cite{Sp} and \cite{AFH} that it
is easy and natural to work with DG algebras without
boundedness conditions.

\subsection{Lemma} \label{DGresol}
{\it For every DG category $\cA$ there exists a semi-free DG
category $\tilde\cA$ with $\Ob\tilde\cA =\Ob\cA$ and a
functor $\Psi:\tilde\cA\to\cA$ such that $\Psi (X)=X$ for
every $X\in\Ob\tilde\cA$ and $\Psi$ induces a surjective
quasi-isomorphism $\Hom (X,Y)\to\Hom (\Psi (X),\Psi (Y))$
for every $X,Y\in\tilde\cA$.} 

\medskip

The proof is same as for DG algebras (\S\S2, 4 of
\cite{H}) and similar to that of Lemma \ref{free}. 
$(\tilde\cA ,\Psi )$ is constructed as the direct limit of
$(\tilde\cA_i,\Psi_i)$ where $\Ob\tilde\cA_i =\Ob\cA$, 
$\tilde\cA_0\mono\tilde\cA_1\mono$~\dots, $\Psi_i:\tilde\cA_i\to\cA$, 
$\Psi_i|_{\tilde\cA_{i-1}}=\Psi_{i-1}$, and the following
conditions are satisfied:

i) $\tilde\cA_0$ is the discrete $k$-category;

ii) for every $i\ge 1$ $\tilde\cA_i$ as a graded $k$-category
is freely generated over $\tilde\cA_{i-1}$ by a family of
homogeneous morphisms $f_{\alpha}$ such that
$df_{\alpha}\in\Mor\tilde\cA_{i-1}$;

iii) for every $i\ge 1$ and  $X,Y\in\Ob\cA$ the morphism
$\Hom_{\tilde\cA_i}(X,Y)\to\Hom_{\cA}(\Psi (X),\Psi (Y))$
induces a surjective map between the sets of the cocycles;

iv) for every $i\ge 2$ and  $X,Y\in\Ob\cA$ the morphism
$\Hom_{\tilde\cA_i}(X,Y)\to\Hom_{\cA}(\Psi (X),\Psi (Y))$ is surjective;

v)  for every $i\ge 1$ and  $X,Y\in\Ob\cA$ every cocycle 
$f\in\Hom_{\tilde\cA_i}(X,Y)$ whose image in 
$\Hom_{\cA}(\Psi (X),\Psi (Y))$ is a coboundary becomes a
coboundary in $\Hom_{\tilde\cA_{i+1}}(X,Y)$.

One constructs $(\tilde\cA_i,\Psi_i)$ by induction.
Note that  if property iii) or iv) holds for some $i$ then it holds for $i+1$,
so after $(\tilde\cA_2,\Psi_2)$ is constructed one only has
to kill cohomology classes by adding new morphisms.

\subsection{Lemma} \label{closedmodel}
{\it If a DG functor $\pi :\tilde\cC\to\cC$ is a surjective
quasi-equivalence (i.e., if $\pi$ induces a
surjection $\Ob\tilde\cC\to\Ob\cC$ and surjective
quasi-isomor\-phisms between the\, $\Hom$ complexes) then
every DG functor from a semi-free DG category $\cA$ to $\cC$
lifts to a DG functor $\cA\to\tilde\cC$. More generally,
for every commutative diagram
\begin{displaymath} 
        \begin{array}{ccc}
\cK&\buildrel{\Phi}\over{\longrightarrow}&\tilde\cC\\
\!\!\!\! \nu\downarrow& &\, \, \downarrow\pi\\
\cR&\buildrel{\Psi}\over{\longrightarrow}&\cC\\
     \end{array}
\end{displaymath}
such that $\cR$ is semi-free over $\cK$ and 
$\pi$ is a surjective quasi-equivalence there exists a
DG functor $\tilde\Psi :\cR\to\tilde\cC$ such that 
$\pi\tilde\Psi=\Psi$ and $\tilde\Psi\nu=\Phi$.}

\medskip

{\bf Remark.} This is one of the closed model
category axioms checked in \cite{H}.

\begin{proof} 
Use the following fact: if $f:A\to B$ is a surjective
quasi-isomorphism of complexes, $a\in A$, $b\in B$,
$f(a)=db$ and $da=0$ then there is an $a'\in A$ such that
$f(a')=b$ and $a=da'$.
\end{proof}

\section{Appendix III: DG modules over DG categories}
\label{DGapp}

Additive functors from a preadditive category $\cA$ to the
category of abelian groups are often called $\cA$-modules
(see \cite{Mi2}). We are going to introduce a similar
terminology in the DG setting. The definitions below are
similar to those of Mitchell \cite{Mi1}. 

\subsection{}
Let $\cA$ be a DG category. A {\it left DG $\cA$-module\,}
is a DG functor from $\cA$ to the DG category of complexes of
$k$-modules. Sometimes left  DG $\cA$-modules will be called
simply DG $\cA$-modules. If $\cA$ has a single object $U$
with $\End_{\cA}U=R$ then a DG $\cA$-module is the same as a
DG $R$-module. A {\it right DG $\cA$-module\,} is a left DG
module over the dual DG category $\cA^{\circ}$. The DG
category of DG $\cA$-modules is denoted by $\cA\DGmod$. In
particular, $\Cold$ is the DG category of complexes of
$k$-modules.

\subsection{}   \label{Alg}
Let $\cA$ be a DG category. Then the complex
\begin{displaymath} 
\Alg_{\cA}:=\bigoplus_{X,Y\in\Ob\cA}\Hom (X,Y)
\end{displaymath}
has a natural DG algebra structure (interpret elements of
$\Alg_{\cA}$ as matrices $(f_{XY})$, $f_{XY}\in\Hom (Y,X)$,
whose rows and columns are labeled by $\Ob\cA$). The DG
algebra $\Alg_{\cA}$ has the following property: every finite
subset of $\Alg_{\cA}$ is contained in
$e\Alg_{\cA}e$ for some idempotent $e\in\Alg_{\cA}$ such
that $de=0$ and $\deg e=0$. We say that a module $M$ over
$\Alg_{\cA}$ is {\it quasi-unital} if every element of $M$
belongs to $eM$ for some idempotent $e\in\Alg_{\cA}$ (which
may be assumed closed of degree $0$ without loss of
generality). If $\Phi$ is a DG $\cA$-module then
$M_{\Phi}:=\bigoplus_{X\in\Ob\cA}\Phi (X)$ is a DG module
over $\Alg_{\cA}$ (to define multiplication write
elements of $\Alg_{\cA}$ as matrices and elements of
$M_{\Phi}$ as columns). Thus we get a DG equivalence between
the DG category of DG $\cA$-modules and that of quasi-unital
DG modules over $\Alg_{\cA}$.

\subsection{} \label{tensorproduct}
Let $F:\cA\to \Cold$ be a left DG $\cA$-module and $G:\cA\to
\Cold$ a right DG $\cA$-module. A {\it DG pairing} $G\times
F\to C$,
$C\in\Cold$, is a DG morphism from the DG bifunctor 
$(X,Y)\mapsto\Hom (X,Y)$ to the DG bifunctor 
$(X,Y)\mapsto\Hom (G(Y)\otimes F(X),C)$. It can
be equivalently defined as a DG morphism $F\to\Hom (G,C)$ or
as a DG morphism $G\to\Hom (F,C)$, where $\Hom (G,C)$ is the
DG functor $X\mapsto\Hom (G(X),C)$, $X\in\cA$. There is a
universal DG pairing $G\times F\to C_0$. We say that $C_0$
is the {\it tensor product} of $G$ and $F$, and we write
$C_0=G\otimes_{\cA}F$. Explicitly, $G\otimes_{\cA}F$ is the
quotient of $\bigoplus_{X\in\cA}G(X)\otimes F(X)$ by the
following relations: for every morphism $f:X\to Y$ in $\cA$
and every $u\in G(Y)$, $v\in F(X)$ one should identify
$f^*(u)\otimes v$ and $u\otimes f_*(v)$. In terms of \S IX.6
of \cite{ML}, $G\otimes_{\cA}F=\int\limits^X G(X)\otimes
F(X)$, i.e.,
$G\otimes_{\cA}F$ is the coend of the functor
$\cA^{\circ}\times\cA\to\Cold$ defined by 
$(Y,X)\mapsto G(Y)\otimes F(X)$. In terms of \ref{Alg},
a DG pairing $G\times F\to C$ is the same as a DG pairing
$M_G\times M_F\to C$, so
$G\otimes_{\cA}F=M_G\otimes_{\Alg_{\cA}}M_F$.

\subsection{Example} For every $Y\in\cA$ one has the
right DG $\cA$-module $h_Y$ and the left DG $\cA$-module
$\tilde h_Y$ defined by $h_Y(Z):=\Hom (Z,Y)$,
$\tilde h_Y(Z):=\Hom (Y,Z)$, $\;Z\in\cA$. One has the
canonical isomorphisms
\begin{equation} \label{tensorfree1}
G\otimes_{\cA}\tilde h_Y=G(Y)\, ,
\end{equation}
\begin{equation} \label{tensorfree2}
h_Y\otimes_{\cA}F=F(Y)
\end{equation}
induced by the maps $G(Z)\otimes\Hom (Y,Z)\to G(Y)$,
$\Hom (Z,Y)\otimes F(Z)\to F(Y)$, $Z\in\cA$.

\subsection{} \label{parameters}
Given DG categories $\cA,\cB,\bar\cB$, a DG
$\cA\otimes\cB$-module $F$, and a DG
$(\cA^{\circ}\otimes\bar\cB )$-module $G$, one defines the DG
$\bar\cB\otimes\cB$-module $G\otimes_{\cA}F$ as follows. We
consider $F$ as a DG functor from $\cB$ to the DG category
of DG $\cA$-modules, so $F(X)$ is a DG $\cA$-module for every
$X\in\cB$. Quite similarly, $G(Y)$ is a DG 
$(\cA )^{\circ}$-module for every
$Y\in\bar\cB$. Now $G\otimes_{\cA}F$ is the DG functor
$Y\otimes X\mapsto G(Y)\otimes_{\cA}F(X)$, $X\in\cB$,
$Y\in\bar\cB$. 

\subsection{} \label{Hom}
Denote by $\HHom_{\cA}$ the DG
$\cA\otimes\cA^{\circ}$-module $(X,Y)\mapsto\Hom (Y,X)$,
$X,Y\in\cA$. E.g., if $\cA$ has a single object and $R$ is
its DG algebra of endomorphisms then $\HHom_{\cA}$ is the DG
$R$-bimodule $R$. For any DG category $\cA$ the isomorphisms
(\ref{tensorfree1}) and (\ref{tensorfree2}) induce
canonical isomorphisms
\begin{equation} \label{unitproperty}
\HHom_{\cA}\otimes_{\cA}F=F,\quad G\otimes_{\cA}\HHom_{\cA}=G
\end{equation}
for every left DG $\cA$-module $F$ and right DG $\cA$-module
$G$ (the meaning of $\HHom_{\cA}\otimes_{\cA}F$ and
$G\otimes_{\cA}\HHom_{\cA}$ was explained in
\ref{parameters}). The isomorphisms (\ref{unitproperty}) are
clear from the point of view of \ref{Alg} because
$M_{\HHom_{\cA}}$ is $\Alg_{\cA}$ considered as a DG bimodule
over itself.

\subsection{} \label{flatness}
A left or right DG $\cA$-module $F:\cA\to\Cold$ is said to be
{\it acyclic} if the complex $F(X)$ is acyclic for every
$X\in\cA$. A left DG $\cA$-module $F$ is said to
be {\it homotopically flat} if $G\otimes_{\cA}F$ is acyclic
for every acyclic right DG $\cA$-module $G$. A right DG
$\cA$-module  is said to be homotopically flat if it is
homotopically flat as a left DG $\cA^{\circ}$-module. It
follows from (\ref{tensorfree1}) and (\ref{tensorfree2})
that $h_Y$ and $\tilde h_Y$ are homotopically flat.

\subsection{} \label{semifreedom}
Let $\cA$ be a DG category. A DG $\cA$-module is said to be 
{\it free} if it is isomorphic to a direct sum of
complexes of the form $\tilde h_X[n]$, $X\in\cA$, $n\in\BZ$.
The notion of semi-free DG $\cA$-module is quite similar to
that of semi-free module over a DG algebra (see \ref{2.1}):
an $\cA$-module $\Phi$ is said to be {\it semi-free} if it
can be represented as the union of an increasing sequence of
DG submodules $\Phi_i$, $i=0,1,\ldots$, so that $\Phi_0=0$
and each quotient $\Phi_i/\Phi_{i-1}$ is free. Clearly a
semi-free DG $\cA$-module is homotopically flat. For every
DG $\cA$-module $\Phi_i$ there is a quasi-isomorphism
$F\to\Phi$ such that $F$ is a semi-free DG $\cA$-module; this
is proved just as in the case that $\cA$ has a single object
(see Lemma \ref{free}). Just as in \ref{rema} one shows that
a semi-free DG $\cA$-module is homotopically projective
(i.e., the complex $\Hom (F,N)$ is acyclic for every
acyclic DG $\cA$-module $N$) and that the functor from the
homotopy category of semi-free DG
$\cA$-modules to the derived category $D(\cA^{\circ})$ of
$\cA$-modules is an equivalence.

\subsection{}  \label{induction}
Let $F:\cA\to \cA'$ be a DG functor between DG categories.
Then we have the {\it restriction\,} DG functor
$\Res_F:\cA'\DGmod\to\cA\DGmod$, which maps a DG
$\cA$-module $\Psi :\cA'\to\Cold$ to $\Psi\circ F$. Sometimes
instead of $\Res_F\Psi$ we write $\Psi$ or ``$\Psi$
considered as a DG $\cA$-module''.

We define the {\it induction\,} functor 
$\Ind_F:\cA\DGmod\to
\cA'\DGmod$ by 
\begin{equation} \label{basechange}
\Ind_F\Phi (Y)=(\Res_Fh_Y)\otimes_{\cA}\Phi\, ,\quad
Y\in\cA'\, .
\end{equation}
or equivalently by 
\begin{equation} \label{basechange2}
\Ind_F\Phi :=\HHom_{\cA'}\otimes_{\cA}\Phi
\end{equation}
(according to
\ref{Hom} $\HHom_{\cA'}$ is a DG
$\cA'\otimes(\cA')^{\circ}$-module, but in
(\ref{basechange2}) we consider it as a DG
$\cA'\otimes \cA^{\circ}$-module). Usually we write 
$\cA'\otimes_{\cA}\Phi$ instead of
$\HHom_{\cA'}\otimes_{\cA}\Phi=\Ind_F\Phi$. 

The DG functor $\Ind_F$ is left adjoint to $\Res_F$. Indeed,
for every DG $\cA'$-module $\Psi$ the complex 
$\Hom_{\cA'\dgmod}(\HHom_{\cA'}\otimes_{\cA}\Phi ,\Psi)$
is canonically isomorphic to
$\Hom_{\cA\dgmod}(\Phi ,
\Hom_{\cA'\dgmod}(\HHom_{\cA'},\Psi))$, and the 
DG $\cA'$-module $\Hom_{\cA'\dgmod}(\HHom_{\cA'},\Psi))$ is
canonically isomorphic to~$\Psi$.

In terms of \ref{Alg} the DG functors $\Res_F$ and $\Ind_F$
correspond to the usual restriction and induction for the
DG algebra morphism $\Alg_{\cA}\to\Alg_{\cA'}$ corresponding
to $\Phi$.

Similar definitions and conventions apply to right DG modules
(in this case 
$\Ind_F\Phi (Y)=\Phi\otimes_{\cA}(\Res_F \tilde h_Y)$,
$\Phi\otimes_{\cA}\cA':=\Phi\otimes_{\cA}\HHom_{\cA'}=
\Ind_F\Phi$).

\subsection{Example}  \label{tensorfree4}
There is a canonical isomorphism
\begin{equation} \label{tensorfree3}
\Ind_F\tilde h_X=\tilde h_{F(X)}, \quad X\in\cA ,
\end{equation}
where $\tilde h_X(Y):=\Hom_{\cA}(X,Y)$, $Y\in\cA$.
This follows either from (\ref{basechange}) and
(\ref{tensorfree1}) or equivalently from (\ref{basechange2})
and (\ref{unitproperty}) (or from the fact that $\Ind_F$ is
the DG functor left adjoint to $\Res_F$). Quite similarly,
there is a canonical isomorphism $\Ind_Fh_X=h_{F(X)}$, which
means that the following diagram is commutative up to
isomorphism:
\begin{equation} \label{commutation2}
        \begin{array}{ccc}
\cA&\longrightarrow&\cA^{\circ}\DGmod\\
\downarrow&&\downarrow\\
\cA'&\longrightarrow&(\cA' )^{\circ}\DGmod\\
     \end{array}
\end{equation}
The horizontal arrows of (\ref{commutation2}) are the
Yoneda embeddings defined by $X\mapsto h_X$, the left
vertical arrow is $F$, and the right one is the induction
functor.

\subsection{Example}  \label{Moritaeq}
Let $\cA$ be a DG category and $F:\cA\to\cA^{\pretr}$ the
embedding. Then $\Res_F:\cA^{\pretr}\DGmod\to\cA\DGmod$ is a
DG equivalence. So $\Ind_F:\cA\DGmod\to\cA^{\pretr}\DGmod$ is
a quasi-inverse DG equivalence.

\subsection{Derived induction} 
As explained, e.g., in \S10 of \cite{BL}, in the situation of
\ref{induction} the functor 
$\Ind_F:\Ho (\cA^{\circ}\DGmod )\to
\Ho ((\cA')^{\circ}\DGmod )$ has a left derived functor 
$L\Ind_F:D(\cA )\to D(\cA' )$, which is
called {\it derived induction.} Derived induction is left
adjoint to the obvious restriction functor 
$D(\cA')\to D(\cA )$. 

By \ref{semifreedom} one can identify $D(\cA )$ with 
$\Ho (\cAto )$, where $\cAto$ is the DG category of semi-free
DG $\cA^{\circ}$-modules. Derived induction viewed as a
functor $\Ho (\cAto )\to\Ho (\cAto' )$ is the obvious
induction functor. Restriction viewed as a functor
$\Ho (\cAto' )\to\Ho (\cAto )$ sends a semi-free DG
$(\cA')^{\circ}$-module to a semi-free resolution of its
restriction to $\cA^{\circ}$.

\subsection{}    \label{catovercattensorcat}
Given DG algebras $A,C,A'$ and DG
morphisms $C\leftarrow A\to A'$ one has the DG
$C\otimes( A')^{\circ}$-module $C\otimes_AA'$.
Quite similarly, given DG categories $\cA ,\cC ,\cA'$
and DG functors $F:\cA\to\cA'$, $G:\cA\to\cC$ one
defines the DG $\cC\otimes(\cA')^{\circ}$-module
$\cC\otimes_{\cA}\cA'$ by
\begin{eqnarray} \label{long} 
&\cC\otimes_{\cA}\cA':=
\HHom_{\cC}\otimes_{\cA}\HHom_{\cA'}=&\\
&=\cC\otimes_{\cA}\HHom_{\cA'}=
\HHom_{\cC}\otimes_{\cA}\cA'=
\cC\otimes_{\cA}\HHom_{\cA}\otimes_{\cA}\cA'\,,
\nonumber
\end{eqnarray}
where $\HHom_{\cC}$ is considered as a
$\cC\otimes \cA^{\circ}$-module and $\HHom_{\cA'}$ as
an $\cA\otimes (\cA')^{\circ}$-module. In other words,
$\cC\otimes_{\cA}\cA'$ is the DG functor
$\cC\times (\cA')^{\circ}\to\Cold$ defined by
\begin{displaymath} 
(X,Y)\mapsto\int^Z\Hom (F(Z),Y)\otimes\Hom (X,G(Z)),\quad
X\in\Ob\cC\,,\;Y\in\Ob\cA'\,,
\end{displaymath}
where the $\int$ symbol denotes the coend (see
\ref{tensorproduct}), so the above ``integral'' is the
tensor product of the right $\cA$-module 
$Z\mapsto\Hom (F(Z),Y)$ and the left $\cA$-module 
$Z\mapsto\Hom (X,G(Z))$. In terms of \ref{Alg}, the DG
module over 
$\Alg_{\cC}\otimes (\Alg_{\cA'})^{\circ}$
corresponding to $\cC\otimes_{\cA}\cA'$ equals
$\Alg_{\cC}\otimes_{\Alg_{\cA}}\Alg_{\cA'}$.

\subsection{}
Given a DG functor $F:\cA\to\cA'$ we say that $\cA'$ is
{\it right $F$-flat\,} (or {\it right homotopically flat
over $\cA$\,}) if the right $\cA$-module $\Res_F h_X$ is
homotopically flat for all $X\in\cA'$; here 
$h_X(Y):=\Hom (Y,X)$, $\;X,Y\in\cA'$. We say that $\cA'$ is
{\it right module-semifree over $\cA$\,} if the right DG
$\cA$-modules $\Res_F h_X$, $X\in\cA'$, are semi-free.
$\cA'$ is said to be {\it left $F$-flat\,} (or {\it left
homotopically flat over $\cA$\,}) if the left $\cA$-module
$\Res_F\tilde h_X$ is homotopically flat for all $X\in\cA'$;
here $\tilde h_X(Y):=\Hom (X,Y)$,
$\;X,Y\in\cA'$. If $\cA'$ is right homotopically flat over
$\cA$ then the induction functor $\Ind_F$ maps acyclic left
DG $\cA$-modules to acyclic left DG $\cA'$-modules. The
previous sentence remains true if ``left'' and ``right'' are
interchanged.

\subsection{Lemma}     \label{flatness2}
{\it Let $\cA$ be a DG category and $\cB\subset\cA$ a full DG
subcategory. 

(i) If (\ref{FL2}) holds then $\cA/\cB$ is right
homotopically flat over $\cA$.

(ii) If (\ref{SF}) holds then $\cA/\cB$ is right
module-semifree over $\cA$.}

\begin{proof}
We will only prove (i) (the proof of (ii) is similar). We
have to show that for every $Y\in\cA$ the functor
$\Psi_Y:\cA^{\circ}\to\Cold$ defined by
$\Psi_Y(X)=\Hom_{\cA/\cB}(X,Y)$ is a homotopically flat right
$\cA$-module. By (\ref{0gr}), there is a
filtration $\Psi_Y=\bigcup_n\Psi_Y^n$, 
$\Psi_Y^n\subset \Psi_Y^{n+1}$, such that $\Psi_Y^0=h_Y$ and
$\Psi_Y^n/\Psi_Y^{n-1}=\bigoplus_{U\in\cB}C_U^n\otimes h_U$
for every $n>0$, where $C_U^n$ is the direct sum of complexes
\begin{displaymath}
\Hom_{\cA}(U_1,U_2)\otimes\ldots\Hom_{\cA}(U_{n-1},U_n)
\otimes\Hom_{\cA}(U_n,Y),\; U_i\in\cB, U_1=U
\end{displaymath}
It remains to notice that for every $Y\in\cA$ the right
$\cA$-module $h_Y$ is homotopically flat (see
\ref{flatness}) and by (\ref{FL2}) the complexes $C^n_U$ are
homotopically flat.
\end{proof}

\subsection{Quasi-representability}  \label{qrep}

Let $\cA$ be a DG category. We have the DG functor from
$\cA$ to the DG category of DG $\cA^{\circ}$-modules defined
by $X\mapsto h_X$.

\subsubsection{\bf Definition}
A DG $\cA^{\circ}$-module $\Phi$ is {\it quasi-representable}
if there is a quasi-isomorphism $f:h_X\to\Phi$ for some
$X\in\cA$. 

\medskip

\noindent {\bf Remark.} By \ref{semifreedom}, for every DG
$\cA^{\circ}$-module $\Phi$ there exists a semi-free
resolution $\pi :\bar\Phi\to\Phi$ (i.e., $\bar\Phi$ is
semi-free and $\pi$ is a quasi-isomorphism), and
the homotopy class of $\bar\Phi$ does not depend
on the choice of $(\bar\Phi,\pi )$. So $\Phi$ is 
quasi-representable if and only if this class contains $h_X$
for some $X\in\cA$.

\subsubsection{\bf Lemma}  \label{stupid}
{\it $\Phi$ is quasi-representable if and only if the
graded functor $H^{\bcdot}\Phi:(\Hodot (\cA
))^{\circ}\to$\{graded $k$-modules\} is representable.}

\begin{proof}
We only have to prove the ``if'' statement. Suppose
$H^{\bcdot}\Phi$ is represented by $(X,u)$, $X\in\Ob\cA$,
$u\in H^0\Phi (X)$. Our $u$ is the cohomology class of some 
$\tilde u\in\Phi (X)$ such that $d\tilde u=0$, 
$\deg\tilde u=0$. Then $\tilde u$ defines a closed morphism
$f:h_X\to\Phi$ of degree 0 such that for every $Y\in\cA$ the
morphism $H^{\bcdot}h_X(Y)\to H^{\bcdot}\Phi (Y)$ is an
isomorphism, so $f$ is a quasi-isomorphism.
\end{proof}

\subsubsection{}           \label{fixquasirep}
Let $\cA'\subset\cA^{\circ}\DGmod$ be the full DG subcategory
of quasi-repre\-sent\-able DG modules. We have the DG
functors
$\cA\leftarrow\cA''\buildrel{\pi}\over{\longrightarrow}\cA'$,
where $\cA''$ is the DG category whose objects are triples
consisting of an object $Y\in\cA$, a DG
$\cA^{\circ}$-module $\Psi$, and a
quasi-isomorphism $h_Y\to\Psi$ (more precisely, 
$\cA''$ is the full DG subcategory of the DG category
$\cA^{\circ}\mbox{-resDGmod}$ from \ref{resolvedDG} which is
formed by these triples). Clearly $\pi$ is a surjective
quasi-equivalence.

\subsubsection{\bf Quasi-corepresentability}
We say that a DG $\cA$-module $\Phi$ is {\it
quasi-corepresentable} if there is a quasi-isomorphism
$f:\tilde h_X\to\Phi$ for some
$X\in\cA$, i.e., if $\Phi$ is representable as
a DG $(\cA^{\circ})^{\circ}$-module

\section{Appendix IV: The diagonal DG categories}
\label{diagapp}

\subsection{}     \label{diag}
Given topological spaces $M',M''$ mapped to a space $M$, one
has the ``homotopy fiber product'' 
$(M'\times M'')\times_{M\times M}\Delta_M^h$, where
$\Delta_M^h$ is the ``homotopy diagonal'', i.e., the space of
paths $[0,1]\to M$ ($\gamma\in\Delta_M^h$ is mapped to
$(\gamma (0),\gamma (1))\in M\times M$). In the same spirit,
given a DG category $\cC$ it is sometimes useful to replace
the naive diagonal $\Delta_{\cC}\subset\cC\times\cC$ by one
of the following DG categories $\rightDelta_{\cC}$,
$\leftDelta_{\cC}$, $\leftrightDelta_{\cC}$, each of them
equipped with a DG functor to $\cC\times\cC$. We define 
 $\rightDelta_{\cC}$ to be the full DG subcategory of the DG
category $\Mmor\cC$ from \ref{Mmor} that consists of triples
$(X,Y,f)$ such that $f$ is a homotopy equivalence; the DG
functor $\rightDelta_{\cC}\to\cC\times\cC$ is defined by
$(X,Y,f)\mapsto (X,Y)$. We define $\leftDelta_{\cC}$ to be
the same full DG subcategory of $\Mmor\cC$, but the DG
functor $\leftDelta_{\cC}\to\cC\times\cC$ is defined by
$(X,Y,f)\mapsto (Y,X)$. 

Finally, define $\leftrightDelta_{\cC}$ to be the DG
category $A_{\infty}\mbox{\rm -funct} (\bI_2,\cC)$  of
{\it $A_{\infty}$-functors\,} $\bI_2\to\cC$,
where $\bI_n$ denotes the $k$-category freely
generated by the category $\bJ_n$ with objects $1,\ldots ,n$
and precisely one morphism with any given source and target.
Here the word ``$A_{\infty}$-functor'' is understood in the
``strictly unital'' sense (cf. \S3.5 of \cite{Ke3}
or \S3.1 of \cite{L-H}; according to \cite{K2, KS, L-H, L}
there are several versions of the notion of
$A_{\infty}$-functor which differ on how an $A_{\infty}$
analog of the axiom $F(\id )=\id$ in the definition of usual
functor is formulated; the difference is inessential for our
purposes and for any reasonable purpose). So an
$A_{\infty}$-functor $\bI_2\to\cC$ is a DG functor
$\cD_2\to\cC$, where $\cD_2$ is a certain DG category with
$\Ob\cD_2=\{ 1,2\}$, which is freely generated (as a graded
$k$-category, i.e., after one forgets the
differential) by morphisms
$f_{12}:1\to 2$ and $f_{21}:2\to 1$ of degree 0, morphisms
$f_{121}:1\to 1$ and $f_{212}:2\to 2$ of degree -1,
morphisms $f_{1212}:1\to 2$ and $f_{2121}:2\to 1$ of degree
-2, etc. One has $df_{12}=0=df_{21}$,
$df_{121}=f_{21}f_{12}-1$, $df_{212}=f_{12}f_{21}-1$,
and  we do not need explicit formulas for the differential
of $f_{1212}$, $f_{2121}$, etc.

\subsection{} \label{leftrighttoright}
Let $e_{ij}$ be the unique
$\bJ_2$-morphism $i\to j$, $i,j\in\{ 1,2\}$. Let
$\bI'_2\subset\bI_2$ denote the $k$-subcategory
generated by $e_{12}$. Then 
$A_{\infty}\mbox{\rm -funct} (\bI'_2,\cC)$ identifies with
$\Mmor\cC$, so we get a canonical DG functor
$\leftrightDelta_{\cC}\to\rightDelta_{\cC}\subset\Mmor\cC$.
There is a similar DG functor
$\leftrightDelta_{\cC}\to\leftDelta_{\cC}$.

\subsection{Lemma} \label{diagolemma}
{\it For every DG category $\cK$ equipped with a DG functor
$\cK\to\rightDelta_{\cC}$ the DG functor 
$\cK\times_{\rightDelta_{\cC}}\leftrightDelta_{\cC}\to\cK$ is
a quasi-equivalence . Same is true if 
$(\leftrightDelta_{\cC},\rightDelta_{\cC})$ is replaced by
$(\rightDelta_{\cC},\cC )$, 
$(\leftrightDelta_{\cC},\leftDelta_{\cC})$
$(\leftrightDelta_{\cC},\rightDelta_{\cC})$, or 
$(\leftDelta_{\cC},\cC )$.}

In other words, the lemma says that the DG functors
$\leftrightDelta_{\cC}\to\rightDelta_{\cC}\to\cC$ are
quasi-equivalences and this remains true after any ``base
change'' in the sense of \ref{fiberprod}.

\begin{proof}
The DG functors
$\leftrightDelta_{\cC}\to\rightDelta_{\cC}\to\cC$ 
induce surjections of Hom complexes (this follows from the
definition of these complexes, see \cite{K2, KS, L-H, L}).
So it suffices to show that they are quasi-equivalences and
induce surjections 
$\Ob\leftrightDelta_{\cC}\to\Ob\rightDelta_{\cC}\to\Ob\cC$.
Both statements are clear for $\rightDelta_{\cC}\to\Ob\cC$.
The DG functor $F:\leftrightDelta_{\cC}\to\cC$ is the
DG functor 
\begin{displaymath} 
A_{\infty}\mbox{\rm -funct} (\bI_2,\cC)\to
A_{\infty}\mbox{\rm -funct} (\bI_1,\cC)
\end{displaymath}     
that comes from a functor $i:\bI_1\to\bI_2$ induced
by an embedding $\bI_1\mono\bI_2$. $F$ is a quasi-equivalence
because $i$ is an equivalence (more generally, if all the
Hom complexes of DG categories
$\cA_1,\cA_2$ are semi-free DG $k$-modules then a
quasi-equivalence $\cA_1\qeq\cA_2$ induces a
quasi-equivalence 
$A_{\infty}\mbox{\rm -funct} (\cA_2,\cC)\qeq
A_{\infty}\mbox{\rm -funct} (\cA_1,\cC)$: this follows
from \ref{Konts3} because the functor 
$T(\cA_2 ,\cC)\to T(\cA_1 ,\cC )$ is an equivalence).

Finally, let us prove the surjectivity of the map
$\Ob\leftrightDelta_{\cC}\to\Ob\rightDelta_{\cC}$ essentially
following \cite{K2} (where a slightly weaker
statement is formulated). We will prove a formally
more general statement. Let $e_{ij}$ and
$\bI'_2\subset\bI_2$ have the same meaning as in
\ref{leftrighttoright}. Suppose that the embedding
$\bI'_2\mono\bI_2$ (considered as a DG functor between DG
categories) is decomposed as
$\bI'_2\mono\cR\to\bI_2$, where
$\Ob\cR=\Ob\bI_2=\bI'_2=\{ 1,2\}$ and $\cR$ is semi-free over
$\bI'_2$ (see \ref{semifree}). Let $F:\bI'_2\to\cC$ be a DG
functor such that $F(e_{12})$ is a homotopy equivalence.
Then we will show that $F$ extends to a DG functor
$G:\cR\to\cC$ (to prove the surjectivity of the  map
$\Ob\leftrightDelta_{\cC}\to\Ob\rightDelta_{\cC}$ put
$\cR=\cD_2$). We will do this by decomposing $F$ as
\begin{equation} \label{decompos}
\bI'_2\buildrel{\Phi}\over{\longrightarrow}\cR'\to\cC,
\quad      \Ho^{\bcdot}(\cR')=\bI_2
\end{equation}
(here the equality $\Ho^{\bcdot}(\cR')=\bI_2$ means that the
functor $\bI'_2=\Ho^{\bcdot}(\bI'_2 )\to\cR'$ extends to an
isomorphism $\bI_2\iso\cR'$). Such a decomposition allows to
extend $F$ to a DG functor $G:\cR\to\cC$: first reduce to
the case that all $\Ext^n$ groups in $\cR'$ vanish for $n>0$
(otherwise replace $\cR'$ by a suitable DG subcaregory), then
one has a commutative diagram
\begin{displaymath} 
        \begin{array}{ccc}
\bI'_2&\buildrel{\Phi}\over{\longrightarrow}&\cR'\\
\!\!\!\! \nu\downarrow& &\, \, \downarrow\pi\\
\cR&\longrightarrow&\bI_2\\
     \end{array}
\end{displaymath}
with $\pi$ being a surjective quasi-equivalence, and it
remains to decompose $\Phi$ as
$\bI'_2\buildrel{\nu}\over{\longrightarrow}\cR\to\cR'$ by
applying \ref{closedmodel}.

Here are two ways to construct a decomposition
(\ref{decompos}). The first way is, essentially, to
construct an $\cR'$ independent on $\cC$ and $F:\bI\to\cC$ by
slightly modifying
$\bI'_2$. The second construction seems simpler to me, but
it gives an $\cR'$ which depends on $\cC$ and $F:\bI\to\cC$.

(i) Our $\bI'_2$ equals the DG category $\cA_0$ from 
\ref{vogt1}. Let $\cR'$ be the DG category 
$(\cA/\cB )_0\subset\cA/\cB$ from \ref{vogt1}. One gets a DG
functor $\cR':=(\cA/\cB )_0\to\cC$ and, in fact, a DG functor
$\cA/\cB\to\cC^{\pretr}$ as follows. First extend
$F:\cA_0:=\bI'\to\cC$ to a DG functor
$F^{\pretr}:\cA:=(\bI'_2)^{\pretr}\to\cC$. Then 
$F^{\pretr}$ sends the unique object of $\cB$ to a
contractible object $Y\in\cC^{\pretr}$. A choice of a
homotopy between $\id_Y$ and 0 defines a DG functor 
$\cA/\cB\to\cC^{\pretr}$. By Lemma
\ref{contractiblequot}, $\Ho^{\bcdot}(\cR')=\bI'_2$.

(ii) Notation: given a DG category $\cA$ and
$a\in\Ob\cA$ one defines $\cA/a$ to be the fiber
product in the Cartesian square
\begin{displaymath} 
        \begin{array}{ccc}
\cA/a&\to&\Mmor\cA\\
\downarrow& &\; \downarrow t\\
*&\buildrel{i_a}\over{\longrightarrow}&\cA\\
     \end{array}
\end{displaymath}
where $\Mmor\cA$ is the DG category from \ref{Mmor}, $t$
sends an $\cA$-morphism to its target, $*$ is the DG category
with one object whose endomorphism algebra equals $k$ and
$i_a:*\to\cA$ maps the object of $*$ to $a$. Decompose
$F:\bI'_2\to\cC$ as $F=s\bar F$, where $s:\cC/F(2)\to\cC$
sends a $\cC$-morphism to its source and 
$\bar F:\bI'_2\to\cC/F(2)$ is the composition of the DG
functor $\bI'_2\to \bI_2/2$ that sends $i\in\{ 1,2\}$ to the
unique $\bJ'_2$-morphism $e_{i2}:i\to 2$ and the DG
functor $\bI_2/2\to\cC/F(2)$ corresponding to $F:\bI_2\to\cC$
(here $\bI_2$ is considered as a DG category). Now define
$\cR'$ from (\ref{decompos}) as follows: 
$\Ob\cR':=\Ob\bI'_2=\{ 1,2\}$, 
$\Hom (j_1,j_2)=\Hom (\bar F(j_1),\bar F(j_2))$ for
$j_1=j_2\in\Ob\cR':=\Ob\bI'_2$, and composition in $\cR'$
comes from composition in $\cC/F(2)$. We have a canonical
decomposition of $\bar F$ as $\bI'_2\to\cR'\to\cC/F(2)$, and
to get (\ref{decompos}) one uses $s:\cC/F(2)\to\cC$. To show
that $\Ho^{\bcdot}(\cR')=\bI_2$ use that $F(e_{i2})$ is a
homotopy equivalence.
\end{proof}

\section{Appendix V: The 2-category of DG categories}
\label{2catofDG}

In \ref{flatcase}-\ref{trianggen} 
we recall the definition of the 2-category of DG categories
used by Keller in \cite{Ke4}, and in
\ref{Konts1}-\ref{Konts3}  we mention a different approach
used by Kontsevich. We prefer to work with the weak notion
of 2-category due to B\'enabou.  The definition and basic
examples of 2-categories can be found in \cite{Be} or
Ch.~XII of
\cite{ML}, where they are called ``bicategories''. Let us
just recall that we have to associate to each two DG
categories $\cA_1, \cA_2$ a  category $T(\cA_1,\cA_2)$ and
to define the composition functors 
$T(\cA_1,\cA_2)\times T(\cA_2,\cA_3)\to T(\cA_1,\cA_3)$. The
2-category axioms say that composition should be weakly
associative and for every DG category $\cA$ there is a weak 
unit object in $T(\cA ,\cA )$. The meaning of ``weak'' is
clear from the following example: a 2-category with one
object is the same as a monoidal category.

The 2-category of DG categories is only the tip of the
``iceberg'' of DG categories. In \ref{DGmodels} we make some
obvious remarks regarding the whole iceberg, but its detailed
description is left to the experts (see \ref{experts}).

\subsection{Flat case} \label{flatcase}
First let us construct the 2-category {\bf FlatDGcat} of flat
DG categories (``flat'' is a shorthand for ``homotopically
flat over $k$'', see \ref{hoflat}). Define 
$T(\cA_1,\cA_2)\subset D(\cA_1^{\circ}\otimes\cA_2)$ to be
the full subcategory of quasi-functors in the sense of
\S7 of \cite{Ke4} (see also \cite{Ke4'}). According to
\cite{Ke4}, a {\it quasi-functor} from $\cA_1$ to $\cA_2$ is
an object $\Phi\in D(\cA_1^{\circ}\otimes\cA_2)$ such that
for every $X\in\cA_1$ the object $\Phi (X)\in D(\cA_2)$
belongs to the essential image of the Yoneda embedding
$\Ho (\cA_2 )\to D(\cA_2)$ (here $\Phi (X)$ is the
restriction of $\Phi :\cA_1\otimes\cA_2^{\circ}\to k\DGmod$
to $\{ X\}\otimes\cA_2=\cA_2$). In other words, an object
of $D(\cA_1^{\circ}\otimes\cA_2)$ is a quasi-functor if it
comes from a DG functor from $\cA_1$ to the full subcategory
of quasi-representable DG $\cA_2^{\circ}$-modules
(``quasi-representable'' means ``quasi-isomorphic to a
representable DG $\cA_2^{\circ}$-module'', see \ref{qrep}).
The composition of 
$\Phi\in D(\cA_1^{\circ}\otimes\cA_2)$ and
$\Psi\in D(\cA_2^{\circ}\otimes\cA_3)$ is defined to be
$\Phi\Ltensor_{\cA_2}\Psi$, and the associativity isomorphism
is the obvious one.

$D(\cA_1^{\circ}\otimes\cA_2)$ is a graded $k$-category (the
morphisms $\Phi_1\to\Phi_2$ of degree $n$ are the elements
of $\Ext^n (\Phi_1 ,\Phi_2 )$). This structure induces a
structure of graded $k$-category on $T(\cA_1,\cA_2)$.

\subsection{Remark}  \label{triangflat}
If $\cA_2$ is pretriangulated in the sense of \ref{pretr}
then the subcategory 
$T(\cA_1,\cA_2)\subset D(\cA_1^{\circ}\otimes\cA_2)$ is
triangulated.

\subsection{General case}   \label{generalcase}
It suffices to define for every DG category $\cA$ a
2-functor $\daleth :S_{\cA}\to {\bf FlatDGcat}$,
where ${\bf FlatDGcat}$ is the 2-category of flat DG
categories and $S_{\cA}$ is a non-empty 2-category such that
for every $s_1,s_2\in\Ob S_{\cA}$ the category of 1-morphisms
$s_1\to s_2$ has one object and one morphism (``$\daleth$''
is the Hebrew letter Dalet). We define 
$\Ob S_{\cA}$ to be the class of all flat resolutions of
$\cS$ (by \ref{DGresol}, $\Ob S_{\cA}\ne\emptyset$).
$\daleth$ sends each $\tilde\cA\in\Ob S_{\cA}$ to itself
considered as an object of ${\bf FlatDGcat}$. The unique
1-morphism from
$\tilde\cA_1\in\Ob S_{\cA}$ to $\tilde\cA_2\in\Ob S_{\cA}$ is
mapped by $\daleth$ to 
$\HHom_{\tilde\cA_1,\tilde\cA_2}\in
T(\tilde\cA_1,\tilde\cA_2)\subset
D(\tilde\cA_1^{\circ}\otimes\tilde\cA_2)$, 
where the DG $\tilde\cA_1\otimes\tilde\cA_2^{\circ}$-module
$\HHom_{\tilde\cA_1,\tilde\cA_2}$ is defined by 
\begin{equation}     \label{HHom}
 (X_1,X_2)\mapsto\Hom (\pi_2(X_2),\pi_1(X_1)),\quad 
X_i\in\tilde\cA_i
\end{equation}
and $\pi_i$ is the DG functor $\tilde\cA_i\to\cA$. To
define $\daleth$ one also has to specify a quasi-isomorphism
\begin{equation}     \label{HomHomtoHom}
\HHom_{\tilde\cA_1,\tilde\cA_2}\Ltensor_{\cA_2}
\HHom_{\tilde\cA_2,\tilde\cA_3}\to
\HHom_{\tilde\cA_1,\tilde\cA_3}
\end{equation}
for every three resolutions
$\tilde\cA_i\to\cA$. It comes from the composition morphism
$\HHom_{\tilde\cA_1,\tilde\cA_2}\otimes_{\cA_2}
\HHom_{\tilde\cA_2,\tilde\cA_3}\to
\HHom_{\tilde\cA_1,\tilde\cA_3}$.

\subsection{}  \label{trianggen}
Each $T(\cA_1,\cA_2)$ is equipped with a graded $k$-category
structure, and if $\cA_2$ is pretriangulated then
$T(\cA_1,\cA_2)$ is equipped with a triangulated structure.
We already know this if $\cA_1$ and $\cA_2$ are flat (see
\ref{flatcase}-\ref{triangflat}), and in the general case we
get it by transport of structure via the equivalence
$T(\tilde\cA_1,\tilde\cA_2)\to T(\cA_1,\cA_2)$ corresponding
to flat resolutions $\tilde\cA_1\to\cA_1$ and
$\tilde\cA_1\to\cA_2$.

\subsection{Remarks}    \label{triangrem}
(i) $T(\cA_1,\cA_2)$ is a full
subcategory of the following triangulated category
$D(\cA_1^{\circ}\Ltensor\cA_2)$ equipped with a triangulated
functor $R:D(\cA_1^{\circ}\otimes\cA_2)\to
D(\cA_1^{\circ}\Ltensor\cA_2)$, which is an equivalence if
$\cA_1$ or $\cA_2$ is flat. The objects of
$D(\cA_1^{\circ}\Ltensor\cA_2)$ are triples
$(\tilde\cA_1,\tilde\cA_2, M)$, where $\tilde\cA_i$ is a
flat resolution of $\cA_i$ and $M\in
D(\tilde\cA_1^{\circ}\otimes\tilde\cA_2)$. Morphisms of
degree $n$ from $(\tilde\cA_1,\tilde\cA_2, M)$ to
$(\tilde\cA'_1,\tilde\cA'_2, M')$ are elements of
$\Ext^n_{\tilde\cA'_1\otimes (\tilde\cA'_2)^{\circ}}
((\HHom_{\tilde\cA'_1,\tilde\cA_1}\otimes
\HHom_{\tilde\cA_2,\tilde\cA'_2})
\otimes_{\tilde\cA_1\otimes\tilde\cA_2^{\circ}}M,M')$. One
defines composition in $D(\cA_1^{\circ}\Ltensor\cA_2)$ and
$R:D(\cA_1^{\circ}\otimes\cA_2)\to
D(\cA_1^{\circ}\Ltensor\cA_2)$
in the obvious way.

(ii) $D(\cA^{\circ}\Ltensor\cA )$ equipped with the functor
$\Ltensor_{\cA}$ is a monoidal category.
$\HHom_{\cA}:=\HHom_{\cA ,\cA}$ viewed as an object of
$D(\cA^{\circ}\Ltensor\cA )$ is a unit object.

\subsection{Ind-version and duality} \label{Induality} We
are going to define an involution $\circ$ of the 2-category
${\bf DGcat}$ which  preserves the composition of
1-morphisms, reverses that of 2-morphisms, and sends each
$\cA\in {\bf DGcat}$ to $\cA^{\circ}$. 

To define it at the level of 1-morphisms and 2-morphisms
consider the 2-category ${\bf DGcat}_{\ind}$ whose objects
are DG categories, as before, but the category
$\Tto (\cA ,\cK )$ of 1-morphisms from a DG category $\cA$
to a DG category $\cK$ equals $D (\cA^{\circ}\Ltensor\cK )$
(1-morphisms are composed in the obvious way). Clearly 
${\bf DGcat}\subset {\bf DGcat}_{\ind}$. The DG category
${\bf DGcat}_{\ind}$ has a canonical involution $\bullet$
which  reverses the composition of 1-morphisms and preserves
that of 2-morphisms: at the level of objects one has
$\cA^{\bullet}:=\cA^{\circ}$, and to define $\bullet$ at the
level of 1-morphisms and 2-morphisms one uses the obvious
equivalence between $\Tto (\cA ,\cK )$ and 
$\Tto (\cK^{\circ} ,\cA^{\circ})$.

Now it is easy to see that each 
$F\in T(\cA ,\cK )\subset\Tto (\cA ,\cK )$ has a right
adjoint $F^*\in\Tto (\cK ,\cA )$ and 
$(F^*)^{\bullet}\in T(\cA^{\circ},\cK^{\circ})\subset
\Tto (\cA^{\circ},\cK^{\circ})$. So putting
$F^{\circ}:=(F^*)^{\bullet}$ one gets the promised involution
of ${\bf DGcat}$.

\medskip

\noindent {\bf Remarks.} (i) It is easy to show that if 
$\cK\in {\bf DGcat}$ is pretriangulated and $\Ho (\cK)$ is
Karoubian then $F\in\Tto (\cA ,\cK )$ has a right adjoint if
and {\it only\,} if $F\in T(\cA ,\cK )$.

(ii) At the 2-category level the definitions of the right
derived DG functor from \ref{RFnew} and \ref{RF} amount to
the following one. Suppose that in the situation of
\ref{univ} we are given $F\in T(\cA ,\cA')$. Then
$RF\in\Tto (\cC ,\cA')$ is the composition of 
$F\in T(\cA ,\cA')\subset\Tto (\cC ,\cA')$ and the right
adjoint $\ppi^*\in\Tto (\cC ,\cA )$ of 
$\ppi\in T (\cA ,\cC )$.

\subsection{Relation with Kontsevich's approach}
\subsubsection{}
\label{Konts1}
Let $\cA ,\cK$ be DG categories and suppose that $\cA$ is
flat. Given a DG functor $F:\cA\to\cK$ denote by $\Phi_F$ the
DG $\cA\otimes\cK^{\circ}$-module 
$(X,Y)\mapsto\Hom (Y,F(X))$. Clearly 
$\Phi_F\in D(\cA^{\circ}\otimes\cK )$ belongs to 
$T(\cA ,\cK )$. Let us describe the full subcategory of
$T(\cA ,\cK )$ formed by the DG
$\cA\otimes\cK^{\circ}$-modules $\Phi_F$. One has
$\Phi_F=\Ind_{\id_{\cA}\otimes F^{\circ}} (\HHom_{\cA})$,
where $F^{\circ}$ is the DG functor
$\cA^{\circ}\to\cK^{\circ}$ corresponding to $F:\cA\to\cK$
and $\HHom_{\cA}$ is the $\cA^\circ\otimes\cA$-module
$(X,Y)\mapsto\Hom (X,Y)$. As $\cA$ is homotopically flat over
$k$ the morphism
$L\Ind_{\id_{\cA}\otimes F^{\circ}} (\HHom_{\cA})\to
\Ind_{\id_{\cA}\otimes F^{\circ}} (\HHom_{\cA})$ is a
quasi-isomorphism. Therefore the adjunction between derived
induction and restriction yields a canonical isomorphism
\begin{equation}  \label{Extisoext}
\Ext^n(\Phi_F,\Phi_G)= \Ext^n(L\Ind_{\id_{\cA}\otimes
F^{\circ}} (\HHom_{\cA}),\Phi_G)\iso\Ext^n (F,G), 
\end{equation}
where 
$\Ext^n (F,G):=\Ext^n _{\cA\otimes\cA^{\circ}} (\HHom_{\cA}
\, ,\, \mathcal Hom (F,G))$ and 
$\mathcal Hom (F,G):=$ $\Res_{\id_{\cA}\otimes
F^{\circ}}(\Phi_G)$, \, i.e., \,
$\mathcal Hom (F,G)$ \, is \, the \, DG
$\cA\otimes\cA^{\circ}$-module \,
$(X,Y)\mapsto$ $\Hom (F(Y),G(X)),\,$ $X,Y\in\cA$. The
morphism 
$\Ext^m (F_2,F_3)\otimes\Ext^n (F_1,$ $F_2)\to
\Ext^{m+n}(F_1,F_3)$ coming from (\ref{Extisoext}) is, in
fact, induced by the morphism 
$\mathcal Hom (F_2,F_3)\otimes\mathcal Hom (F_1,F_2)
\to\mathcal Hom (F_1,F_3)$ and the
quasi-isomomor\-phism 
$(\HHom_{\cA})\otimes_{\cA}(\HHom_{\cA})\to\HHom_{\cA}$. So
we have described the full subcategory of $T(\cA ,\cK )$
formed by the DG $\cA\otimes\cK^{\circ}$-modules $\Phi_F$.
The next statement shows that it essentially equals 
$T(\cA ,\cK )$ if $\cA$ is semi-free.

\subsubsection{\bf Proposition}  \label{Konts2}
{\it If $\cA$ is semi-free over $k$ then every object of
$T(\cA ,\cK )$ is isomorphic to $\Phi_F$ for
some $F:\cA\to\cK$.}

\begin{proof}
An object $\Phi\in T(\cA ,\cK )$ is a
DG $\cA\otimes\cK^{\circ}$-module. Consider $\Phi$ as a DG
functor $\cA\to\cK'\subset\cK^{\circ}\DGmod$, where
$\cK'$ is the full DG subcategory of
quasi-representable DG modules. We have the DG functors
$\cK\leftarrow\cK''\buildrel{\pi}\over{\longrightarrow}\cK'$,
where $\cK''$ is the DG category whose objects are triples
consisting of an object $Y\in\cK$, a DG
$\cA\otimes\cK^{\circ}$-module $\Psi$, and a
quasi-isomorphism $f:h_Y\to\Psi$ (see \ref{fixquasirep} for a
precise definition of $\cK''$). We also have a canonical
DG functor $\Cone :\cK''\to\cK^{\circ}\DGmod$, which sends
$(Y,\Psi ,f)$ to $\Cone (f)$ (the
definition of the $\Cone$ functor on morphisms is clear from
\ref{Mmor}). $\cA$ is semi-free and
$\pi$ is a surjective quasi-equivalence, so by
\ref{closedmodel} our DG functor $\cA\to\cK'$ lifts to a DG
functor
$\cA\to\cK''$. Let $F:\cA\to\cK$ be the composition
$\cA\to\cK''\to\cK$. One has an exact sequence of DG
$(\cA\otimes\cK^{\circ})$-modules $0\to\Phi\to
M\to\Phi_F[1]\to 0$, where $M$ corresponds to the 
composition
$\cA\to\cK'\buildrel{\Cone}
\over{\longrightarrow}\cK^{\circ}\DGmod$.
As
$M$ is acyclic we get a $T(\cA ,\cK )$-isomorphism
$\Phi_F\iso\Phi$.
\end{proof}

\subsubsection{\bf The standard resolution}  \label{stand}
Consider the category $\DGalg$ of (non-unital) associative DG
algebras and the category $\DGcoalg$ of
(non-counital) cocomplete coassociative coalgebras
( a coalgebra $U$ is {\it cocomplete} if for every $u\in U$ 
there exists $n\in\BN$ such that $u$ is annihilated by the
$n$-fold coproduct $\Delta_n:U\to U^{\otimes n}$). If
$U\in\DGcoalg$ and $A\in\DGalg$ then $\Hom (U,A)\in\DGalg$
(the product of $f:U\to A$ and $g:U\to A$ is defined to be
the composition of the coproduct $U\to U\otimes U$, the map 
$f\otimes g:U\otimes U\to A\otimes A$, and the product
$m:A\otimes A\to A$). Define the {\it Maurer--Cartan
functor\,} $\MC:\DGcoalg^{\circ}\times\DGalg\to$Sets as
follows: $\MC (U,A)$ is the set of elements 
$\omega\in\Hom (U,A)$ of degree 1 such that
$d\omega+\omega^2=0$. There exist functors 
$B:\DGalg\to\DGcoalg$ and $\Omega:\DGcoalg\to\DGalg$ such
that $\MC (U,A)=\Mor (U,BA)=\Mor (\Omega U,A)$ (they are
called ``bar construction'' and ``cobar construction''). As
$\Omega$ is left adjoint to $B$ we have the adjunction
morphisms $\Omega BA\to A$ and $U\to B\Omega U$. In fact,
they are quasi-isomorphisms. The above statements are
classical (references will be given in \ref{hist}).

Caution: while $B$ sends quasi-isomorphisms to
quasi-isomorphisms this is {\it not\,} true for $\Omega$.
Indeed, consider the morphism $\varphi :0\to k$, where $k$ is
equipped with the obvious DG algebra structure. Then $B
(\varphi )$ is a quasi-isomorphism but $\Omega B(\varphi )$
is not.

It is easy to see that if $A$ is a
semi-free DG $k$-module then $\Omega BA$ is a semi-free DG
algebra (in the non-unital sense), so $\Omega BA$ is a
semi-free resolution of $A$.
$\Omega BA$ is non-unital even if $A$ is unital. The DG
algebra one gets by adding the unit to a DG algebra $B$ will
be denoted by $u(B)$. If $A$ is unital then $u(A)$ is
the Cartesian product of DG algebras $A$ and $k$, so we get
a quasi-isomorphism $u(\Omega BA )\to u(A)=A\times k$. Let us
call it the {\it standard resolution\,} of $A\times k$. It is
semi-free (in the unital sense) if $A$ is a semi-free DG
$k$-module.

As explained in \cite{K2,Ke3,KS,L-H}, there is a similar
construction in the more general setting of DG categories.
Given a DG category $\cA$ let $\cA_{\discr}$ denote the DG
category with $\Ob\cA_{\discr}=\Ob\cA$ such that the
endomorphism DG algebra of each object of $\cA_{\discr}$
equals $k$ and $\Hom_{\cA_{\discr}}(X,Y)=0$ if $X,Y$ are
different objects of $\cA_{\discr}$. Let 
$u (\cA )\subset\cA\times\cA_{\discr}$ be the full DG
subcategory formed by objects $(a,a)$,
$a\in\Ob\cA =\Ob\cA_{\discr}$. There is a {\it standard 
resolution\,} $\Stand (\cA )\to u(\cA )$. If
all Hom complexes of $\cA$ are semi-free over $k$ then 
$\Stand (\cA )$ is semi-free.

\subsubsection{\bf $A_{\infty}$-functors} \label{Konts3}
If $\cA$ is any DG category and $\tilde\cA$ is a semi-free
resolution of $\cA$ then $T(\cA ,\cK )=T(\tilde\cA ,\cK )$,
so \ref{Konts1}-\ref{Konts2} give a graded $k$-category
equivalent to $T(\cA ,\cK )$ whose objects are DG functors
$\tilde\cA\to\cK$. In particular, if all Hom complexes of
$\cA$ are semi-free (or, more generally, homotopically
projective) over $k$ we get a  category equivalent to
$T(u(\cA ),\cK )$ whose objects are DG functors $\Stand (\cA
)\to\cK$. Notice that if $k$ is a field (and if you believe
in the axiom  of choice, which ensures that modules over a
field are free) then every DG
$k$-module is semi-free. The functor 
$T(\cA ,\cK )\to T(u(\cA ),\cK )$ corresponding to the
canonical projection $u(\cA )\to\cA$ is fully faithful (this
follows from the decomposition 
$D(u(\cA )^{\circ}\otimes\cK )=D(\cA^{\circ}\otimes\cK )
\oplus D(\cA_{\discr}^{\circ}\otimes\cK )\,$). 
DG functors $\Stand (\cA )\to\cK$ such that the corresponding
object of $T(u(\cA ),\cK )$ is in 
$T(\cA ,\cK )\subset T(u(\cA ),\cK )$ are
called {\it $A_{\infty}$-functors.} More precisely, this is
{\it one\,} of the versions of the notion of
$A_{\infty}$-functor $\cA\to\cK$. They differ on how an
$A_{\infty}$ analog of the axiom $F(\id )=\id$ in the
definition of usual functor is formulated (the difference is
inessential from the homotopy viewpoint). The above notion
is as ``weak'' as possible. 

According to Kontsevich, the structure of graded
$k$-category on $T(\cA ,\cK )$ comes from a canonical
DG category $A_{\infty}\mbox{\rm -funct} (\cA,\cK )$ whose
objects are $A_{\infty}$-functors $\cA\to\cK$. Here is its
definition if $\cA$ and $\cK$ have one
object (the general case is similar). Let $A,K$ be the
endomorphism DG algebras of these objects. Then an
$A_{\infty}$-functor $\cA\to\cK$ is a DG algebra
morphism $\Omega BA\to K$ satisfying a certain condition
(see \ref{stand}). So it remains to construct a DG category
whose objects are elements of $\Mor (\Omega BA,K)=\MC
(BA,K)$, i.e., elements $\omega$ of the DG algebra
$R:=\Hom (BA,K)$ such that $\deg\omega =1$ and 
$d\omega +\omega^2=0$. Such $\omega$ defines a DG
$R^{\circ}$-module $N_{\omega}$: it equals $R$ as a
graded $R^{\circ}$-module, and the differential
in $N_{\omega}$ maps $r$ to $\nabla r:=dr+\omega r$. Now put
$\Hom (\omega ,\omega'):=\Hom (N_{\omega},N_{\omega'})$ and
define the composition map $\Hom (\omega ,\omega')\times
\Hom (\omega' ,\omega'')\to\Hom (\omega ,\omega'')$ in the
obvious way.

\medskip

\noindent {\bf Remark.} According to \cite{KS,L-H}, in the
more general case that $\cK$ is an $A_{\infty}$-category
$A_{\infty}$-functors $\cA\to\cK$ form an
$A_{\infty}$-category. Kontsevich informed me that if 
$\cK$ is a DG category then the $A_{\infty}$-category of
$A_{\infty}$-functors $\cA\to\cK$ is a DG category. I do not
know if this DG category equals the above DG category
$A_{\infty}\mbox{\rm -funct} (\cA,\cK )$. 

\subsection{DG models of $T(\cA_1,\cA_2)$} \label{DGmodels}
{\it Kontsevich's model\,} has already been mentioned in
\ref{Konts3}:  if the Hom complexes of $\cA_1$ are semi-free
(or, more generally, homotopically projective) over $k$ then
$T(\cA_1,\cA_2)$ is the graded homotopy category of the DG
category $A_{\infty}\mbox{\rm -funct} (\cA_1,\cA_2)$.

Keller's model is easier to define. If $\cA_1$ or $\cA_2$ is
flat then
$D(\cA_1^{\circ}\Ltensor\cA_2)=D(\cA_1^{\circ}\otimes\cA_2)=
\Hodot (\cRto )$, where $\cR:=\cA_1^{\circ}\otimes\cA_2$ and
$\cRto$ is the DG category of semi-free
DG $\cR^{\circ}$-modules. This identifies
$T(\cA_1,\cA_2)\subset D(\cA_1^{\circ}\otimes\cA_2)$ with the
graded homotopy category of a certain full DG subcategory
$\dg(\cA_1,\cA_2)\subset\cRto$, which will be called
{\it Keller's model}.

One also has the {\it dual Keller model} $(\dg
(\cA_1^{\circ},\cA_2^{\circ}))^{\circ}$: its graded homotopy
category is 
$T(\cA_1^{\circ},\cA_2^{\circ})^{\circ}=T(\cA_1,\cA_2)$.
The equality
$T(\cA_1,\cA_2)=T(\cA_1^{\circ},\cA_2^{\circ})^{\circ}$
identifes $T(\cA_1,\cA_2)$ with the graded homotopy category
of the DG category
$(\dg (\cA_1^{\circ},\cA_2^{\circ}))^{\circ}$, which is  a
full DG subcategory of the DG category $\cRfrom:=$\{the dual
of the DG category of semi-free DG $\cR$-modules\}.

If the Hom complexes of $\cA_1$ are homotopically projective
over $k$ there is a canonical quasi-equivalence 
$A_{\infty}\mbox{\rm -funct}
(\cA_1,\cA_2)\to\dg(\cA_1,\cA_2)$, which is not
discussed here.

\medskip

\noindent {\bf Remark.}  Let $\cA ,\cC_1 ,\cC_2$ be DG
categories and suppose that $\cC_1 ,\cC_2$ are flat.
Then $\dg (\cA ,\cC_1)$,  $\dg (\cC_1 ,\cC_2)$, and
$\dg (\cA ,\cC_2)$ are defined, but in general (if $\cC_1$
is not semi-free) the image of 
\begin{displaymath} 
\otimes_{\cC_1}:\dg (\cA ,\cC_1)\otimes \dg (\cC_1
,\cC_2)\to (\cA\otimes\cC_2^{\circ})\DGmod
\end{displaymath}
is not contained in $\dg (\cA ,\cC_2)$ or even in
$\cRto$, where $\cR:=\cA^{\circ}\otimes\cC_2$. So we do not
get a composition DG functor 
$\dg (\cA ,\cC_1)\otimes \dg (\cC_1,\cC_2)\to\dg (\cA
,\cC_2)$
but rather a DG functor
\begin{equation} \label{quasitensor2}
\Psi :\dg (\cA ,\cC_1)\times \dg (\cC_1 ,\cC_2)\times 
\dg (\cA ,\cC_2)^{\circ}\to k\dgmod ,
\end{equation}
which lifts the graded functor
\begin{displaymath} 
T(\cA ,\cC_1)\times T(\cC_1 ,\cC_2)\times 
T(\cA ,\cC_2)^{\circ}\to\mbox{\{Graded }
k\mbox{-modules\}}
\end{displaymath}
defined by $(F_1,G,F_2)\mapsto\bigoplus_n\Ext^n(F_2,GF_1)$.
One defines (\ref{quasitensor2}) by
\begin{displaymath} 
(M_1,N,M_2)\mapsto\Hom(M_2,M_1\otimes_{\cC_1}N).
\end{displaymath}

\subsection{Some historical remarks}  \label{hist}
As explained in \cite{M}, the functors $B$ and $\Omega$
from \ref{stand} go back to Eilenberg -- MacLane and
J.~F.~Adams. It was E.~H.~Brown \cite{Br0} who introduced
$MC(U,A)$; he called its elements ``twisting cochains''. The 
fact that the morphism $\Omega BA\to A$ is a
quasi-isomorphism appears as Theorem 6.2 on p.~7-28 of
\cite{M1}. All the properties of $B$ and $\Omega$ from
\ref{stand} were formulated in \cite{M} and proved in
\cite{HMS}; their analogs for Lie algebras and commutative
coalgebras were proved in \S7 of Appendix B of \cite{Q}. In
these works DG algebras and DG coalgebras were assumed to
satisfy certain boundedness conditions. The general case was
treated in \cite{H',L-H}.

\end{document}